

\input harvmac
\input epsf

\catcode`@=12

\baselineskip16pt
\noblackbox

\font\tabfont=cmr7 scaled\magstep0

\def\undertext#1{$\underline{\smash{\hbox{#1}}}$}

\def\A{{\cal A}}
\def\B{{\cal B}}
\def\IP{{\bf P}}
\def\IR{{\bf R}}
\def\IC{{\bf C}}
\def\IQ{{\bf Q}}
\def\ZZ{{\bf Z}}
\def\pd{\partial}

\def\D{\Delta}
\def\Dp{{\Delta'}}
\def\Ds{{\Delta^*}}
\def\Dps{{\Delta'^*}}
\def\({ \left(  }
\def\){ \right) }
\def\S{\Sigma}
\def\SD{{\Sigma(\D)}}
\def\SDs{{\Sigma(\Ds)}}
\def\SDps{{\Sigma(\Dps)}}
\def\bn#1{\bar\nu_{#1}^*}
\def\ns#1{\nu_{#1}^*}
\def\ta#1{\theta_{a_#1}}
\def\da#1{{\pd \; \over \pd {a_{#1}}}}
\def\g#1{{\gamma_{#1}}}
\def\en#1{{e_{\bn{#1}}}}
\def\tx{\theta_x}
\def\ty{\theta_y}
\def\tz{\theta_z}
\def\tw{\theta_w}
\def\sqbox{{\cal D}}
\def\drho#1{{\pd_{\rho_{#1}}}}
\def\X#1#2(#3)#4#5{ {$X_{#1}^#2(#3)_{#5}^{#4}$} }

\def\I{{\rm I}}
\def\II{{{\rm I}\hskip-1.5pt{\rm I}} }
\def\III{{{\rm I}\hskip-1.5pt{\rm I}\hskip-1.5pt{\rm I}} }

\Title{alg-geom/9511001}
{\vbox{
 \centerline{GKZ-Generalized Hypergeometric Systems}
 \vskip2pt
 \centerline{in Mirror Symmetry of Calabi-Yau Hypersurfaces} } }

\centerline{S. Hosono$^1$\footnote{$^\dagger$}
{email: hosono@sci.toyama-u.ac.jp},
B.H. Lian$^2$\footnote{$^\ddagger$}{email: lian@max.math.brandeis.edu}
and S.-T. Yau$^3$\footnote{$^\diamond$}{email: yau@math.harvard.edu}}

\bigskip\centerline{
\vbox{
\hbox{ $^1$ Department of Mathematics}
\hbox{ \hskip30pt Toyama University}
\hbox{ \hskip28pt Toyama 930, Japan} }
\hskip0.5cm
\vbox{
\hbox{ $^2$ Department of Mathematics}
\hbox{ \hskip30pt Brandeis University}
\hbox{ \hskip27pt Waltham, MA 02154} } }
\bigskip
\centerline{
\vbox{
\hbox{ $^3$ Department of Mathematics}
\hbox{ \hskip30pt Harvard University}
\hbox{ \hskip23pt Cambridge, MA 02138} } }

\bigskip 

We present a detailed study of the generalized hypergeometric system
introduced by Gel'fand, Kapranov and Zelevinski (GKZ-hypergeometric system)
in the context of toric geometry. GKZ systems arise naturally in
the moduli theory of Calabi-Yau toric varieties, and play
an important role in applications of the mirror
symmetry. We find that the Gr\"obner basis for the so-called
 toric ideal determines a finite set of differential operators
for the local solutions of the GKZ system. At the
 special point called the large
radius limit, we find a close relationship between the principal parts of the
operators in the GKZ system and the intersection ring of a toric
variety.  As applications, we analyze general three dimensional
 hypersurfaces of Fermat
and non-Fermat
types with Hodge numbers up to
$h^{1,1}=3$. We also find and analyze several non Landau-Ginzburg
models which are related to singular models.
\Date{Oct.1995}

\lref\HKTYI{S.Hosono, A.Klemm, S.Theisen and S.-T.Yau,
            Commun.Math.Phys.{\bf 167}(1995)301.}
\lref\HKTYII{S.Hosono, A.Klemm, S.Theisen and S.-T.Yau,
             Nucl.Phys.{\bf B433}(1995)501.}
\lref\BatyrevI{V.Batyrev, J.Algebraic Geometry {\bf 3}(1994)493.}
\lref\BatyrevII{ V.Batyrev, Duke Math.J.{\bf 69}(1993)349.}
\lref\GKZ{ I.M.Gel'fand, A.V.Zelevinski and M.M.Kapranov, Funct. Anal. Appl.
            {\bf 23}, 2. }
\lref\HLY{ S.Hosono,B.H.Lian and S.-T.Yau,
           optional appendix to alg-geom/9511001 }
\lref\CandelasETAL{P.Candelas, X. de la Ossa, P.Green and L.Parks,
                  Nucl.Phys.{\bf B359}(1991)21.}
\lref\CandelasIIa{P.Candelas, X. de la Ossa, A.Font, S.Katz and D.Morrison,
                  Nucl.Phys.{\bf B416}(1994)481.}
\lref\CandelasIIb{P.Candelas, X. de la Ossa, A.Font, S.Katz and D.Morrison,
                  Nucl.Phys.{\bf B429}(1994)629.}
\lref\Morrison{D.Morrison,
      {\it Picard-Fuchs Equations and Mirror Maps For Hypersurfaces},
      in {\it Essays on Mirror Manifolds},
      Ed. S.-T.Yau (1992) International Press.}
\lref\KTa{A.Klemm and S.Theisen, Nucl.Phys.{\bf B389}(1993)153.}
\lref\KTb{A.Klemm and S.Theisen, Theor.Math.Phys.{\bf 95}(1993)583.}
\lref\Font{A.Font, Nucl.Phys.{\bf B391}(1993)358.}
\lref\KM{M.Kontsevich and Yu.Manin, Commun.Math.Phys.{\bf 164}(1994)525.}
\lref\Kont{M.Kontsevich,
      {\it Enumeration of rational curves via torus action},
      hep-th/9405035.}
\lref\LYa{ B.H.Lian and S.-T.Yau,
      {\it Arithmetic properties of the mirror map and quantum coupling},
      hep-th/9411234.}
\lref\LYb{ B.H.Lian and S.-T.Yau,
      {\it Mirror Maps, Modular Relations and Hypergeometric Series I,II},
      hep-th/9507151, 9507153.}
\lref\KLM{A.Klemm, W.Lerche and P.Myer,
       {\it K3-Fibrations and Heterotic-Type II String Duality},
      hep-th/9506112.}
\lref\KachruVafa{ S.Kachru and C.Vafa,
       {\it Exact Results for N=2 Compactifications of Heterotic Strings},
      hep-th/9505105.}
\lref\LVW{ W.Lerche, C.Vafa and N.Warner, Nucl.Phys.{\bf B329}(1990)163.}
\lref\DistlerGreene{ J.Distler and B.Greene, Nucl.Phys.{\bf B309}(1988)295.}
\lref\Roan{S.-S.Roan, Int.J.Math.{\bf 2}(1991)439.}
\lref\KKLMV{S. Kachru, A. Klemm, P. Mayr, W. Lerche and C. Vafa,
      {\it Nonperturbative Results on the Point Particle Limit of N=2
       Heterotic String Compactification}, hep-th/9508155.}
\lref\SkarkeETAL{M.Kreuzer and H.Skarke, Nucl.Phys.{\bf B388}(1993)113.}
\lref\lt{A. Libgober and J. Teitelboim, {\sl Duke Math. Journ., Int.
         Res. Notices} 1 (1993) 29.}
\lref\AGM{P.Aspinwall, B.Greene and D.Morrison,
          Phys.Lett.{\bf B303}(1993)249.}
\lref\BKK{P.Berglund, S.Katz and A.Klemm, {\it Mirror Symmetry and
        the Moduli Space for Generic Hypersurfaces in Toric Varieties},
        hep-th/9506091.}
\lref\yonemura{ T.Yomemura, T\^ohoku Math. J. {\bf 42}(1990),351.}
\lref\BH{P.Berglund and T.H\"ubsch, Nucl.Phys.{\bf B393}(1993)377.}
\lref\BatyrevP{Batyrev,
 {\it Quantum Cohomology Ring of Toric Manifolds}, alg-geom/9310004.}
\lref\AspinwallMorr{ Aspinwall and D. Morrison, Commun.Math.Phys.{\bf 151}245.}
\lref\topWitten{E.Witten, {\it Mirror Manifolds and Topological Field Theory}
      in {\it ``Essays on Mirror Symmetry''}, ed. S.-T.Yau, (1992)
      International Press.}
\lref\CandelasKatz{P.Candelas,X.de la Ossa and S.Katz, Nucl.Phys.{\bf B450}
                   (1995)267.}
\lref\GH{P.Griffiths and J.Harris, {\it Principles of Algebraic Geometry},
         Wiley-Interscience(1978)}

\newsec{Introduction}

Recent studies on nonperturbative aspects of string
theory have made remarkable progress in understanding the structure of
moduli spaces in string theory.
Applications of mirror symmetry, for example, in type II string
compactification to studying the geometry of moduli spaces is one of
the most successful developments. Starting from the pioneering work by
Candelas {\it et al}\CandelasETAL, and subsequently by others,
 the quantum geometry of the moduli spaces for many Calabi-Yau
models\Morrison\KTa\lt\Font\KTb\CandelasIIa\HKTYI\CandelasIIb\HKTYII\BKK\
have now been well understood via mirror symmetry.
At the same time, there is parallel progress in studying the axiomatic
framework of quantum geometry and its application to
enumerative geometry\KM.
Also in explicit constructions of the geometry of concrete Calabi-Yau
models, it is now understood that for a large class of Calabi-Yau
varieties, the mirror maps have remarkable modular and integrality properties
\LYa\KLM\LYb.
These models present strong and even beautiful evidence for
the recent proposal for the so-called type II-heterotic string
duality\KachruVafa.
These Calabi-Yau models continue to provide fruitful testing ground for
string duality \KKLMV.

Mirror symmetry was first recognized in the local operator algebra of
the N=2 string theory\LVW. Soon after the introduction of the
framework of toric geometry into the study of Calabi-Yau models
\BatyrevI\Roan, mirror symmetry has since been widely checked
for many Calabi-Yau hypersurfaces and complete intersections
in toric varieties.  Mirror
symmetry relates two moduli spaces with apparently very
different properties -- one moduli space is described by purely
classical geometry, while the other is described by
quantum geometry which receives nonperturbative corrections from worldsheet
instanton\DistlerGreene. Mirror symmetry thus gives us a powerful
means for studying quantum geometry of one moduli space via
classical means such as
the theory of variation of
Hodge structures.

Variation of Hodge structures allows us to
study the period integrals for Calabi-Yau varieties.
It is known that the
period integrals satisfy differential equations with regular
singularities, known as Picard-Fuchs differential equations. A
general technique  for constructing Picard-Fuchs equations is the reduction
method of Dwork-Griffith-Katz. For Calabi-Yau
toric varieties, it was remarked in \BatyrevII\ that the period integrals
satisfy a generalized hypergeometric system introduced by
Gel'fand-Kapranov-Zelevinski\GKZ.
It has been observed\HKTYI\ in solving several examples that
the GKZ system is not generic and is reducible. Moreover
there is an irreducible part in which
the period integrals live. In this paper we study the GKZ hypergeometric
system for general Calabi-Yau hypersurfaces, and discuss the
previous observations in a different light but with
much greater generality. As  applications, we
determine the Picard-Fuchs differential equations for all hypersurfaces
with Hodge numbers $h^{1,1}\leq 3$
in weighted projected spaces.

In section 1, we review the toric description of mirror
symmetry, due to Batyrev.
We introduce period integrals in the language of toric geometry, and
introduce a GKZ system which we call $\Ds$-hypergeometric system.
The system is extended by incorporating the symmetry
coming from the automorphism group of the ambient space\HKTYI.
We classify according to the toric data \HKTYI\ Calabi-Yau hypersurfaces
into three classes: types I, II and III.

In section 2, we analyze local solutions to the
$\Ds$-hypergeometric system. We  construct a finite set of
differential operators for local solutions by relating the system to an
algebro-combinatorial object, known as a toric ideal. We find that the local
properties near the so-called large radius limit are determined
completely by the
intersection ring of the
ambient space. In the case of type I and type II models,
we prove in general the existence of the large radius limit, hence establish
the existence of the point of maximally unipotent monodromy.
We give a natural explanation for the reducibility of our
$\Ds$-hypergeometric system in terms of certain aspects of
the intersection ring of the ambient space.  We also extend
our arguments to type III models.

In section 3, we will apply our general framework to three dimensional
Calabi-Yau
hypersurfaces with $h^{1,1}\leq 3$. Detailed analyses are given
for a few typical models. For others, we will append a list of the
Picard-Fuchs equations to the source file of this article\HLY\ for
interested readers.

In the final section we will discuss some relationships among different
Calabi-Yau manifolds which come from the inclusion relations among
reflexive polyhedra.

\vfill\eject

\newsec{Toric Geometry and Generalized Hypergeometric Differential Equation}

In this section we analyze the differential equations,
known as Picard-Fuchs equations, satisfied by the periods of a toric variety.
Applications of toric geometry to the description
of the Picard-Fuchs equation was first initiated in \BatyrevII\ and
further developed in \HKTYI.
Here we summarize some of the analyses in \HKTYI\ and extract some
combinatorial aspects of the Picard-Fuchs equations.

\subsec{A construction of mirror manifolds}

In order to fix some notations, we review Batyrev's construction of
the mirror manifolds, which is applicable to the list of 7,555
hypersurfaces of
\ref\KlemmSchimmrigk{A.Klemm and R.Schimmrigk,Nucl.Phys.{\bf B411}(1994)559.}
\SkarkeETAL\ as well as  complete intersections
\ref\cicylist{
P.Candelas
, A.M.Dale, C.A.L\"utken and R.Schimmrigk,
Nucl.Phys.{\bf B298}(1988)493\semi
P.Candelas, C.A.L\"utken and R.Schimmrigk, Nucl.Phys.{\bf B306}(1988)113.}
in a product of (weighted) projective spaces.
In the following we restrict our attention to
hypersurfaces, although generalization to complete intersections
\ref\Hdiff{V.Batyrev and D.van Straten,{\it Generalized Hypergeometric
Functions and Rational Curves on Calabi-Yau Complete Intersections in
Toric Varieties}, Alg-geom/9307010}
\ref\BatyrevBorisov{
V.Batyrev and L.Borisov,{\it On Calabi-Yau Complete Intersections in
Toric Varieties},alg-geom/9412017}\ can be done.

Let us consider a weighted projective space $\IP^n(w)$ and a
hypersurface $X_d(w)$
with (weighted) homogeneous degree $d=w_1+\cdots+w_{n+1}$.
Without loss of generality, we
may assume that the weight $w$
is normalized \ref\Dais{D.I.Dais,
{\it Enumerative Combinatorics of Invariants of Certain Complex Threefolds
with Trivial Canonical Bundle},MPI-preprint(1994)}, i.e.,
$gcd(w_1,\cdots,\hat w_i,\cdots,w_{n+1})=1 \,, (i=1,\cdots,n+1)$.
(See also \ref\weight{I.Dolgachev,
In {\it Group actions and vector fields}, Lecture Notes in Math.{\
bf Vol.956}(1991),Springer-Verlag,34.}.)
For $n=4$, the list of \KlemmSchimmrigk\SkarkeETAL\ exhausts all hypersurfaces
$X_d(w)$ defined by weighted homogeneous polymonials satisfying
the transversality condition. Now let
\eqn\pot{
W(z)=\sum_{(w,m)=d} a_m z^m =
     \sum_{(w,m)=d} a_{m_1,\cdots, m_{n+1}}
      z_1^{m_1} \cdots z_{n+1}^{m_{n+1}} \;.   }
For generic $a_{m}$, the zero locus $\{W(z)=0\}$ defines a
 hypersurface $X_d(w)$
in general position. Its intersection with
singular locus of the ambient space $\IP(w)$ gives the singular locus
of the hypersurface.
We denote the Newton polyhedron of $W(z)$ as $\D(w)$. It is the
convex hull of the exponents of \pot\ $m$ in $\IR^{n+1}$, shifted by
$(-1,\cdots,-1)$.  If we take into account the condition
$d=w_1+\cdots+w_{n+1}$, it is easy to deduce that the shifted polyhedron
can be written as
\eqn\Dw{
\D(w)={\rm Conv.}\(
   \{\;  x \in \ZZ^{n+1} \;
     \vert \; (w,x)=0 \;,\; x_i \geq -1 \; (i=1,\cdots,n+1) \; \}
           \) \;.  }

An $n+1$-dimensional polyhedron $\D$ in $\IR^{n+1}$ is called
integal if all its vertices are integral (with
respect to the lattice $\ZZ^{n+1}$). A {\it reflexive}
polyhedron is an integral polyhedron with exactly one integral interior
point, the origin.
The polar dual of $\D$,
\eqn\polard{
\D^* := \{ \; y \in \IR^{n+1} \;\vert\;
          (y,x)\geq -1 \;\; (\forall x \in \D) \;\}   }
is again integral and reflexive.
 If we consider the set of cones over the faces of a
polyhedron, we will obtain a complete fan which covers
$\IR^{n+1}$.
Thus to each pair of reflexive polyhedra $(\D,\Ds)$, we can
associate a pair of complete fans $( \SD, \SDs )$ and in turn
a pair of the $n+1$ dimensional toric varieties
$(\IP_{\SDs},\IP_{\SD})$.
In each of the
toric varieties, there is a family of Calabi-Yau hypersurfaces
given by the zero loci of certain sections of the anticanonical bundle.
The toric variety $\IP_{\SD}$ contains a canonical Zariski open torus
$(\IC^*)^{n}$ whose coordinates we denote as
$X=(X_1,\cdots,X_{n})$. In these coordinates, the sections are
\eqn\defeq{
f_{\D^*}(X,a)=\sum_{\nu^*_i \in \D^* \cap \ZZ^{n+1}} a_i X^{\nu^*_i}  \;\;.
}
For generic values of the $a_i$'s in \defeq, the
$X_{\Ds}$ in $\IP_{\SD}$ admits a minimal resolution to a Calabi-Yau manifold
(which we also denote $X_{\Ds}$).
 Similarly there is a corresponding family of
hypersurfaces $X_\D$ in $\IP_\SDs$. Batyrev showed that a
pair of the Calabi-Yau manifolds $(X_{\D},X_{\Ds})$ is mirror
symmetric to each other in the sense that we have the following
relations for their Hodge numbers $(n\geq4)$;
\eqn\Hodge{
\eqalign{
h^{1,1}(X_\D)
&= h^{n-2,1}(X_{\D^*})   \cr
&=l(\D^*)-(n+1)-\sum_{{\rm codim} S^* =1} l'(S^*)
  +\sum_{{\rm codim} S^*=2} l'(S^*)l'(S)   \;\;, \cr
h^{n-2,1}(X_\D)
&= h^{1,1}(X_{\D^*})   \cr
&=l(\D)-(n+1)-\sum_{{\rm codim} S =1} l'(S)
  +\sum_{{\rm codim} S=2} l'(S)l'(S^*)   \;\;, \cr  } }
where the $S$ are faces of $\D$, $S^*$ the polar dual face of $S$.
The functions $l$ and $l'$ count the numbers of integral
points in a face and in the interior of a face respectively.

When $W(z)$ is Fermat, the toric
variety $\IP_{\SDs}$ is isomorphic to the weighted projective
space $\IP^n(w)$, with $X_{\D}$ isomorphic to some
$X_d(w)$. Then the mirror hypersurface $X_{\Ds}$
can be understood\BatyrevI\
as an orbifold of the $X_{\D}$ in
$\IP_{\S(\D^*)}$, giving the orbifold construction of Greene and
Plesser\ref\GreenePl{B.Greene and M.Plesser, Nucl.Phys.{\bf B338}(1990)15.}
based on conformal field theory. For general hypersurfaces of non-Fermat
type, $\IP(w)$ and $\IP_{\SDs}$ are only birational.
In fact the fan $\S(\D^*)$ is a refinement of the fan of $\IP(w)$.
The hypersurfaces $X_d(w)$ and $X_{\D}$
are related by flop operations on the ambient
spaces. It has been shown\CandelasKatz ,
in this way, that Batyrev's constructions
applies to all 7,555 hypersurfaces and reproduces the generalized mirror
constructions known to\BH. In addition, there are
several mirror pairs $(X_\D,X_{\D^*})$ which do not come from
hypersurfaces in weighted projective spaces.

\vskip0.5cm

The quantity most relevant to the applications of the mirror symmetry
to the quantum geometry of $X_\D$ are the period integrals
for its mirror $X_\Ds$. For example,
\eqn\period{
\Pi(a)={1\over (2\pi i)^n}
\int_{C_0} {1\over f_\Ds (X,a)} \prod_{i=1}^n {d X_i \over X_i}
\quad ,
}
is the period integral over the torus cycle $C_0=\{ |X_1|=|X_2|= \cdots
=|X_n|=1\}$ in $(\IC^*)^n$. For other periods, we will analyze the
differential equation satisfied by \period.

\subsec{$\A$-hypergeometric system for the periods}

In \Hdiff,
it is remarked that the period integral \period\ satisfies an
$\A$-hypergeometric system introduced by Gel'fand, Kapranov and Zelevinski
\GKZ. In \HKTYI, it is found that the hypergeometric system is
not generic but reducible, and the period integrals can be extracted
from the system as the irreducible part of its solution space.
Furthermore, for most of the hypersurface models, it is noted that
the hypergeometric system must be generalized in order to extract
the irreducible part of the solutions. We reproduce here an extension
which is called an extended $\Ds$-hypergeometric system, from purely
combinatorial data of the polyhedron. We note that
for type I models (see below),
the extended $\Ds$-hypergeometric system coincides with the GKZ system. For
type
II or III, the extended $\Ds$-hypergeometric system incorporates additional
differential operators
associated with the action of an automorphism group.

An $\A$-hypergeometric system is described by a finite set
$\A$ in a lattice $\{1\}\times \ZZ^{n}$
with the property that $\A$ linearly spans $\IR^{n+1}$.
In our case of the $\Ds$-hypergeometric system,
the finite set is given by the set of all integral points in the polyhedron
$\Ds$. Namely we have $\A
=\{ \bn0,\bn1,\cdots,\bn{p} \; \vert \;
\bn{i}=(1,\nu_i^*),\; \nu_i^* \in \Ds\cap \ZZ^n \;\}$. Here we let
$\bn0=(1,\nu_0^*)$ for the origin $\nu_0^*$ in $\Ds$.
We consider a lattice $L$ of affine dependencies on $\A$:
\eqn\latticeL{
L=\{ (l_0,l_1,\cdots,l_p)\in \ZZ^{p+1} \; \vert \;
l_0\bn0+l_1\bn1+\cdots l_p \bn{p} =0 \}.
}
Then it is found in \Hdiff\ that the period integral \period\ satisfies
the following set of differential equations, ($\A$-hypergeometric
system with exponents $\beta=(-1,0,\cdots,0)\in \IR^{n+1}$),
\eqn\GKZeq{
\sqbox_l \Pi(a)=0 \;\; (l\in L) \;\;, \;\;
{\cal Z}_i \Pi(a)=0 \;\; (i=0,1,\cdots,p) ,
}
where the differential operators $\sqbox_l$ and ${\cal Z}_i$
are defined to be
\eqn\BoxOp
{
\eqalign{
\sqbox_l&=\prod_{l_i >0} \( {\pd \; \over \pd a_i }\)^{l_i} -
\prod_{l_j<0} \({\pd \; \over \pd a_j}\)^{-l_j} \quad (l\in L) \cr
{\cal Z}_j &= \sum_{i=0}^p \bn{{i,j}} \ta{i} -\beta_j \quad
(j=0,1,\cdots,n)\;.        \cr
}}
The solution space of
\GKZeq\ is typically too large -- it contains more than the period
integrals of the Calabi-Yau manifolds $X_\Ds$. It turns out that
the period integrals satisfy additional differential equations.

\vskip0.5cm

\subsec{Automorphism of $\IP_\SD$ }

It is easy to recognize the origin of the linear differential
operators ${\cal Z}_j \; (j=1,\cdots,n)$ as the invariance of the
period integral \period\ under the canonical torus action on a
toric variety, $X_i \rightarrow \lambda_i X_i \; (\lambda_i \in \IC^*)$.
Since the algebraic torus acts by a subgroup of the
automorphism group of the toric variety $\IP_{\SD}$, it is natural
to incorporate into the PDE system the invariance under infinitesimal action of
the full automorphism group.
To describe this action in full generality, we will introduce the root
system for a toric variety.

Let us consider a compact nonsingular toric variety $\IP_\S$
based on a regular
fan $\S$ in the scalar extension $N_{\IR}$ of a lattice $N\;(\cong\ZZ^r)$
of rank $r$.
Let $M\;(\cong \ZZ^r)$ be the lattice dual to $N$. We
choose a basis $\{ n_1,\cdots,n_r\}$ for $N$ and a dual basis
$\{ m_1,\cdots,m_r\}$ for $M$. There is
a canonical algebraic torus $T_N:={\rm Hom}_\ZZ (M,\IC^*)=(\IC^*)^r$
in $\IP_\S$ whose coordinate ring is
$\IC[M]=\displaystyle{\oplus_{m\in M}}\IC {\bf e}(m)$. We write it as
$\IC[X_1^{\pm1},\cdots,X_r^{\pm 1}]$ with $X_i={\bf e}(m_i)$.
Define the derivations $\delta_n \; (n\in N)$ on
$\IC[M]$ by $\delta_n {\bf e}(m)=\langle m,n \rangle {\bf e}(m)$.
These derivations describe the natural action of $Lie(T_N)$ on $T_N$.
We may write $\{ \delta_{n_1}, \cdots, \delta_{n_r} \}
= \{ X_1 {\pd \; \over \pd X_1}, \cdots, X_r {\pd \; \over \pd X_r} \}$.
The Lie algebra of the full automorphism group of $\IP_\S$ is described by the
root system
$R(\S)$ in addition to the torus action.
The root system $R(\S)$ is determined by the data of the fan $\S$ as follows.
We denote the subset of one dimensional cones in the fan as
$\S(1)$. In each one dimensional cone $\sigma^{(1)} \in \S(1)$, there
is a primitive element $n(\sigma^{(1)})$ in $N$. Let
\eqn\rootsys{
\eqalign{
R(\S)=\{ \alpha\in M \; \vert \;
\exists \sigma_\alpha^{(1)} \in \S(1)\;
& {\rm with} \langle \alpha, n(\sigma_\alpha^{(1)}) \rangle =-1 \;  \cr
&
{\rm and } \;
\langle \alpha, n(\sigma^{(1)}) \rangle \geq 0  \; {\rm for} \, {\rm all}\,
\sigma^{(1)}\not= \sigma_\alpha^{(1)} \; \}. \cr}
}
In terms of the root system, the Lie algebra of the automorphism group can be
expressed by
\eqn\auto{
Lie(Auto(\IP_\S)) =
Lie(T_N)\oplus \( \oplus_{\alpha\in R(\S)}
\IC {\bf e}(\alpha) \delta_{n(\sigma_\alpha^{(1)})} \).
}
The linear differential operators ${\cal Z}_1,\cdots, {\cal Z}_n$ in \GKZeq\
express the invariance of the period integral $\Pi(a)$ under the action of
$Lie(T_N)$. In fact it is easy to check that
\eqn\lieaction{
{\cal Z}_i \Pi(a) =
\int_{C_0} \delta_{n_i} \( {1\over f_\Ds (X,a)} \)
\prod_{k=1}^n {dX_k \over X_k} \;\;\; (i=1,\cdots,n).
}
The operator ${\cal Z}_0$ represents the change of the period under the
overall scaling of the Laurent polynomial $f_\Ds(X,a) \rightarrow
\lambda f_\Ds(X,a)$. We can now clearly extend the formula \lieaction~ to
define ${\cal Z}_Y\Pi(a)$ for every $Y\in Lie(Auto(\IP_\SD))$
by replacing $\delta_{n_i}$ by $Y$.
We thus arrive at the definition of the extended $\Ds$-hypergeometric system
\eqn\extendedGKZ{
\sqbox_l \Pi(a)=0 \;\; (l \in L) \;\;,\;\;
{\cal Z}_Y \Pi(a)=0 \;\; (Y\in Lie(Auto(\IP_\SD)) ).
}
This extended system was first introduced in \HKTYI\ and was used
successfully  to determine the complete set of the period integrals.

Because of the special value of the exponent $\beta=(-1,0,\cdots,0)\in
\IR^{n+1}$, the following gauge for the period
\eqn\tgauge{
\tilde \Pi(a)=a_0 \Pi(a) \;\; ,
}
will be useful. We will denote the hypergeometric system in this gauge
as $\tilde \sqbox_l \tilde \Pi(a)=0 \;,\; \tilde {\cal Z}_i\tilde
\Pi(a)=0$.
Especially the first order differential operators $\tilde{\cal Z}_0,
\tilde{\cal Z}_1,\cdots,\tilde{\cal Z}_n$ may be written concisely as
\eqn\PSdiff{
\tilde{\cal Z}_u=\sum \langle u, \bn{i} \rangle \ta{i} \quad
(u\in \IR^{n+1} ).
}

\vskip0.5cm

In ref.\HKTYI\ , several Calabi-Yau hypersurfaces with $h^{1,1}=2$ and
3 have been studied. There hypersurfaces in a  weighted projective
space have been classified into three types depending on the properties
of the fan $\SDs$. Type I models are those which do
not have any integral points in the interior of
codimension-one faces of $\Ds(w)$ and
for which we have a regular fan $\SDs$ after taking into account
subdivisions of the cones resulting from the integral points on the
lower dimensional faces.
Type II models are those which have integral points in the interior of
codimension-one faces of $\Ds(w)$ but for which we still have a
regular fan $\SDs$ after subdivisions of the cones resulting
from the integral points on the faces.
Type III models are those which we do not have a regular fan $\SDs$
even if we subdivide the cones by incorporating all the integral points
on the faces.  In this sense type III models may be called 'singular'.
According to this classification, we reproduce here
the models
analyzed in \HKTYI\
\eqn\models{
\eqalign{
&{\rm Type}\;{\rm I}: X_8(2,2,2,1,1) \cr
&{\rm Type}\;{\rm II}:
X_{12}(6,2,2,1,1) \;,\;
X_{14}(7,2,2,2,1) \;,\;
X_{18}(9,6,1,1,1) \;,\;
X_{12}(6,3,1,1,1) \;,\; \cr
& \hskip1cm
X_{24}(12,8,2,1,1) \;,\; \cr
&{\rm Type}\;{\rm III}:
X_{12}(4,3,2,2,1) \;,\;
X_{12}(3,3,3,2,1) \;,\;
X_{15}(5,3,3,3,1) \;,\;
X_{18}(9,3,3,2,1) \;.\; \cr
}}
It was found that for a model of type I or II, the extended
$\Ds$-hypergeometric system is sufficient to determine the complete
set of the period integrals. Whereas for models of type III,
one needs to consider additional (non-toric)
differential operator(s) whose form can be determined from the
Jacobian ring of the hypersurface. If we supplement these additional
operators to the extended $\Ds$-hypergeometric system,
we can derive the Picard-Fuchs differential equations.
Thus for type III models, the combinatorial data of the polyhedron
$\Ds$ alone do not seem sufficient for the explicit construction of
 the full system of differential operators.
Nevertheless we will find in the next section
that the local solutions are determined purely by the combinatorial
data of the polyhedron, and this property is shared by all three
types of the Calabi-Yau hypersurfaces.

\vskip0.5cm

\leftline{\undertext{Example: $X_{14}(7,2,2,2,1)$ }}

This is a typical model with non-trivial automorphism group. The polyhedron
$\D(w)=\{ x\in \IR^5 \;\vert\;
w_1x_1+\cdots+w_5x_5=0\;,\; x_i\geq -1 \; (i=1,\cdots,5) \; \}$
is simplicial and is given by the convex hull of the vertices
\eqn\exImus{
\eqalign{
&
\nu_1=(1,-1,-1,-1) \;,\;
\nu_2=(-1,6,-1,-1) \;,\;
\nu_3=(-1,-1,6,-1) \;,\; \cr
&
\nu_4=(-1,-1,-1,6) \;\;,
\nu_5=(-1,-1,-1,-1) \;,  \cr
}}
where we fix a basis $\{ \Lambda_1,\cdots,\Lambda_4\}$ for
the lattice $H(w)=\{ x\in \ZZ^5 \;\vert\; w_1x_1+ \cdots + w_5x_5=0 \; \}$,
with
$\Lambda_1=(1,0,0,0,-w_1), \Lambda_2=(0,1,0,0,-w_2),
\Lambda_3=(0,0,1,0,-w_3)$ and $ \Lambda_4=(0,0,0,1,-w_4) $.
The integral points in the dual polyhedron
$\Ds(w)$ are
\eqn\exInus{
\eqalign{
&
\nu_0^*=(0,0,0,0) \;,\;
\nu_1^*=(1,0,0,0) \;,\;
\nu_2^*=(0,1,0,0) \;,\;
\nu_3^*=(0,0,1,0) \;,\; \cr
&
\nu_4^*=(0,0,0,1) \;,\;
\nu_5^*=(-7,-2,-2,-2) \;,\;
\nu_6^*=(-3,-1,-1,-1) \;,\; \cr
&
\nu_7^*=(-4,-1,-1,-1) \;,\;
\nu_8^*=(-1,0,0,0) \;.\; \cr
}}
The points $\nu_1^*, \cdots, \nu_5^*$
are the vertices of the simplicial polyhedron $\Ds(w)$ and all other points
(except the origin)
appear on some faces of the polyhedron. The point $\nu_6^*={1\over2}(
\nu_1^*+\nu_5^*)$ appears on the edge (one dimensional face) and corresponds to
an exceptional divisor in $X_\D$. The point
$\nu_7^*={1\over7}(\nu_2^*+\nu_3^*+\nu_4^*+7\nu_5^*)$ and
$\nu_8^*={1\over7}(2\nu_2^*+2\nu_3^*+2\nu_4^*+\nu_5^*)$ are both
in the interior of the
codimension-one face dual to the corner $\nu_1$ of $\D(w)$. Hence they
describe the automorphism of $\IP_\SD$ and of the family of
hypersurfaces $X_\Ds$.
In fact the two points describe
the root system for the fan $\SD$ and generate the nontrivial part
of the automorphism,
\eqn\exIxis{
\IC \xi_1\oplus \IC \xi_2 :=
\IC {\rm e}(\nu_7^*)\delta_{\nu_1}\oplus
\IC {\rm e}(\nu_8^*)\delta_{\nu_1} \;.
}
These infinitesimal actions on the coordinate ring can be expressed in
terms of the natural basis for $N=\ZZ^4$ and $M=\ZZ^4$ as
$\xi_i=X^{\nu_{6+i}^*} (\delta_{X_1}-\delta_{X_2}-
\delta_{X_3}-\delta_{X_4}) \; (i=1,2)$ and have the expressions
\eqn\exIliexi{
\eqalign{
\xi_1
&={1\over X_1^4X_2X_3X_4} \( X_1{\pd \; \over \pd X_1} -
   X_2 {\pd \; \over \pd X_2} - X_3{\pd \; \over \pd X_3} -
   X_4{\pd \; \over \pd X_4} \) \;\;, \cr
\xi_2
&={1\over X_1} \( X_1{\pd \; \over \pd X_1} -
   X_2 {\pd \; \over \pd X_2} - X_3{\pd \; \over \pd X_3} -
   X_4{\pd \; \over \pd X_4} \) \;. \cr
}}
We may verify the algebra $[\xi_1,\xi_2]=0$. The linear
differential operators ${\cal Z}_{\xi_1}$ and ${\cal Z}_{\xi_2}$, which
follows from \lieaction\ , turns out to be
\eqn\exIzxis{
\eqalign{
{\cal Z}_{\xi_1}&=a_0{\pd \;\over \pd a_7} +2a_1{\pd \; \over \pd a_6}
+ a_6 {\pd \; \over \pd a_5} \cr
{\cal Z}_{\xi_2}&=a_0{\pd \;\over \pd a_8}+2a_1{\pd\; \over \pd a_0}
+a_6{\pd \; \over \pd a_7}\;. \cr
}}
These linear operators together with ${\cal Z}_0, \cdots, {\cal Z}_5$
and the higher order operator $\sqbox_l \; (l\in L)$ constitute the
full extended $\Ds$-hypergeometric system.

\vskip0.5cm

\newsec{ Secondary fan, Gr\"obner fan and local solutions }

In this section, we analyze the local solutions of the
$\Ds$-hypergeometric system. We find that the local properties of
the $\Ds$-hypergeometric system are determined purely by
an algebro-combinatoric object, known as a toric ideal. At a special point,
called "large radius limit", the toric ideal is related
to an ideal which determines the cohomology ring of the
toric variety $\IP_\SDs$.

\vskip0.3cm

\subsec{Convergent series solutions for $\A$-hypergeometric system }

Here we will summarize, with some modification,
the general results in \GKZ\ about the convergent
series solutions of the $\A$-hypergeometric system. We set
$\A=(1,\Ds)\cap \ZZ^{p+1} =\{ \bn{0},\cdots,\bn{p} \}$ for our case
of the $\Ds$-hypergeometric system. The description here is brief
and is meant to fix notations and to prepare for later discussions.
We refer the reader to the original paper \GKZ\ for details.

{}From the definition of the $\A$-hypergeometric system \GKZeq\BoxOp, it is
easy to check that a formal solution to the $\A$-hypergeometric system with
exponent $\beta\in \IR^{n+1}$ is given by
\eqn\fseries{
\Pi(a,\gamma)=\sum_{l\in L}
{ 1 \over
  \prod_{0\leq i\leq p} \Gamma(l_i+\gamma_i+1) }  a^{l+\gamma} \;\;,
}
where $\beta=\sum_i \gamma_i\bn{i}$. Evidently the formal
solution is invariant under $\gamma
\rightarrow \gamma+v \; (v\in L)$. Define the affine subspace
$\Phi(\beta):=\{\g{}\in \IR^{p+1} \vert \beta=\sum \g{i}\bn{i}\}$.
If we choose a basis $l^{(1)},..,l^{(p-n)}$ for $L$,
the formal series \fseries\
takes the form $\Pi(a,\g{})=a^\g{}\sum_{m_1,\cdots,m_{p-n}\in\ZZ}c_m x^m$,
where $x_k=a^{l^{(k)}}$.
The relevant solutions are those with $c_m(\gamma)=0$ unless
$m_i\geq0$. One must therefore restrict the choices of the basis and of
$\gamma$.

A subset $I\in \{0,1,\cdots,p\}$ is a base if
$\{\bn{i} \vert i\in I\}$ form a basis of $\IR^{n+1}$. Given
a base $I$ and $\g{j} \;(j\notin I)$,
 we can solve for $\g{j} \;(j\in I)$
using the linear relation $\sum_{j\in I}\g{j}\bn{j}=
\beta-\sum_{j\notin I}\g{j}\bn{j}$. Consider
$\Phi_\ZZ(\beta,I):= \{\g{} \in \Phi(\beta) \vert \g{j}\in \ZZ \; (j\notin
I)\}$, and $\Phi_\ZZ^A(\beta,I):= \{ \g{} \in
\Phi_\ZZ(\beta,I) \vert \g{j}=\sum_{k=1}^{p-n} \lambda_k l^{(k)}_j \;
(0\leq \lambda_k < 1, j\notin I) \}$. It is clear that
$\Phi_\ZZ^A(\beta,I)$ is a set of representatives of
$\Phi_\ZZ(\beta,I)/L$.
Consider the cone ${\cal K}(\A,I)=\{ l \in L_\IR \vert l_i\geq 0
(i\notin I)\}$ where $L_\IR=L\otimes \IR$. A $\ZZ$-basis $A\subset L$ is
said to be compatible with the base $I$ if the cone generated by the basis $A$
contains the cone ${\cal K}(\A,I)$.

\vskip0.5cm
\leftskip=0.8cm \rightskip=1.5cm
\vbox{
\noindent
{\it
If $A=\{l^{(1)},\cdots,l^{(p-n)}\}$ is compatible with the base $I$,
then the formal series \fseries\
takes the form $\Pi(a,\g{})=a^\g{}\sum_{m_1,\cdots,m_{p-n}\geq 0}c_m x^m$ for
each $\g{}\in \Phi_\ZZ^A(\beta,I)$ with $x_k=a^{l^{(k)}}$, and this power
series converges for sufficiently small $\vert x_k\vert$. }
\vskip-60pt
\eqn\convs{ }
\vskip30pt
}
\leftskip=0cm \rightskip=0cm
\noindent
By definition we may write the formal series \fseries\ as above with
$c_m=c_m(\gamma):=\prod_{i=0}^{p} 1/\Gamma(\sum m_k
l_i^{(k)}+\g{i}+1)$. For $\g{}\in \Phi_\ZZ^A(\beta,I)$,
we have $\sum m_k l_j^{(k)} +\g{j}+1 \in \ZZ$ for $j\notin
I$. It follows that if $c_m\neq0$, then $\sum m_k l_j^{(k)}
+\g{j} = \sum (m_k+\lambda_k)l_j^{(k)} \geq 0 (j\notin I)$ where we use
$\g{j}=\sum \lambda_k l_j^{(k)} \; (0\leq \lambda_k <1, j\notin I)$.
Since the basis $A$ is compatible with the base $I$, we have
$m_k+\lambda_k \geq 0$ for all $k$, implying $m_k\geq 0$.
Thus given a basis $A$ compatible with the base $I$,
if for every $\gamma\in \Phi_\ZZ^A(\beta,I)$ there is $c_m(\gamma)\neq0$
for some $m$, then we have $\vert \Phi_\ZZ(\beta,I)/L\vert =\vert {\rm
det}(\bn{j,i})_{1\leq  i \leq n+1, j\in I} \vert $ linearly independent
power series solutions \GKZ.

However it can happen that $c_m(\gamma)=0$ for all $m$, i.e., the
series solution becomes trivial
$\Pi(a,\g{}) \equiv 0$ when $\sum m_k l_i^{(k)} +\g{i} +1
\in \ZZ_{\leq 0}\; (m_k\geq 0)$ for some $i \in I$.
In this case, we multiply $c_m$ by a constant infinite
renormalization $\Gamma(\g{i}+1)$. More precisely,
we assume that $\gamma$ is such that the following limit exists:
\eqn\climit{
\eqalign{
{\Gamma(\g{i}+1) \over \Gamma(l_i+\g{i}+1) }
& :=
 {{\rm lim} \atop {}^{\varepsilon \rightarrow 0}}
{\Gamma(\g{i}+1+\varepsilon) \over
 \Gamma(l_i+\g{i}+1+\varepsilon) }
}}
for all $l\in L$.

All linearly independent power series solutions are constructed from a
set of bases $\{I\}$ which form a triangulation of the polyhedron
$P:={\rm Conv.}(\{0,\bn{0},\bn{1},\cdots,\bn{p}\})$, where $0$ is the
origin in $\IR^{n+1}$. We call a collection of bases $T=\{I\}$ a
triangulation of $P$ if $\cup_{I\in T} \langle \bn{I}\rangle =P$ and $
\langle \bn{I_1}\rangle \cap \langle \bn{I_2}\rangle \; (I_1,I_2 \in T)$ is
a lower dimensional common face. Here $\langle
\bn{I}\rangle$ a $n+1$ dimensional simplex with vertices
$\bn{i} (i\in I)$ and the origin. Because the $n+1$-simplices in $P$
are in 1-1 correspondence with the $n$-simplex in $\Ds$, there is a
notion of a triangulation of $\Ds$( or $\A$). We use the two notions
interchangeably.
A triangulation $T$ is called maximal if $T$ gives
 the maximum number of
$n$-simplices in $\Ds$ and $0\in I$ for all $I\in T$.
A $\ZZ$-basis $A$ of $L$ is called compatible with a
triangulation $T$ if $A$ is compatible with every $I\in T$.

For a base $I$ and a point $\eta \in \IR^{p+1}$, we consider a linear
function $h_{I,\eta}$ on $\IR^{p+1}$ such that
$h_{I,\eta}(\bn{i})=\eta_i \; (i\in I)$. We define a cone ${\cal
C}(\A,I)$ by $\{ \eta\in\IR^{p+1} \vert h_{I,\eta}(\bn{i})\leq \eta_i \;
(i\notin I) \}$. For a triangulation $T$, we define the cone ${\cal
C}(\A,T):=\cap_{I\in T} {\cal C}(\A,I)$. We may associate with $\eta \in
\IR^{p+1}$ and a triangulation $T$, a piecewise linear continuous function
$h_{T,\eta}$ on the polyhedron $P$ defined by 1)
$h_{T,\eta}(\bn{i})=\eta_i$ for each vertex $\bn{i}$ of the
triangulation $T$, 2) the restriction $h_{T,\eta}\vert_{\langle \bn{I}
\rangle} \; (I\in T)$ is a linear function. Then the cone ${\cal
C}(\A,T)$ consists of $\eta \in \IR^{p+1}$ for which the function
$h_{T,\eta}$ is convex and $h_{T,\eta}(\bn{i})\leq \eta_i$ for
$\bn{i}$ not a vertex of $T$ \GKZ. A regular triangulation is a
triangulation for which we have interior points in the cone ${\cal
C}(\A,T)$. For every regular triangulation $T$,
there are infinitely many $\ZZ$-basis of $L$ compatible with $T$.
We set $\Phi_\ZZ^A(\beta,T):=\cup_{I\in T} \Phi_\ZZ^A(\beta,I)$.
Now we may
state the result (theorem 3) in \GKZ;

\vskip0.5cm
\leftskip=0.8cm \rightskip=1.5cm
\vbox{
\noindent
{\it
For a regular triangulation $T$ of the polyhedron $P$, and a $\ZZ$-basis
$A=\{ l^{(1)},\cdots,l^{(p-n)} \}$ of $L$ compatible with $T$, we have
integral power series in the variables $x_k=a^{l^{(k)}}$ for
$a^{-\g{}}\Pi(a,\g{}) \; (\g{}\in \Phi_\ZZ^A(\beta,T))$, which converge
for sufficiently small $\vert x_k\vert$. If the exponents $\beta$ is
$T$-nonresonant, the series $\Pi(a,\g{})\; (\g{}\in \Phi_\ZZ^A(\beta,T))$
constitute vol($P$) linearly
independent solutions for \GKZeq.
}
\vskip-80pt
\eqn\GKZth{ }
\vskip50pt
}
\leftskip=0cm \rightskip=0cm
\noindent
In the above theorem,
the exponent $\beta$ is called $T$-nonresonant if the sets
$\Phi_\ZZ(\beta,I) \;(I\in T)$ are pairwise disjoint. It turns out that
in our $\Ds$-hypergeometric system there are many regular triangulations
for which the exponent $\beta=(-1,0,\cdots,0)$ is $T$-resonant.
In particular, if $T$ is a maximal triangulation and the polyhedron $\Ds$
is of type I or II, then
$\beta$ is 'maximally $T$-resonant', i.e., $\Phi_\ZZ(\beta,I)$ consists
of a unique element $\g=(-1,0,\cdots,0)$ modulo $L$ for all $I\in T$.
(Note that each simplex $I\in T$ has volume
$\vert {\rm det}(\bn{j,i})_{1\leq  i \leq n+1, j\in I} \vert=1$.)
In this case, we will obtain only one power series solution \fseries,
 and all other solutions contain logarithms, whose forms we will
determined by the Frobenius method.

Given a regular triangulation $T$, a compatible $\ZZ$-basis
$A=\{l^{(1)},..,l^{(p-n)}\}$ and
$\gamma\in \Phi_\ZZ^A(\beta,T)$,
we define a power
series $w_0(x,\rho)=a_0\Pi(a,\g{})$ where $\rho=(\rho_1,\cdots,\rho_{p-n})$
is defined by $\gamma=\sum\rho_kl^{(k)}+(-1,0,\cdots,0)$ and
$x_k=(-1)^{l_0^{(k)}} a^{l^{(k)}}$. Explicitly we have
\eqn\wnot{
w_0(x,\rho)=\sum_{m_1,\cdots,m_k \geq 0}
  {\Gamma(-\sum ((m_k+\rho_k) l^{(k)}_0 +1) \over
   \prod_{1\leq i \leq p} \Gamma(\sum (m_k+\rho_k) l^{(k)}_i +1) }
   x^{m+\rho} \;
}
The $\rho$ can also be determined by the indicial equations of
the hypergeometric system.

Given a regular triangulation $T$, we shall now construct
a compatible $\ZZ$-basis $A$ with the criterion that
the cone generated by $A$ in $L_\IR$ contains the
cone ${\cal K}(\A,T):=\cup_{I\in T}{\cal K}(\A,I)$.
First we introduce the Gale transformation. Consider
exact sequence
\eqn\gale{
0 \longrightarrow \; \IR^{n+1} \;
\overrightarrow{\quad \;_{\bf A} \quad} \;  \IR^{p+1} \;
\overrightarrow{\quad \;_{\bf B} \quad} \;  \IR^{p-n} \;
\longrightarrow 0 \;\;,
}
where we let $\IR^{n+1}$ be the span of the integral points
$\bn{i}$'s and $\IR^{p+1}$ in the middle is the span of a basis $\{
\en{0},\en{1},\cdots, \en{p} \}$ labeled by the points. The linear map
${\bf A}$ sends $v \in \IR^{n+1}$ to $\sum (v,\bn{i})\en{i}$ and ${\bf
B}$ is the natural map onto $\IR^{p+1}/\IR^{n+1}=\IR^{p-n}$.
The Gale transform of a point
configuration $\A$ in $\IR^{n+1}$, which we denote $\{ \A,\IR^{n+1}\}$,
is defined to be a point configuration $\{ \B,\IR^{p-n}\}$ with $\B:= \{
{\bf B}(\en{0}),\cdots,{\bf B}(\en{p}) \}$. Now we consider a cone in
$\IR^{p-n}$,
\eqn\cp{
{\cal C}'(\A,T)={\cap_{I\in T}}
 \( \sum_{i\notin I} \IR_{\geq 0}
{\bf B}(\en{i}) \) \;\;.
}
Then it is shown in \ref\strumsecfan{L.Billera,P.Filliman and B.Strumfels,
Adv.in Math.{\bf 83}(1990),155-179}
\ref\OdaPark{T.Oda and H.S.Park, T\^ohoku Math. J. {\bf 43}(1991),375-399.}
that the cone ${\cal C}(\A,T)$
decomposes into $\IR^{n+1}\oplus {\cal C}'(\A,T)$. The secondary fan
${\cal F}(\A)$ is defined as
\eqn\secfan{
{\cal F}(\A)=\{ {\cal C}'(\A,T) \vert T: {\rm regular}\; {\rm
triangulation} \}.
}
It is known that the secondary fan is complete and strongly polytopal
polyhedral fan\strumsecfan\OdaPark.

In our point configuration $\{\A,\IR^{n+1}\}$, the set $\A$ consists of
integral points. Therefore the sequence \gale\ can be endowed with
an integral structure:
$
0 \rightarrow \; \ZZ^{n+1} \;
\overrightarrow{\; \;_{\bf A} \;} \;  \ZZ^{p+1} \;
\overrightarrow{\; \;_{\bf B} \;} \;  $ $ \ZZ^{p-n} \;
\rightarrow 0 \;\;
$.
The dual of this sequence is
$
0 \leftarrow \; (\ZZ^{n+1})^* \;$$
\overleftarrow{} \;  (\ZZ^{p+1})^* \;$$
\overleftarrow{} \;  L \;$$
\leftarrow 0 \;\;
$, where $L$ is the lattice of the affine relations among $\A$. The
cone dual to ${\cal C}'(\A,T)\subset \IR^{p-n}$ is the cone
${\cal K}(\A,T)\subset L_\IR$. In general ${\cal C}'(\A,T)$ is strongly
convex but not necessarily regular. There is a canonical refinement
of the secondary fan known as the Gr\"obner fan (see next subsection).
However even a cone in this refinement is not necessarily regular.
By suitably subdividing the cone, we obtain a regular subcone and
hence a $\ZZ$-basis of this subcone.
The dual basis $A=\{ l^{(1)},..,l^{(p-n)}\}$
thus generates a cone containing ${\cal K}(\A,T)$. This gives
us a $\ZZ$-basis of $L$ compatible with $T$.
Note that when ${\cal C}'(\A,T)$ is already regular, the basis $A$ is uniquely
determined by $T$.

Suppose now the polyhedron $\Ds$ is of type I or II. Then endowed with
a maximal subdivision, it defines a regular fan
$\SDs$ and $\IP_\SDs$ is smooth. It is known that
the Gale transform $\{ {\cal B},\ZZ^{p-n}\}$
generates the Picard group
\ref\Oda{T.Oda, {\it
Convex bodies and Algebraic Geometry, An Introduction to the Theory of
Toric Varieties}, A Series of Modern Surveys in Mathematics (1985)Springer-
Verlag.}
\ref\Fulton{
W.Fulton,{\it Introduction to Toric Varieties}, Ann.of Math.Studies 131,
Princeton University Press(1993).}.
Now associated with $\SDs$
is a maximal triangulation $T$. In this case,
${\cal C}'(\A,T)$ is the K\"ahler cone, and ${\cal K}(\A,T)$
is the Mori cone of $\IP_\SDs$.
In particular ${\cal C}'(\A,T)$ is a maximal cone (hence has interior points),
implying that that $T$ is regular.

\vskip0.3cm

\subsec{ Toric ideal and universal Gr\"obner basis }

Here we will focus on the differential operators $\sqbox_l \; (l\in L)$ in
\BoxOp. Although there are infinitely many operators, we can describe
the system by a finite set of the operators. The problem is
how to construct such a finite set. We will see
that the so-called toric ideal in the theory of Gr\"obner basis
gives us a powerful  tool for this purpose.

Let $\A$ be the finite set in the previous subsection and $L$ be the lattice
representing the integral relations among the vertices in $\A$. We may
decompose any element $l \in L$ uniquely into $l_+-l_-$ with two nonnegative
vectors $l_+, l_-$ having disjoint support, where the support $m\in L$
is defined by supp$(m):=\{\;i\;\vert\; m_i\not= 0\;\}$. Toric ideal
${\cal I}_\A$ is defined as the ideal in $\IC[y_0,\cdots,y_p]$ which is
generated by $y^{l_+}-y^{l_-}$, i.e.,
\eqn\tideal{
{\cal I}_\A = \langle y^{l_+}-y^{l_-} \; \vert \; l\in L \rangle \;\;.
}
Let $\omega$ be a term order on
in $\IC[y_0,\cdots,y_p]$. It is a vector
$\omega=(\omega_0,\cdots, \omega_p) \in \IR^{p+1}$ which defines
an monomial ordering by the weights: the weight of
$y_0^{\alpha_0}\cdots y_p^{\alpha_p}$ being $\omega_0 \alpha_0 +
\cdots + \omega_p \alpha_p$. With respect to this term order,
we consider an ideal $LT_\omega ({\cal I}_\A)=
\langle LT_\omega(f) \; \vert \; f \in {\cal I}_\A \rangle$
of the leading terms of
${\cal I}_\A$. Two different weights $\omega$ and $\omega'$ may give the
same ideal. The equivalence class in $\IR^{p+1}$
\eqn\gfan{
{\cal C}({\cal I}_\A,\omega):=
\{ \omega' \in \IR^{p+1} \; \vert \;
LT_\omega({\cal I}_\A)=LT_{\omega'}({\cal I}_\A) \; \} \;\;,
}
is an open convex polyhedral cone.
The collection of cones $\{ {\cal C}({\cal I}_\A,\omega) \}$ is known to
be finite and
defines a polyhedral fan called the Gr\"obner fan ${\cal F}({\cal I}_\A)$
of ${\cal I}_\A$.

The Gr\"obner basis of ${\cal I}_\A$ with respect to a term order
$\omega$ is a finite generating set
${\cal B}_\omega$ of ${\cal I}_\A$, with the property that the ideal
$\langle LT_\omega(g) \; \vert \; g \in {\cal B}_\omega \rangle$
is equal to $LT_\omega({\cal I}_\A)$.
By Hilbert's basis theorem, ${\cal I}_\A$ is generated by
a finite set binomials $y^{l_+}-y^{l_-}$ with $l\in L$.
Starting from such a finite set,
the (reduced) Gr\"obner basis ${\cal B}_\omega$ obtained by
Buchberger's algorithm \ref\cox{see for example, D. Cox, J. Little and
D. O'shea, {\it Ideals, Varieties, and Algorithms}, UTM(1991),
Springer-Verlag.} is also a set of binomials. This is because the algorithm
consists of forming the $S$-polymonials for the generators and the
reductions of the minimal Gr\"obner basis and both processes close in the
set of binomials. Moreover the elements of the reduced Gr\"obner basis take
the form $y^{l_+}-y^{l_-}$ ($l\in L$) of binomials.

Next given a term order $\omega$, we shall obtain a regular triangulation
$T_\omega$ and hence a compatible $\ZZ$-basis $A$ of $L$ (last section).
The elements of the toric ideal ${\cal I}_\A$ may be identified
as differential operators which annihilate the formal
series $\Pi(a,\g{})$ with $\g{}\in\Phi^A_\ZZ(\beta,T_\omega)$.
The ideal $LT_\omega({\cal I}_\A)$ is then a set of
of 'leading' terms of the operators which determine the indices for the series
$w_0(x,\rho)$.
Therefore the Gr\"obner basis ${\cal B}_\omega$ with respect to
 $\omega$ gives a finite set of the differential
operators $\{ \sqbox_\l \}$ which suffices to describe
the local solutions.
A finite set which contains the
Gr\"obner basis ${\cal B}_\omega$ for all term orders is known as a universal
Gr\"obner basis ${\cal U}_\A$. This basis is useful to describe global
property of the system.

A nonzero integral relation $l\in L$ is called elementary if 1) $l$ is
primitive, i.e., gcd$(l_0,l_1,\cdots,l_p)=1$, 2) supp$(l)$ is minimal with
respect to inclusion. It is known that the set $\{l^{(1)},
l^{(2)}, \cdots,l^{(m)}\}$ of all elementary integral relations
 generates
a $(p+1)$-dimensional zonotope ${\cal P}_\A:=\langle 0,l^{(1)}\rangle+
\langle 0,l^{(2)}\rangle+\cdots + \langle 0,l^{(m)}\rangle$, where $\langle
0,l^{(k)}\rangle$ represents a one-dimensional simplex and the sum means
the Minkowski sum.  The universal Gr\"obner
basis is then given by
\ref\strum{B. Strumfels, T\^ohoku Math. J.{\bf 43}(1991)249.}\
\eqn\ugbasis{
{\cal U}_\A =
\{ y^{l_+}-y^{l_-} \; \vert \; l\in {\cal P}_\A \cap \ZZ^{p+1} \; \} \;\; .
}

Given a term order $\omega$ the notion of a regular triangulation
can in fact be recast as follows. Consider the polytope
$P_\omega:=$Conv.$\{ (\omega_0,\ns{0}),\cdots,(\omega_p,\ns{p})\}$ in
$\IR^{n+1}$, which is a lifting of $\A$ by assigning the weights
$\omega_i$ as height to each point $\ns{i}$.
For sufficiently generic $\omega$, the lower envelope of
$P_\omega$ naturally induces a triangulation $T_\omega$ of $\A$.

It turns out that a triangulation $T$ of $\A$ is regular if and only
$T=T_\omega$ for some generic
weight $\omega$. Also the interior points of the cone
${\cal C}(\A,T)$ consists of all weights $\omega \in \IR^{p+1}$ such that
$T_\omega =T$ \strum.
The Stanley-Reisner ideal $SR_T$ for a triangulation $T$ of $\A$ is
the ideal in $\IC[y_0,\cdots,y_p]$ generated by all monomials
$y_{i_1}y_{i_2}\cdots y_{i_k}$ where $\{i_1,i_2,\cdots,i_k\} \notin T$.
Then the following is shown in \strum,

\vskip0.5cm
\leftskip=0.8cm \rightskip=1.5cm
\vbox{
\noindent
{\it If a weight vector $\omega$ defines a term order for the
toric ideal ${\cal I}_\A$, then the corresponding subdivision $T_\omega$
is a regular triangulation.
The Stanley-Reisner  ideal $SR_{T_\omega}$ is
equal to the radical of the ideal
$LT_\omega({\cal I}_\A)$. }

\vskip-50pt
\eqn\thmI{ }
\vskip20pt
}
\leftskip=0cm \rightskip=0cm

As an immediate corollary to \thmI,  the Gr\"obner fan
${\cal F}({\cal I}_\A)$ is a refinement of the fan $\{{\cal C}(\A,T)\}$.
Since each cone ${\cal C}(\A,T)$ has a decomposition $\IR^{n+1}\oplus
{\cal C}'(\A,T)$, we have a similar decomposition
${\cal C}({\cal I}_\A,\omega)=\IR^{n+1}\oplus{\cal C}'({\cal I}_\A,\omega)$.
We will call the collection  $\{{\cal C}'({\cal I}_\A,\omega)\}$ the
Gr\"obner fan which we also denote by ${\cal F}({\cal I}_\A)$.

\vskip0.3cm

\leftline{\undertext{Example: $\IP(2,2,2,1,1)$} }

This is a simple example of a toric variety in which we can define a
Calabi-Yau hypersurface with $h^{1,1}(X_\D)=2$.
The polyhedron $\D(w)$ is given by the convex
hull of the following integral points,
\eqn\exIImus{
\eqalign{
&
\nu_1=(3,-1,-1,-1) \;,\;
\nu_2=(-1,3,-1,-1) \;,\;
\nu_3=(-1,-1,3,-1)  \cr
&
\nu_4=(-1,-1,-1,7) \;,\;
\nu_5=(-1,-1,-1,-1) \;,\; \cr } }
where the vector components are those with respect to a fixed basis
$\Lambda_1=(1,0,0,0,0,-w_1), \cdots, \Lambda_4=(0,0,0,1,-w_4)$ for
the lattice $H(w)$ (see the example in the previous section).
The integral points in the dual $\Ds(w)$ are
\eqn\exIInus{
\eqalign{
&
\nu_1^*=(1,0,0,0) \;,\;
\nu_2^*=(0,1,0,0) \;,\;
\nu_3^*=(0,0,1,0) \;,\;  \cr
&
\nu_4^*=(0,0,0,1) \;,\;
\nu_5^*=(-2,-2,-2,-1) \;,\;
\nu_6^*=(-1,-1,-1,0) \;,\;  \cr } }
The point $\nu_6^*={1\over2}(\nu_4^*+\nu_5^*)$ in a codimension 3
face of $\Ds$ corresponds to a $A_1$-type
Du Val singularity in the affine subvariety determined by the
cone $\IR_{\geq}\nu_4^*+\IR_{\geq}\nu_5^*$ in the fan $\SDs$.
 We can find three elementary
relations (up to sign) in $\A=(1,\Ds(w))\cap \ZZ^5$ expressed by
\eqn\exIIels{
\eqalign{
&
l^{(1)}=(-4,1,1,1,0,0,1) \;,\;
l^{(2)}=( 0,0,0,0,1,1,-2) \;,\; \cr
&
l^{(3)}=(-8,2,2,2,1,1,0) \;.\;  \cr }}
Then the zonotope ${\cal P}_\A=\langle 0, l^{(1)}\rangle +
\langle 0,l^{(2)} \rangle + \langle 0,l^{(3)} \rangle$ determines the
universal Gr\"obner basis
\eqn\exIIugbasis{
{\cal U}_\A=\{ y_1y_2y_3y_6-y_0^4 ,\; y_4y_5-y_6^2 , \;
y_1^2y_2^2y_3^2y_4y_5-y_0^8, \; y_1y_2y_3y_4y_5-y_0^4y_6 \; \}
}

It is straightforward to find all possible regular triangulations of the
set $\A$, or equivalently the polyhedron $\Ds(w)$.
We find the following four regular triangulations;
\eqn\exIItriang{
\eqalign{
&
T_0=\{ \langle 0,2,3,5,6 \rangle, \langle 0,1,3,5,6 \rangle,
       \langle 0,1,2,5,6 \rangle, \langle 0,2,3,4,6 \rangle,   \cr
&\hskip1cm
       \langle 0,1,3,4,6 \rangle, \langle 0,1,2,4,6 \rangle,
       \langle 0,1,2,3,5 \rangle, \langle 0,1,2,3,4 \rangle  \} \;\;, \cr
&
T_1=\{\langle 1,2,3,4,5 \rangle \} \;\;,\;\;
T_2=\{ \langle 1,2,3,4,6 \rangle, \langle 1,2,3,5,6 \rangle \} \;\;, \cr
&
T_3=\{ \langle 0,2,3,4,5 \rangle, \langle 0,1,3,4,5 \rangle,
       \langle 0,1,2,4,5 \rangle, \langle 0,1,2,3,5 \rangle,
       \langle 0,1,2,3,4 \rangle \}, \cr}
}
where, for example, $\langle 0,2,3,4,5 \rangle$ represents a simplex
with vertices $\bn{0},\bn{2}, \cdots, \bn{5}$.

\goodbreak\midinsert
\centerline{\epsfxsize 2truein\epsfbox{fig1.eps}}
\smallskip
\noindent{\tabfont
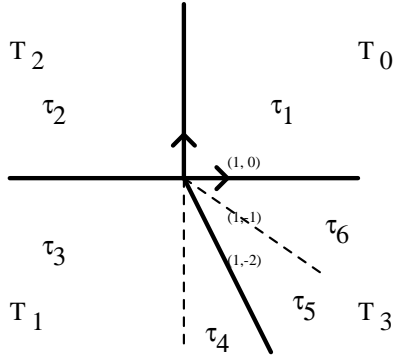
Fig.1 The secondary fan and the Gr\"obner fan for $\IP(2,2,2,1,1)$
\par\noindent
The secondary fan consists of the polyhedral cones parametrized by the
regular triangulations $T_0,\cdots,T_3$ in the text. The Gr\"obner fan
provides a refinement consisting of  $\tau_1,\cdots,\tau_6$
represented by the typical weights in the table1. }
\smallskip\endinsert

The Gr\"obner fan consists of six two-dimensional cones, together with
lower dimensional cones as their faces. We list the typical weight
with the corresponding ideal
$LT_\omega({\cal I}_\A)$ and its radical in the table 1.
We draw, in the figure 1, the secondary fan and the Gr\"obner fan
as its refinement using a $\ZZ$-basis $\{ \tilde l^{(1)},
\tilde l^{(2)} \}$ which is dual to a $\ZZ$-basis $\{ l^{(1)},
l^{(2)} \}$ in \exIIels\ of the lattice $L$.

\vskip0.2cm
\vbox{{\offinterlineskip
\hrule
\halign{
& \vrule# & \quad # \hfill &  \quad \hfil # \hfill  &   \quad\hfil # \hfill & #
\cr
height2pt
& \omit    & \omit    & \omit &  \omit &\cr
 & cone
  & weight $\omega$
   & $LT_\omega({\cal I}_\A)$
    & rad($LT_\omega({\cal I}_\A)$)
     &\cr
height2pt & \omit & \omit & \omit & \omit &\cr
\noalign{\hrule}
height2pt & \omit & \omit & \omit & \omit &\cr
 & $\tau_1$
  & $(0,1,1,1,1,1,0)$
   & $\langle y_1y_2y_3y_6,  y_4y_5 \rangle $
    & $\langle y_1y_2y_3y_6,  y_4y_5 \rangle $
     & \cr
height2pt & \omit & \omit & \omit & \omit &\cr
 & $\tau_2$
  & $(1,0,0,0,1,1,0)$
   & $\langle y_0^4,  y_4y_5 \rangle $
    & $\langle y_0,  y_4y_5 \rangle $
     & \cr
height2pt & \omit & \omit & \omit & \omit &\cr
 & $\tau_3$
  & $(1,0,0,0,0,0,1)$
   & $\langle y_0^4,  y_6^2 \rangle $
    & $\langle y_0,  y_6 \rangle $
     & \cr
height2pt & \omit & \omit & \omit & \omit &\cr
 & $\tau_4$
  & $(1,0,0,0,0,0,5)$
   & $\langle y_1y_2y_3y_6,  y_6^2, y_0^8, y_0^4y_6 \rangle $
    & $\langle y_0,  y_6 \rangle $
     & \cr
height2pt & \omit & \omit & \omit & \omit &\cr
 & $\tau_5$
  & $(0,1,1,1,0,0,4)$
   & $\langle y_1y_2y_3y_6,y_6^2,y_1^2y_2^2y_3^2y_4y_5, y_0^4y_6 \rangle $
    & $\langle y_6,  y_1y_2y_3y_4y_5 \rangle $
     & \cr
height2pt & \omit & \omit & \omit & \omit &\cr
 & $\tau_6$
  & $(0,1,1,1,0,0,1)$
   & $\langle y_1y_2y_3y_6,y_6^2, y_1y_2y_3y_4y_5 \rangle $
    & $\langle y_6, y_1y_2y_3y_4y_5 \rangle $
     & \cr
height2pt & \omit & \omit & \omit &  \omit & \cr }
\hrule}}
\vskip0.2cm
\vbox{\baselineskip=12pt
\noindent
{\tabfont
Table 1. Gr\"obner cones with typical weights.  \par\noindent
Each cone determines the ideal $LT_\omega({\cal I}_\A)$ and its radical.
The radical coincides with the Stanley-Reisner ideal $SR_{T_\omega}$
according to \thmI. }}
\vskip0.2cm

\vskip0.5cm

\subsec{Cohomology ring of $\IP_{\SDs}$ and the local solutions
    --- when $\SDs$ is regular --- }

In this subsection we will study the local solutions of the
$\Ds$-hypergeometric system. Since the Gr\"obner
fan is a refinement of the secondary fan, each cone of the Gr\"obner
fan naturally defines a convergent series for \fseries.
Namely we consider a cone $\tau$ with typical weight $\omega$. If $\tau$
is simplicial and regular we consider a $\ZZ$-basis
$\{\tilde l_\tau^{(1)}, \cdots, \tilde l_\tau^{(p-n)}\}$ of $\tau$, and
if not we subdivide $\tau$ into simplicial and regular cones and take a
$\ZZ$-basis for one of these cones. Then the dual basis
$\{ l_\tau^{(1)}, \cdots, l_\tau^{(p-n)}\}$ gives us a $\ZZ$-basis
compatible with the regular triangulation $T_\omega$ and the series
\wnot. Even though the choice of $\ZZ$-basis of $L$ is not unique,
once a choice is made we refer to it as a $\ZZ$-basis of $L$ for
the cone $\tau$ with typical weight $\omega$.

Since the exponent $\beta$ is $T$-resonant for some
regular triangulations, we do not have vol$({\Ds(w)})$ linearly
independent power series solutions, and so we need to search for the
logarithmic solutions.
The type of logarithmic solutions which arises
for a given triangulation depends on the type of degenerations
of the hypersurface hence on the monodromy of its period integrals.
In general, the differential equations
satisfied by the period integrals have regular singularities
\ref\Griffiths{P.Griffiths, Ann. of Math.(2){\bf 80}(1964)227.}.
Therefore we can determine the local solutions from the data of the
leading terms of the differential equations -- the so-called indicial
equations. We expect that among the singularities,
the quotient singularities can be resolved
by toric method via a refinement -- such as the Gr\"obner
fan -- of the secondary fan(\BatyrevII\ Conjecture 13.2).
Thus near these singularities,
we should recover the data for our
local solutions from the structure of the cones in the fan.

Let us consider a power series solution \wnot\ determined
by a $\ZZ$-basis $\{ l_\tau^{(1)},\cdots,l_\tau^{(p-n)} \}$ of $L$
for a cone $\tau$ with typical weight $\omega\in \tau$.
We identify the toric ideal ${\cal I}_\A$ in $\IC[y_1,\cdots,y_p]$
with the ideal generated by $\{ \sqbox_l \}$ in
$\IC[{\pd \; \over \pd a_1},\cdots, {\pd\; \over \pd a_p} ]$.
While we will consider a multiplication of $\sqbox_l$ by the rational
functions of $a_k$'s extending the coefficient, we need to
be careful with the noncommutativity resulting from this extension.
Now let us consider an operator in
the Gr\"obner basis ${\cal B}_\omega$. If $\sqbox_l \in {\cal B}_\omega$
and $\omega\cdot l_+ -\omega\cdot l_- > 0 \; (<0)$ then
$({\pd \; \over \pd a})^{l_+} \;\;
(\; ( {\pd \; \over \pd a} )^{l_-} \;)$ is one of the generators for
the ideal $LT_\omega({\cal I}_\A)$. For the case
$\omega\cdot(l_+-l_-) >0$, we multiply ${\cal D}_l$ by $a^{l_+}$ to
obtain
\eqn\msqbox{
a^{l_+} \sqbox_l = a^{l_+}\( {\pd \; \over \pd a} \)^{l_+} -
a^{l_+ - l_-} a^{l_-} \( {\pd \; \over \pd a} \)^{l_-}
}
Since $\omega\cdot(l_+-l_-)>0$ we have $l_+-l_- \in \tau^\vee$.
Since $\tau^\vee$ is generated by $\{l^{(1)},..,l^{(p-n)}\}$, it follows that
$a^{l_+ -l_-}$ in the second term
can be expressed by a monomial of $\{ x_\tau^{(k)} \}$ which vanishes when
$x_\tau^{(k)} \rightarrow 0$. Other parts in \msqbox\ are
'homogeneous' and can be rewritten in terms of the log derivatives
$\theta_{a_i}=a_i {\pd \; \over \pd a_i }$. The same argument applies
to the case $\omega\cdot(l_+-l_-)<0$.
Therefore the principal part of $\sqbox_l$ which determines the local
properties about $x_\tau^{(k)}=0$ are given, through
the generators of ${\cal B}_\omega$, by
\eqn\mapind{
\(\da{}\)^{l_+}-\(\da{}\)^{l_-}
\mapsto a^{l_\pm} \( {\pd \; \over \pd a} \)^{l_\pm} \;\;,
}
where $l_+$ and $l_-$ in the right hand side for $\omega\cdot(l_+-l_-)>0$
and $\omega\cdot(l_+-l_-)<0$, respectively. Clearly we can
express the principal part
\mapind\ as a polynomial $I_l(\ta{1},\cdots,\ta{p})$.
 In the gauge $\tilde \Pi =\tilde \Pi(x_\tau^{(1)},\cdots,
x_\tau^{(p-n)})$ \tgauge, using $x_k=a^{l^{(k)}_\tau}$ we can write
\eqn\ppart{
I_l(\ta{1},\cdots,\ta{p})=J_l(\theta_{x_\tau^{(1)}},\cdots,
\theta_{x_\tau^{(p-n)}}) \;,
}
where the right hand side is a polynomial in the
$\theta_{x_\tau^{(k)}}$. Note also that $J_l$ is homogeneous
if the entries of $l_\pm$ are 0 or 1.
Due to the property of the Gr\"obner
basis, the principal parts \ppart\ for the elements in ${\cal B}_\omega$
give us a complete set of the indicial equations for $\rho$.
Summarizing our results;

\vskip0.5cm
\leftskip=0.8cm \rightskip=1.5cm
\vbox{
\noindent
{\it
Consider a local solution $w_0(x_\tau,\rho)$\wnot\ with a $\ZZ$-basis
of $L$ for a cone $\tau$ with typical weight $\omega$.
Then the indices $\rho$ are
determined from the finite set of indicial equations
$J_l(\rho_\tau^{(1)},\cdots,\rho_\tau^{(p-n)})=0$,
with $y^{l_+}-y^{l_-}\in{\cal B}_\omega$.
}

\vskip-50pt
\eqn\wnotsol{ }
\vskip25pt
}
\leftskip=0cm \rightskip=0cm
\noindent
This results combined with the Frobenius method enables us to construct
missing solutions in the general theorem\GKZth\ for the case of
$T$-resonant.
\vskip0.3cm

Now let us turn to the description of the intersection ring
$A^*(\IP_\SDs,\ZZ)$, which is isomorphic to the cohomology
ring, $H^{2*}(\IP_\SDs, \ZZ)$ of the nonsingular projective toric variety
$\IP_\SDs$. In the following we assume $\IP_\SDs$ is nonsingular, which
means that we take one of the finest subdivisions of the fan $\SDs$.
Note that for the models of type I
and II, such finest subdivision comes from a maximal triangulation
$T_0$ of the polyhedron $\Ds$. We have seen that $T_0$ is
also regular.
In the next section, we will find that our results apply
also to the singular models of type III with some modifications.

In toric geometry, each integral point $\nu_i^* \;
(i=1,\cdots,p)$ in $\Ds\cap \ZZ^n$ corresponds to an irreducible $T$-invariant
divisor $D_i$. It is known that if $\ns{i_1},\cdots,\ns{i_k}$ generate
a cone in
$\SDs$, the divisors $D_{i_1}\cdots D_{i_k}$ intersect
transversally with the subvariety determined by the cone. Also there
are linear relations among the toric divisors since we are working
modulo the divisors of rational functions on $\IP_\SDs$.
It is then known that \Oda:

\vskip0.5cm
\leftskip=0.8cm \rightskip=1.5cm
\vbox{
\noindent
{\it
For a compact nonsingular toric variety $\IP_\SDs$, the intersection
ring $A^*(\IP_\SDs, \ZZ)$ is described by $\ZZ[D_1,\cdots,D_p]/{\cal I}$,
where ${\cal I}$ is the ideal generated by
\item{(i)} $D_{i_1}\cdot \cdots \cdot D_{i_k}$ for $\nu_{i_1}^*,
\cdots,\nu_{i_k}^*$ not in a cone of $\SDs$,
\item{(ii)} $\sum_{i=1}^p \langle u, \nu_i^* \rangle D_i$ for $u\in \ZZ^n$.
}

\vskip-50pt
\eqn\intring{ }
\vskip20pt
}
\leftskip=0cm \rightskip=0cm
\noindent
We can fix the normalization of the 'volume form' in the ring by the
property that the Euler
number of $\IP_\SDs$ is equal to the number of the $n$-dimensional cones
in the fan $\SDs$. This is the number of the $n$-simplices
in the corresponding maximal triangulation $T_0$ of the polyhedron
$\Ds$.

To related toric ideals to our previous discussion on
 the $\Ds$-hypergeometric
system, we introduce a formal
variable $D_0$ and rewrite the intersection ring as
$A^*(\IP_\SDs ,\ZZ)=\ZZ[D_0,D_1,\cdots,D_p]/
\bar {\cal I}$, where  we define $\bar {\cal I}$ as the ideal generated
by

\vskip0.5cm
\leftskip=0.8cm \rightskip=1.5cm
\vbox{
\noindent
{\it
\item{$(i)'$} $D_{i_1}\cdots D_{i_k}$ for $\bn{i_1}, \cdots,\bn{i_k}$ not
in a cone of $\Sigma((1,\Ds))$,
\item{$(ii)'$} $\sum_{i=0}^p \langle u,\bn{i}\rangle D_i$ for $u \in
\ZZ^{n+1}$.
}
\vskip-40pt
\eqn\intringext{ }
\vskip20pt
}
\leftskip=0cm \rightskip=0cm
\noindent
The fan $\Sigma((1,\Ds))$ in $(i)'$ is defined to be a set of cones over
the simplices of the triangulations $T_0$ of $(1,\Ds)$.
 If the fan $\SDs$ is regular,
 so is the fan $\Sigma((1,\Ds))$, although the latter fan is not complete.

Now note that the set of the generators $(i)'$ is the same as the generators
of the Stanley-Reisner ideal for the maximal triangulation $T_0$ of $\Ds$.
Note also the similarity of the linear relations
$(ii)'$ and the first order relations \PSdiff\ in the hypergeometric
system. By \thmI\ the Stanley-Reisner ideal $SR_{T_0}$ is
the radical of $LT_\omega({\cal I}_\A)$, where $\omega$ is a
weight with $T_\omega=T_0$. In the following, we will show that
the ideal $LT_\omega({\cal I}_\A)$ is radical. This allows us to
determine  the
ideal $LT_\omega({\cal I}_\A)$, or equivalently
 the principal parts \ppart\ of the
$\Ds$-hypergeometric system that governs the local solutions, via a purely
combinatorial object -- the Stanley-Reisner
ideal. As an immediate consequence we show that for the maximal
triangulation, $\rho=(0,..,0)$
is the unique solution to the indicial equations.
Furthermore we observe that the
local solutions for the maximal triangulation $T_0$ can be
described by the
intersection numbers.
The latter are computable from the intersections ring \intring.
We will also see that the Stanley-Reisner ideal above can be
easily computed in terms of the so-called primitive collections.

To discuss the combinatorial description of the Stanley-Reisner ideal, we
introduce the notions of
{\it a primitive collection} and {\it a primitive relation}
\BatyrevP.
A primitive collection of a complete fan $\SDs$ is a set of
integral vectors ${\cal P}=\{ \ns{i_1}, \ns{i_2}, \cdots, \ns{i_a} \}$
such that if we remove any one of $\ns{i_s}$ from ${\cal P}$, then the
integral vectors in ${\cal P}\setminus \{\ns{i_s}\}$ generate a cone in $\SDs$
while ${\cal P}$ itself does not generate any cone in $\SDs$. It is easy to
prove that $i)$ in \intring\ can be replaced by the monomials
$D_{i_1}\cdots D_{i_a}$ corresponding to the
primitive collections of $\SDs$. Once we fix a triangulation $T_0$ which
underlies the $\SDs$, it is straightforward to read off all primitive
collections. So far we don't need the regularity of $\SDs$.
But to discuss primitive relations, we must assume that $\SDs$ is
regular.
A primitive relation will be a certain element of $L$
attached to each primitive collection.

Consider a primitive collection ${\cal P}=\{
\ns{i_1}, \ns{i_2}, \cdots, \ns{i_a} \}$. There is a unique
cone ${\cal C}\in\SDs$ of minimum dimension such that the integral point
$\ns{i_1}+\ns{i_2}+\cdots+\ns{i_a}$ is in the interior of ${\cal C}$.
By regularity, there is a set of generators $\{\ns{j_1},..,\ns{j_s}\}$
 of the cone
such that for some positive integers $c_k$, we have
\eqn\pcolI{
\ns{i_1}+\ns{i_2}+\cdots+\ns{i_a}=\sum_{k\geq 1} c_k \ns{j_k} \quad.
}
It is easy to translate the above statement about the
fan $\SDs$ into a statement about the fan  $\Sigma((1,\Ds))$,
which is not complete but regular. We get
\eqn\pcolII{
\bn{i_1}+\bn{i_2}+\cdots+\bn{i_a}=\sum_{k \geq 0} c_k \bn{j_k} \quad,
}
where $\bn{j_0}=\bn{0}$ and
$c_0=a-\sum_{k\geq 1} c_k  \; \geq 0$.
Eqn. \pcolII\ defines a primitive relation $\l({\cal P}) \in L$. It is easy
to deduce from the defining property of a primitive collection that
the index sets $\{i_1,..,i_a\}$ and $\{j_0,..,j_s\}$ are disjoint.

Let $\omega$ be a weight vector such that
$T_\omega$ is the maximal triangulation $T_0$.
Recall that the convex polytope $P_\omega$ is
defined by the convex hull of the points $\tilde \ns{k}:=(\omega_k,
\ns{k} ) \; (k=0,\cdots,p)$. Then we can show the following ``height''
inequality
$(\tilde\ns{i_1} +\cdots +\tilde\ns{i_a})_0>(\sum c_k \tilde\ns{j_k})_0$,
i.e.,
\eqn\weight{
\omega_{i_1}+\omega_{i_2}+\cdots +\omega_{i_a} > \sum c_k \omega_{j_k} \;\;.
}
This means that $LT_\omega(y^{l({\cal P})_+}-y^{l({\cal P})_-}) =
y_{i_1}y_{i_2} \cdots y_{i_a}$. Since the Stanley-Reisner ideal
 $SR_{T_\omega}$ is generated by those $y^{i_1}y^{i_2} \cdots y^{i_a}$
with $\{\ns{i_1}, \ns{i_2}, \cdots, \ns{i_a} \}$ primitive,
it follows that $SR_{T_\omega}\subset LT_\omega({\cal I}_\A)$.
Combining this with the property \thmI, we see
that $LT_\omega({\cal I}_\A)$ is radical. Moreover we also have
$SR_{T_\omega}=\langle y^{l({\cal P})_+}\vert \; {\cal P}\;{\rm is}\;{\rm
primitive} \rangle$.

In the example $\IP(2,2,2,1,1)$
discussed in the last subsection, we find two primitive collections
$\{ \ns{1},\ns{2},\ns{3},\ns{6} \}, \{\ns{4},\ns{5} \}$ for the maximal
triangulation $T_0$ and the corresponding primitive relations turn out
to be $l^{(1)}$ and $l^{(2)}$ in \exIIels, respectively. As is evident
in the table 1, these primitive collections give the
generators $LT_\omega(y^{l^{(i)}_+}-y^{l^{(i)}_-})$, ($i=1,2$)
of  the ideal $LT_\omega({\cal I_{\cal A}})$ for the cone $\tau_1$.

Note that
the above generators of $SR_{T_0}$ are nothing but
the leading symbols of a generating set of operators ${\cal D}_l$
of the $\Ds$-hypergeometric system. Combining with the
operators $\tilde{\cal Z}_i$, we have a correspondence between
the symbols of the full $\Ds$-hypergeometric system and the ideal $\bar{\cal
I}$ \intringext\ for the intersection ring. This motivates the following map
$m$ from $\ZZ[\ta{0},..,\ta{p}]$ to the
intersection ring $m: \ta{i} \mapsto D_i$.
 Define the following intersection coupling,
\eqn\coupling{
C_{i_1i_2\cdots i_{n}}=\langle
m(\theta_{x_\tau^{(i_1)}})
m(\theta_{x_\tau^{(i_2)}})
\cdots
m(\theta_{x_\tau^{(i_n)}})  \rangle  \;\; ,
}
where the bracket means taking the coefficient of the 'volume element'
in the ring $A^*(\IP_\SDs,\ZZ)$.
Then we {\it observe} the following;

\vskip0.5cm
\leftskip=0.8cm \rightskip=1.5cm
\vbox{
\noindent
{\it
If $\tau$ is the cone in which
$T_\omega \;$ ($\omega\in\tau$) is a maximal triangulation,
then all the indices at the point $x_\tau^{(i)}=0 \; (i=1,\cdots,p-n)$
of the hypergeometric system are identically zero. And the local solutions
near this point are given by
$$
\eqalign{
&
w_0(x_\tau,\rho)\vert_{\rho=0} \;,\;
\drho{i} w_0(x_\tau,\rho)\vert_{\rho=0} \;,\;
\sum_{i_1,i_2} C_{i_1 i_2 \cdots i_{n}} \drho{i_1}\drho{i_2}
w_0(x_\tau,\rho)\vert_{\rho=0} \;,\;  \cr
& \hskip1cm \cdots \cr
&
\sum_{i_1,i_2, \cdots i_{n}} C_{i_1 i_2 \cdots i_{n}}
\drho{i_1}\drho{i_2}
\cdots\drho{i_{n}}
w_0(x_\tau,\rho)\vert_{\rho=0} \;.\;  \cr}
$$
}

\vskip-120pt
\eqn\gsol{ }
\vskip80pt
}
\leftskip=0cm \rightskip=0cm
\noindent
Recall that  rad$(LT_\omega({\cal I}_\A)) =LT_\omega({\cal I}_\A)$ for
the weight $\omega$ such that $T_\omega=T_0$, and that
$LT_\omega({\cal I}_\A)$ is generated by $LT_\omega({\cal B}_\omega)$.
Since an element of the Gr\"obner basis has the form
 $y^{l_+}-y^{l_-}$, it follows that the entries of either vector
$l_\pm$ are 0 or 1. In either case, we see that the corresponding
indicial equation $J_l(\rho_1,..,\rho_{p-n})=0$ is homogeneous.
But the finiteness of the solution set implies that zero is the only
solution.

\vskip50pt

\leftline{\undertext{Example: $\IP(2,2,2,1,1)$}}

As we have seen, there is one maximal triangulation $T_0$ in
\exIItriang\ for the polyhedron $\Ds$. The corresponding cone is
$\tau_1$ in table 1, and thus the Gr\"obner basis ${\cal B}_\omega$
consists of $y^{l_+^{(1)}}-y^{l_-^{(1)}}=y_1y_2y_3y_6-y_0^4$
and $y^{l_+^{(2)}}-y^{l_-^{(2)}}=y_4y_5-y_6^2$.
{}From this we obtain the leading term operators for \ppart;
\eqn\tauI{
\eqalign{
J_{l^{(1)}} &= \ta{1}\ta{2}\ta{3}\ta{6}=\theta_{x_{\tau_1}}^3
(\theta_{x_{\tau_1}}
-2\theta_{y_{\tau_1}})  \;\;,  \cr
J_{l^{(2)}}  &= \;\; \ta{4}\ta{5} \;\;  = \theta_{y_{\tau_1}}^2 \;\;, \cr}
}
with the cprresponding linear operators
\eqn\diff{
\eqalign{
{\cal L}_1&=\theta_{x_{\tau_1}}^3 (\theta_{x_{\tau_1}}-2\theta_{y_{\tau_1}})
-x_{\tau_1}(4\theta_{x_{\tau_1}}+4)(4\theta_{x_{\tau_1}}+3)
           (4\theta_{x_{\tau_1}}+2)(4\theta_{x_{\tau_1}}+1)  \cr
{\cal L}_2&=\theta_{y_{\tau_1}}^2-y_{\tau_1}
(\theta_{x_{\tau_1}}-2\theta_{y_{\tau_1}}-1)
(\theta_{x_{\tau_1}}-2\theta_{y_{\tau_1}}) \; .\cr}
}
Since the generators of the Stanley-Reisner ideal is given
by the primitive collections, it is easy to
determine the intersection ring. The results for the intersection
couplings are
\eqn\exIICxx{
C_{xxxx}=2 \;\;,\;\; C_{xxxy}=1 \;\;,\;\;
C_{xxyy}=C_{xyyy}=C_{yyyy}=0   \;\; ,
}
where $C_{xxxx}=\langle m(\theta_{x_{\tau_1}}) \cdots
m(\theta_{x_{\tau_1}}) \rangle $ for example. From the
indicial equations $J_{l^{(1)}}(\rho)=J_{l^{(2)}}(\rho)=0$,
we see that all indices at the point $x_{\tau_1}=y_{\tau_1}=0$
are zero. In fact we find the following
8 solutions with only one power series solution;
\eqn\exIIsol{
\eqalign{
&
w_0(x,0) \;;\;
\drho{x} w_0(x,0) \;,\;
\drho{y} w_0(x,0) \;;\;
(2\drho{x}^2+2\drho{x}\drho{y}) w_0(x,0) \;,\;
\drho{x}^2 w_0(x,0) \;; \; \cr
&
(2\drho{x}^3+3\drho{x}^2\drho{y}) w_0(x,0) \;,\;
\drho{x}^3 w_0(x,0) \;; \;
(2\drho{x}^4+4\drho{x}^3\drho{y}) w_0(x,0) \;,\; \cr}
}
where
\eqn\www{
w_0(x,\rho)=\sum {\Gamma(4(n+\rho_x)+1) \over
\Gamma(n+\rho_x+1)^3 \Gamma(m+\rho_y+1)^2 \Gamma(n-2m+\rho_x-2\rho_y+1) }
x_{\tau_1}^{n+\rho_x} y_{\tau_1}^{m+\rho_y} \;\;.
}
Because we have ${\cal L}_i w_0(x,\rho)=J_{l^{(i)}}(\rho) x^\rho, \; (i=1,2)$,
it is easy to verify that \exIIsol\ solve the hypergeometric system by
inspecting $(2\drho{x}^4+4\drho{x}^3\drho{y})
(J_{l^{(i)}}(\rho) x^\rho) \vert_{\rho=0}=0 \; (i=1,2)$, for example.

Similarly we obtain for $\tau_2$ the principal parts\ppart\ in the $\tilde
\Pi(a)$
gauge
\eqn\tauII{
(4\theta_{x_{\tau_2}}-4)(4\theta_{x_{\tau_2}}-3)
(4\theta_{x_{\tau_2}}-2)(4\theta_{x_{\tau_2}}-1) \;\; ,\;\;
\theta_{y_{\tau_2}}^2 \;,
}
with $x_{\tau_2}=a^{-l^{(1)}}=1/x_{\tau_1}\;,\; y_{\tau_2}=a^{l^{(2)}}$.
{}From these principal parts, we see that not all solutions to
the indicial equations at
$x_{\tau_2}= y_{\tau_2}=0$ are zero.
Thus the local properties of
$\tau_1$ and $\tau_2$ are quite different. The fact that
$LT_\omega({\cal I}_\A)$ is radical in $\tau_1$ but not in
$\tau_2$ is responsible for this difference.

\vskip0.5cm

\subsec{ Cohomology ring of $X_\D$ and the local solutions }

We consider the restriction map $H^*(\IP_\SDs,\ZZ)\rightarrow H^*(X_\D,\ZZ)$
induced by the inclusion $X_\D\rightarrow \IP_\SDs$, and denote the image by
$H^*_{toric}(X_\D,\ZZ)$.
The restriction map can be realized by
considering the intersection of the elements of $H^*(\IP_\SDs,\ZZ)$
with the divisor $X_\D$.
By construction of the Calabi-Yau hypersurface
$X_\D$, the divisor class $[X_\D]$ coincides with the anti-canonical
class of the ambient space $\IP_\SDs$, namely
\eqn\Hdiv{
[X_\D]= D_1+D_2+\cdots +D_p \;\; ,
}
in the intersection ring. The toric part of the cohomology
$H^*_{toric}(X_\D,\ZZ)$ can then be written as
$A^*_{toric}(X_\D,\ZZ)=A^*(\IP_\SDs,\ZZ) / Ann(D_1+\cdots + D_p)$
(where $Ann(x)$ in a ring $R$ is $Ann(x):=\{y\in R|yx=0\}$)
or equivalently
\eqn\Atoric{
A^*_{toric}(X_\D,\ZZ)=\ZZ[D_1,D_2,\cdots,D_p]/ {\cal I}_{quot} \;\;,
}
where ${\cal I}_{quot}$ is the ideal quotient ${\cal I}_{quot}=
{\cal I}:(D_1+\cdots +D_p)$. (Here $(I:x)=\{y\in R|yx\in I\}$.)

Now recall the close relationship between the ideal $\bar {\cal I}$ in
\intringext\ and the ideal of the symbols for the $\Ds$-hypergeometric
system with respect to the cone $\tau$ of maximal triangulations
$T_\omega$.
First we have $\bar {\cal I}_{quot}= \bar {\cal I}:
(D_1+D_2+\cdots + D_p) =\bar {\cal I}:(-D_0)=\bar {\cal I}:D_0$.
In fact we observe more: suitable  linear combinations of
the differential operators
$a^{l\pm}\sqbox_l$
factorize from the left by the operator $\ta{0}$, implying  that
the hypergeometric system is a reducible system. Factorization
by the operator $\ta{0}$ should be understood as corresponding
to the restriction to
the hypersurface $X_\D$. As we shall see, the
solutions to the factorized system can be
obtained from \gsol\ by a similar restriction of
the intersection couplings (cf. $m(-\ta{0})=-D_0=[X_\D]$).

In cases of type II models, we observed that the quotient \Atoric\
results in setting to zero the
divisors $D_i$ for which the integral points $\ns{i}$ are on a codimension
one face of $\Ds$. This can be understood as follows:
the above divisors come from the
desingularizations of point singularities of the ambient space;
a hypersurface $X_\D$ in general position will not meet these singularities.
 In accordance with this
'decoupling' of the divisors it is natural to
consider the lattice $L'=\{l\in L
\vert l_i=0 $, $\ns{i}$ is on a codimension one face of $\Ds\, \}$.
We define a $\ZZ$-basis $\{l_\tau^{(1)},\cdots,l_\tau^{(p-n')}\} \,
(n'={\rm dim} Auto(\IP_\SDs) )$ of $L'$ of the reduced cone
$(\IR_{\geq0}l_\tau^{(1)} + \cdots + \IR_{\geq0}l_\tau^{(p-n)}) \cap
L'_{\IR}$ as follows. We make subdivisions of the reduced cone if it is not
regular,
in which case the $\ZZ$-basis is not uniquely
determined. However our observations in the following do not depend
on this.  We will call $L'$ the {\it reduced lattice}, or the reduction
of $L$.
%

The decoupling of some divisor $D_i$ in the intersection ring implies
that we can turn off the monomial deformation via $a_i$ ( which
corresponds to the divisor $D_i$ under the monomial-divisor map\AGM).
In fact we observe that
these variables can be eliminated
in the extended $\Ds$-hypergeometric system which is originally
defined to act on functions on $\IC^\A$ as follow.
Recall that the GKZ $\Ds$-hypergeometric system is enlarged by
adjoining $n'-n$ additional
linear differential operators ${\cal Z}_i \; (i=n+1,\cdots, {\rm
dim}Auto(\IP_\SD)=n'$ in \exIzxis ). This creates just enough equations to
eliminate those operators $\da{i}$ corresponding to
points $\ns{i}$ on the codimension one faces, from
the operators $\sqbox_l \; (l\in L)$. We may then
 set $a_i=0$ after the elimination.
We denote  by $\sqbox_l'$ the resulting new operators which act on
functions on $\IC^{\A'}$, where the
set $\A'$ consists of all integral points on the faces with codimension
greater than one. Note that the set
$\{ \sqbox_l' \vert l\in L \}$ is in general larger than
the set $\{ \sqbox_{l'} \vert l' \in L' \}$.

We now define the intersection couplings
on $A^*_{toric}(X_\D,\ZZ)$ by
\eqn\yukcp{
K^{cl}_{i_1i_2\cdots i_{n-1}} = \langle
m(\theta_{x_\tau^{(i_1)}}) m(\theta_{x_\tau^{(i_2)}}) \cdots
m(\theta_{x_\tau^{(i_{n-1})}})\cdot
m(-\ta{0}) \rangle \;\;,
}
then we may state the observation given in \HKTYII as follows:

\vskip0.5cm
\leftskip=0.8cm \rightskip=1.5cm
\vbox{
\noindent
{\it
For a cone $\tau$ with typical weight $\omega$,
some of the operators $a^{l_\pm} \sqbox'_l \; (l\in L)$ or their linear
combinations factorize by
the operator $\ta{0}$ from the left, indicating that the $\Ds$-hypergeometric
system is reducible. If $T_\omega$ is a maximal triangulation,
the local solutions about the point $x_\tau^{(i)}=0 \;
(i=1,\cdots,p-n')$ for the reduced system are given by
$$
\eqalign{
&
w_0(x_\tau,\rho)\vert_{\rho=0} \;,\;
\drho{i} w_0(x_\tau,\rho)\vert_{\rho=0} \;,\;
\sum_{i_1,i_2} K^{cl}_{i_1 i_2 \cdots i_{n-1}} \drho{i_1}\drho{i_2}
w_0(x_\tau,\rho)\vert_{\rho=0} \;,\;  \cr
& \hskip1cm \cdots \cr
&
\sum_{i_1,i_2, \cdots i_{n-1}} K^{cl}_{i_1 i_2 \cdots i_{n-1}}
\drho{i_1}\drho{i_2}
\cdots\drho{i_{n-1}}
w_0(x_\tau,\rho)\vert_{\rho=0} \;.\;  \cr}
$$
}
\vskip-120pt
\eqn\Ksol{ }
\vskip80pt
}
\leftskip=0cm \rightskip=0cm

\noindent

We also observe that in the case of Fermat hypersurfaces, the
operators $a^{l_\pm}\sqbox'_l$ which factorize
whose leading term generate the ideal $\bar {\cal I}_{quot}$
can be constructed from operators $\sqbox_l$ in the Gr\"obner basis
${\cal B}_\omega$. However for general
models of non-Fermat type we need to
consider operators $a^{l_\pm}\sqbox'_l$ outside the basis
${\cal B}_\omega$ as well as their linear combinations with coefficients
in the ring generated by the $\ta{i}$ (see examples in sect. 4).

\vskip0.3cm
\leftline{\undertext{Examples:} }

\noindent
1) $\IP(2,2,2,1,1)$:
We have seen a unique maximal triangulation $T_0$ in \exIItriang\ and
have constructed local solutions for the corresponding cone
in the Gr\"obner fan. Now we
note that $\ta{0}=-(\ta{1}+\cdots +\ta{6})=-4\theta_{x_{\tau_1}}$. It is
easy to observe that the operator $a^{l_+^{(1)}} \sqbox_{l^{(1)}}$
expressed in the $\tilde \Pi(a)$ gauge factorizes from the
left by $\ta{0}=-4\theta_{x_{\tau_1}}$, i.e.,
${\cal L}_1=\theta_{x_{\tau_1}}{\cal O}$ in \diff\ for some
third order operator ${\cal O}$.
If we write the the divisor $m(\theta_{x_{\tau_1}})$ as $J_x$ and similarly
for $J_y$, the topological data for the solutions are summarized as follows:
\eqn\exIItop{
\eqalign{
 K^{cl}_{xxx}=8 & \;\;,  \;\; K^{cl}_{xxy}=4 \;\;, \;\;
K^{cl}_{xyy}=K^{cl}_{yyy}=0 \cr
& c_2\cdot J_x=56 \;\;, \;\; c_2\cdot J_y=24 \;\;,\cr}
}
where the invariants $c_2\cdot J$'s are listed for later use. For their
calculation we use the adjunction formula\GH
;$c(X_\D)=\prod_{i=1}^{p}(1+D_i)/(1+[X_\D])$.

\noindent
2) $\IP(7,2,2,2,1)$:
The toric data of this model have been summarized in the end of the
previous section. Although this model has the same moduli as the above
model, two integral points on the codimension one face make the
combinatorics of this model much more complicated.
It turns out that there are 14 elementary relations which generate
the zonotope ${\cal P}_\A$, and there are 2,154 elements for the universal
Gr\"obner basis.
The secondary fan has 32 four dimensional cones, most of
which are singular.

It seems to be a formidable task to determine the Gr\"obner fan,
however it is easy to find the maximal triangulation of $\Ds$
and the corresponding Stanley-Reisner ideal. As proved in the previous
subsection, for a weight $\omega$ such that $T_\omega$ is the maximal
triangulation, we have $LT_\omega({\cal I}_\A)
= {\rm rad}(LT_\omega({\cal I}_\A)) = SR_{T_\omega}$ and
we can determine the ideal $LT_\omega({\cal I}_\A)$ by
the Stanley-Reisner ideal which is simply described
 by the primitive collections.
In this case, it turns out that the ideal $LT_\omega({\cal I}_\A)=
SR_{T_\omega}$ is generated by
\eqn\exIsr{
\eqalign{
y_1y_5 & \;\;,
\;\; y_1y_7 \;\;,\;\; y_1y_8 \;\;,\;\;
y_5y_8\;\;,\;\; y_6y_7 \;\;, y_6y_8 \;\;,  \cr
&
y_2y_3y_4y_5 \;\;,\;\; y_2y_3y_4y_6 \;\;,\;\; y_2y_3y_4y_7 \;\;. \cr}
}
These generators may be translated into the generators $(i)'$ in
\intringext.
Then together with the linear relations $(ii)'$ in \intringext,
they define the intersection ring of the ambient space.
We find that the ideal $\bar {\cal I}_{quot}=\bar {\cal I}:(D_1+
\cdots+D_8)$ as defined earlier is generated by the monomials
\eqn\exIquot{
D_{1} D_{5} \;\;,\;\;
D_{3} D_{4} D_{6} \;\;,\;\;
D_{7} \;\;, \;\; D_{8} \;\;,
}
with the linear relations $(ii)'$.
The divisors
$D_{7}$ and $D_{8}$ in \exIquot\ being among the generators show that
these divisors decouple from the intersection ring.

We find that a $\ZZ$-basis of $L$ for a cone $\tau$ with the typical weight
$\omega$ is $\{l^{(1)}_\tau, \l^{(2)}_\tau, l^{(3)}_\tau, l^{(4)}_\tau\}$ with
\eqn\exIlbasis{
\eqalign{
&
l^{(1)}_\tau =( -1,0,0,0,0, -1,1,1,\;0) \;,\;
l^{(2)}_\tau =(0,1,0,0,0,1,-2,\;0,0) \;,\;  \cr
&
l^{(3)}_\tau =(\;0,0,1,1,1,\;0,0,1, -4) \;,\;
l^{(4)}_\tau =(0,0,0,0,0,1,\;0,-2,1) \;.\;  \cr }
}
The intersection of the cone generated by \exIlbasis\ with $L_\IR'$ is a
two dimensional regular cone generated by
$ l^{(1)}=7l^{(1)}_\tau+l^{(3)}_\tau+4l^{(4)}_\tau =
(-7,0,1,1,1,-3,7,0,0)$ and $ l^{(2)}=l^{(2)}_\tau $.

The form of our generators of $\bar{\cal I}_{quot}$ suggests
that we should try to factorize
the following differential operators
\eqn\exIbeffac{
\eqalign{
\sqbox_{l_{\{1,5\}}}&=\da{1}\da{5}-(\da{6})^2 \;\;,\cr
\sqbox_{l_{\{2,3,4,6\}}}&=
\da{2}\da{3}\da{4}\da{6}-\da{0}(\da{8})^3 \;\;,\cr }
}
where $l_{\{1,5\}}$ and $l_{\{2,3,4,6\}}$ are the primitive
relations corresponding respectively to the primitive
collections $\{\ns{1},\ns{5}\}$ and $\{\ns{2},\ns{3},\ns{4},\ns{6}\}$.
Although the second operator contains a derivative with respect to
$a_8$, we can eliminate it using the order one operators
${\cal Z}_{\xi_1}$ and ${\cal Z}_{\xi_2}$ corresponding
to the automorphisms \exIzxis.
Defining the local variables $x=-a^{l^{(1)}} \;\; y=a^{l^{(2)}}$,
we observe the factorization of the operator $\ta{0}$
in $a_0a_2a_3a_4a_6 \sqbox_{l_{\{2,3,4,6\}}} {1 \over a_0}$,
and find a complete set of differential equations for
the period integrals:
\eqn\exIdiff{
\eqalign{
&{\cal D}_1=(\ty-3\tx)\ty - y (7\tx-2\ty-1)(7\tx-2\ty) \;,\cr
&{\cal D}_2=
\tx^2(7\tx-2\ty)-7x\(y(28\tx-4\ty+18)+\ty-3\tx-2)\) \cr
\quad & \times
\(y(28\tx-4\ty+10)+\ty-3\tx-1\)
\(y(28\tx-4\ty+2)+\ty-3\tx\), \cr }
}
in the $\tilde \Pi(a)$ gauge \tgauge.
The local solutions of this system are given by \Ksol\ with
the following topological data;
\eqn\exItop{
\eqalign{
K^{cl}_{xxx}=2 &\;\;, \;\; K^{cl}_{xxy}=7 \;\;, \;\; K^{cl}_{xyy}=21
\;\;,\;\; K^{cl}_{yyy}=63 \cr
& c_2\cdot J_x=44 \;\;, \;\; c_2\cdot J_y=126 \;\;.\cr}
}
\vskip0.5cm

\subsec{ Singular models of type III }

In the previous subsections, we have considered the non-singular
models, i.e., models of type I and II in our classification \models.
However in actual applications, singular models
dominate the others. We will see, nevertheless, that several properties
observed in the previous subsections apply with some modifications
even to the singular cases.

Since a complete analysis of the secondary (Gr\"obner)
fan for $\Ds$ is formidable in general
(cf. the example $\IP(7,2,2,2,1)$), we focus only on the
Calabi-Yau phase(s) which corresponds to maximal triangulation(s).
For the nonsingular models of type I and II, we have seen that the ideal
$LT_\omega({\cal I}_\A)$ for a maximal triangulation $T_\omega$ coincides
with the Stanley-Reisner ideal.  For the singular models
of type III however, the ideal $LT_\omega({\cal I}_\A)$ differ from
its radical and from the Stanley-Reisner ideal because
$\SDs$ is no longer regular.


For a singular model, the fan $\SDs$ is singular even relative to the
maximum subdivision incorporating all integral points in $\Ds$. To obtain
a regular fan, which we denote as $\SDs_{reg}$, we subdivide
further the singular cones
taking into account integral points {\it outside} the polyhedron $\Ds$.
Since the polyhedron $\Ds$ is reflexive, the integral points which
generate an $n$-dimensional cone in $\SDs$ are on a hyperplane with
integral distance one from the origin. Moving this hyperplane in
a parallel way to the integral points outside $\Ds$, we can speak of the
integral distance of these points. For the hypersurfaces $X_d(w)$ in \pot,
a point with the integral distance $k>0$ corresponds to a monomial of
the homogeneous degree $kd$.

Let us denote all the integral points generating the one
dimensional cones of $\SDs_{reg}$ as $\{ \ns{1},\cdots,\ns{p},\ns{p+1},
\cdots,\ns{q} \}$ where $\ns{p+1},\cdots,\ns{q}$  are
those new points introduced by the subdivision.  (Note that even though
the new points have distance greater than 1, they are still primitive
vectors of the lattice.) Since we have a nonsingular toric variety
$\IP_{\SDs_{reg}}$, we can describe its
intersection ring according to \intring\ with additional divisors
$D_{p+1},\cdots,D_q$. It turns out
that the divisor class of the Calabi-Yau hypersurface $X_\D$ in this
fully resolved ambient space is given by
\eqn\ndiv{
[X_\D]=D_1+\cdots+D_p+d_{p+1}D_{p+1}+\cdots+d_q D_q \;\; ,
}
where $d_k$ is the integral distance of the point $\ns{k}$
described above.
We should note that the regular fan $\SDs_{reg}$ need not be the fan
associated with a triangulation of the polyhedron $\Dps:= {\rm Conv.}(\{
\ns{1},\cdots,\ns{p}, \ns{p+1},\cdots,\ns{q} \})$.
Therefore  in general, we do not have
a description \intring\ of the intersection ring via
the Stanley-Reisner ideal in terms of a triangulation of $\Dps$.
However in many cases, it happens that the convex hull $\Dps$
is itself a reflexive polyhedron.
In such a case we have another family of Calabi-Yau manifolds $X_\Dp$
in the ambient space $\IP_{\SDps}$. This ambient space is
in general different from $\IP_{\SDs_{reg}}$.
 However if we have the relation
$\SDps=\SDs_{reg}$, then we will have two different families of Calabi-Yau
hypersurfaces in the same ambient space $\IP_{\SDps}=\IP_{\SDs_{reg}}$.
One hypersurface $X_\D$ represents the  divisor class
\ndiv\ and the
other hypersurface $X_\Dp$ represents
\eqn\mdiv{
[X_\Dp]=D_1+\cdots+D_p+D_{p+1}+\cdots+D_q \;\;.
}
We will see an example of this type in sect.4.

Now let us see the detailed analysis in a typical example $X_{12}(4,3,2,2,1)$
which was analyzed in \HKTYI. The polyhedron $\D(w)$ for this model has
vertices
\eqn\singmus{
\eqalign{
&
\nu_1=(2,-1,-1,-1) \;,\;
\nu_2=(-1,3,-1,-1) \;,\;
\nu_3=(-1,-1,5,-1) \;,\;  \cr
&
\nu_4=(-1,-1,-1,5) \;,\;
\nu_2=(-1,-1,-1,-1) \;,\; \cr}
}
with respect to the basis $\Lambda_1,\cdots,\Lambda_4$ for the
lattice $H(w)$ as in \exImus. The  integral points in the dual
polyhedron $\Ds(w)$ are
as
\eqn\singnus{
\eqalign{
&
\ns{1}=(1,0,0,0) \;,\;
\ns{2}=(0,1,0,0) \;,\;
\ns{3}=(0,0,1,0) \;,\;
\ns{4}=(0,0,0,1) \;,\;  \cr
&
\ns{5}=(-4,-3,-2,-2) \;,\;
\ns{6}=(-2,-1,-1,-1) \;,\;  \cr }
}
together with the origin $\ns{0}=(0,0,0,0)$. The maximal triangulation
of the polyhedron $\Ds(w)$ is unique and is given by
\eqn\singT{
\eqalign{
T_0
=\{ &
\langle 0,3,4,5,6 \rangle,
\langle 0,1,3,4,5 \rangle,
\langle 0,2,3,4,6 \rangle,
\langle 0,1,4,5,6 \rangle, \cr
&
\langle 0,1,3,5,6 \rangle,
\langle 0,1,2,4,6 \rangle,
\langle 0,1,2,3,6 \rangle,
\langle 0,1,2,3,4 \rangle  \}  \;. \cr }
}
It is easy to see that  the corresponding fan $\SDs$ is not regular because
the first three simplices in $T_0$
respectively have volumes 2,3, and 2.
We subdivide the first cone by introducing a point
$\ns{7}={\ns{3}+\ns{4} \over2}+{\ns{4}+3\ns{5} \over 4}$.
Similarly by introducing
$\ns{8}={2\ns{1}+\ns{4} \over 3} + {\ns{3}+2\ns{5}\over 3}$ and
$\ns{9}={\ns{1}+2\ns{4} \over 3} + {2\ns{3}+\ns{5}\over 3}$ for the
second cone and
$\ns{10}={\ns{3}+\ns{4} \over 2} + {3\ns{3}+\ns{5}\over 4}$ for
the third cone, we finally obtain the regular fan $\SDs_{reg}$.
All these additional points $\ns{7},\cdots,\ns{10}$
have  integral distance two and correspond
to the charge two monomials $z_2z_3^3z_4^3z_5^9, z_1^2z_3^2z_4^2z_5^8,
z_1z_3^4z_4^4z_5^4$ and $z_2^3z_3^3z_4^3z_5^3$, respectively.
The generators $(i)$ in \intring\ are determined by the
primitive collections for the fan $\SDs_{reg}$,
and there are 20 such generators.
Together with the linear generators $(ii)$ in \intring, these
determine the defining ideal ${\cal I}$ for
the intersection ring $A^*(
\IP_{\SDs_{reg}})$. The ideal quotient by $[X_\D]=D_1+\cdots+D_6+
2(D_7+\cdots+D_{10})$ determines $A_{toric}^*(X_\D)$. It turns out that
${\cal I}_{quot}$ is generated by

\eqn\singIquot{
D_2D_5\;,\; D_1D_3D_6\;,\; D_7\;,\; D_8\;,\; D_9\;,\; D_{10} \;,
}
together with the linear relations $(ii)$. The generators $D_7,D_8,D_9$ and
$D_{10}$ indicate that these divisors decouple from the
intersection ring. This can be understood as follows: the additional points
$\ns{7},\cdots,\ns{10}$ represent point singularities in the ambient
space and the divisors $D_7,\cdots, D_{10}$ resulting from the
desingularization of these points do not intersect with the
hypersurface $X_\D$ in general position.

Now let us turn to the set of the convex piecewise linear functions over
the fan $\SDs_{reg}$, i.e., the K\"ahler cone of $\IP_{\SDs_{reg}}$
(see \OdaPark).
Since the regular fan $\SDs_{reg}$ does not come from any triangulation
of the polyhedron $\Dps$ (In fact we verify $\SDs_{reg}$ has 21 four
dimensional regular cones whereas vol($\Dps)=24$.), the K\"ahler cone so
obtained cannot be interpreted as a cone of the secondary fan for $\Dps$.
It is straightforward to find a
$\ZZ$-basis for the dual cone of K\"ahler cone $\tau$ and we have
\eqn\singIbasetau{
\eqalign{
&
l^{(1)}_\tau=(-1,0,0,1,1,0,0,\,0,1,-2,\,0) \;,\;
l^{(2)}_\tau=(-1,1,0,0,\,0,\,1,\,0,0,-2,1,0) \;,\;  \cr
&
l^{(3)}_\tau=(-2,0,0,1,1,1,1,-2,0,0,0) \;,\;
l^{(4)}_\tau=(\,1,0,0,-1,-1,-1,0,1,\,0,1,0) \;,\;  \cr
&
l^{(5)}_\tau=(-2,0,1,1,1,0,1,\,0,0,0,-2) \;,\;
l^{(6)}_\tau=(\,2,0,0,-1,-1,0,-2,1,\,0,0,1) \;.\;  \cr}
}

The decoupling of the divisors $D_7,\cdots,D_{10}$ in \singIquot\
corresponds to reducing from $L$ to the
the lattice $L'$ generated by
$l^{(1)}=4l^{(1)}_\tau+2l^{(2)}_\tau+3l^{(3)}_\tau+3l^{(4)}_\tau$ and
$l^{(2)}=l^{(3)}_\tau+l^{(5)}_\tau+2l^{(6)}_\tau$ with
\eqn\singIbase{
l^{(1)}=(-6,2,0,1,1,-1,3,0,0,0,0)\;,\;
l^{(2)}=(0,0,1,0,0,1,-2,0,0,0,0) \;. }
We verify that the above basis for the reduced lattice
generates the cone ${\cal K}({\cal A},T_0)$ dual to ${\cal C}'({\cal A},T_0)$
for the maximal triangulation $T_0$ of $\Ds$.
However this is not a general phenomenon
as we will see in the example $X_{14}(7,3,2,1,1)$ presented
in sect.4.

The operators $\sqbox_l$ we deduce from the first two of \singIquot\ are
\eqn\singIbox{
\eqalign{
& \sqbox_{l_{\{2,5\}}}=\da{2}\da{5}-(\da{6})^2 \;\;, \cr
& \sqbox_{l_{\{1,3,4,6\}}}=\da{1}\da{3}\da{4}\da{6} - (\da{0})^4\da{9} \;\;,
\cr}
}
where $l_{\{2,5\}}=l_\tau^{(3)}+l_\tau^{(5)}+2 l_\tau^{(6)}$ and
$l_{\{1,3,4,6\}}=2l_\tau^{(1)}+l_\tau^{(2)}+l_\tau^{(3)}+2l_\tau^{(4)}$
are primitive relations for $\SDs_{reg}$.
These two operators are the analogues of \exIbeffac\ of the
nonsingular model, but with one crucial difference.
In this singular case, we do not have an order one differential
operator in the extended $\Ds$-hypergeometric system
 to eliminate $\da{9}$. In order to eliminate this we must
study the Jacobian ring of the hypersurface in detail. In \HKTYI,
a second order operator was found which has the form
\eqn\singItwo{
\da{0}\da{9}={12a_1a_2 \over a_0^2}(\da{0})^2 -
{24a_1a_2a_6 \over a_0^3}\da{0}\da{6}-
{12a_1a_6^2 \over a_0^3}\da{0}\da{5} \;\;,
}
when acting on the period $\Pi(a)$, see eq.(3.39) in \HKTYI.
Using this relation and the definition $x=a^{l^{(1)}}$
and $y=a^{l^{(2)}}$, a third order differential operator is
derived from $\sqbox_{l_{\{1,3,4,6\}}}$ after a factorization
$\theta_x$ from the left.
As is evident, the linear differential
operators represent relations among the monomial with the homogeneous
degree $d$ or the charge one in the Jacobian ring.
In contrast, the differential operator
\singItwo\ represents a relation among the monomials of charge
two. While the order one differential operators have been related
to the symmetry of the period under automorphisms and thus to
the combinatorial data of the polyhedron $\Ds$, the form of the
operators for the charge two monomials  above do not have a clear
description in terms of the combinatorial data. This is a typical feature
we encounter in the analysis of the singular models.
We observe that despite having to use
 charge two operators to factorize ${\cal D}_l$,
the principal part of the factorized operators still coincide with those
monomial generators of the defining ideal ${\cal I}_{quot}$ for
 $A^*_{toric}(X_\D)$
-- just as in the case of type I, II models.
 This means that the
structure of the local solutions are not affected by the
usage of the charge two operators. That is,
the same properties in \Ksol\ hold for type III models as
they do for type I and II models.
In our example here, the local solutions are described by the following
topological data;
\eqn\singItop{
\eqalign{
K^{cl}_{xxx}=2 &  \;\;,\;\; K^{cl}_{xxy}=3 \;\;,\;\;
K^{cl}_{xyy}=3    \;\;,\;\; K^{cl}_{yyy}=3 \;\;,  \cr
&
c_2\cdot J_x= 32 \;\; ,\;\; c_2\cdot J_y= 42\;\; . \cr }
}

We verify that the convex hull of the points $\{\ns{1},\cdots,\ns{10}\}$
is again reflexive and defines a family of Calabi-Yau manifolds $X_\Dp$ with
Hodge numbers\Hodge\ $h^{1,1}(X_\Dp)=6$ and
$h^{2,1}(X_\Dp)=71$. Thus this is a case in which a polyhedron $\Ds$
results in topologically distinct
Calabi-Yau manifolds $X_\Dp$ and $X_\D$ sitting inside two distinct
ambient spaces (because $\SDs_{reg}\neq\Sigma({\Delta'}^*)$ as we have seen).
In fact $X_\Dp$ is not even in the
list of 7,555 Laudau-Ginzburg models of \KlemmSchimmrigk.

Finally let us calculate the Stanley-Reisner ideal for the
triangulation $T_0$ in \singT. It is straightforward to see that
the ideal is generated by
\eqn\singIsimp{
D_2D_5 \;\;,\;\; D_1D_3D_4D_6 \;\;.
}
Since the model (or the fan $\SDs$) is only simplicial but not regular,
the odd homology groups of the singular toric variety
can have torsion. Thus we consider
the homology groups over $\IQ$. Then the groups are given by
the intersection ring \intring\ $\A^*(\IP_\SDs,\IQ)$ over $\IQ$. Thus
it is $\IQ[D_1,\cdots, D_6]/{\cal I}$
with the ideal ${\cal I}$ generated by \singIsimp\ and the linear
relations among the vertices $\{\ns{1},\cdots,\ns{6}\}$ as in $(ii)$
of \intring \Fulton.
The normalization of the ``volume form'' of this ring
becomes less clear because the
Euler number of the singular $\IP_\SDs$ is not given simply by
the number of the cones with maximal dimensions in $\SDs$.
However we know that the hypersurface $X_\D$ in general position
 does not meet the point
singularities of the $\IP_\SDs$, and the hypersurface divisor class is given by
$[X_\D]=D_1+\cdots+D_6$. Therefore we naturally introduce a
normalization of $\A^*(\IP_\SD,\IQ)$ using the Euler
number of the hypersurface, rather than that of the singular ambient space:
\eqn\normalz{
{ \prod_{i=1}^{p}(1+D_i) \over (1+[X_\D]) } [X_\D] \vert_{top} =
2\( h^{1,1}(X_\D) - h^{2,1}(X_\D) \) \;\;.
}
Here we evaluate the component of the top degree on the left hand
side and we use the Hodge numbers $h^{1,1}(X_\D)$ and $h^{2,1}(X_\D)$
in \Hodge. In the left hand side, we adopt the expression
$\prod (1+D_i)$ for the total Chern class\Fulton\
which is justified
for the nonsingular $\IP_\SDs$, but naively extended to our singular
case. We have verified experimentally that the normalization
\normalz\ indeed results in the right topological couplings and the
linear form $c_2\cdot J$'s.
We may summarize our observation in general,

\vskip0.5cm
\leftskip=0.8cm \rightskip=1.5cm
\vbox{
\noindent
{\it
For a smooth Calabi-Yau models $X_\D$ in a singular toric variety
$\IP_\SDs$, the intersection ring $\A_{toric}^*(X_\D,\IQ)$ is
given by $\A^*(\IP_\SDs,\IQ)/ Ann(D_1+\cdots+D_p)$ with
the normalization determined by \normalz.
}
\vskip-55pt
\eqn\singIint{ }
\vskip25pt
}
\leftskip=0cm \rightskip=0cm
\noindent
The effect of taking quotient by $Ann(D_1+\cdots+D_p)$ may be
replaced by the ideal quotient ${\cal I}_{quot}={\cal I}:(D_1+\cdots+D_p)$
as in the nonsingular case. In our example, it is easy to
derive the first two of \singIquot\ from \singIsimp\ via the ideal quotient.

Finally we note
that all notions in the theory of toric ideals apply to the singular cases
as well. Therefore it would be helpful to compare the Gr\"obner
fan of a singular model with that of a nonsingular model. By
an  analysis similar to \exIIels, we obtain the following elementary
relations for the model $\IP(4,3,2,2,1)$;
\eqn\singIels{
\eqalign{
&
l^{(1)}=(-6, 2, 0, 1, 1, -1, 3) \;\; ,\;\;
l^{(2)}=(0, 0, 1, 0, 0, 1, -2) \;\;,\;\; \cr
&
l^{(3)}=(-6, 2, 1, 1, 1, 0, 1) \;\;,\;\;
l^{(4)}=(-12, 4, 3, 2, 2, 1, 0) \;\;.\cr}
}
The universal Gr\"obner basis are determined from the zonotope
${\cal P}_\A$ as
\eqn\singIugbasis{
\eqalign{
{\cal U}_\A  =\{ &
y_1^8y_2^5y_3^4y_4^4y_5y_6^2-y_0^{24},
y_1^6y_2^4y_3^3y_4^3y_5y_6-y_0^{18}, \cr
&
y_1^6y_2^5y_3^3y_4^3y_5^2-y_0^{18}y_6,
y_1^4y_2^3y_3^2y_4^2y_5-y_0^{12},
y_1^2y_3y_4y_6^3-y_0^6y_5, \cr
&
y_1^2y_2y_3y_4y_6-y_0^6,
y_1^2y_2^2y_3y_4y_5-y_0^6y_6,
y_2y_5-y_6^2 \}. \cr }
}
In table 2, we present the cones in the Gr\"obner fan with the
ideals of the leading terms. There the cone $\tau_1$ corresponds
to the maximal triangulation $T_0$ \singT\ and should be compared
with $\tau_1$ in the table 1. The difference we should note is that the
ideal $LT_\omega({\cal I}_\A)$ is not radical and does not
coincide with the Stanley-Reisner ideal $ST_\omega$.
To see the consequence of this fact, recall that the generators
of the ideal $LT_\omega({\cal I}_\A)$ may be mapped to the symbol
of the differential operators $\sqbox_l$ by \mapind. As we
see in the table 2 explicitly, we simply obtain higher order differential
operators rather than \singIbox.

\vskip0.2cm
\vbox{{\offinterlineskip
\hrule
\halign{
& \vrule# & \quad # \hfill & \quad \hfil # \hfill & \quad\hfil # \hfill& # \cr
height2pt
& \omit    & \omit    & \omit &  \omit &\cr
 & cone
  & weight $\omega$
   & $LT_\omega({\cal I}_\A)$
    & rad($LT_\omega({\cal I}_\A)$)
     &\cr
height2pt & \omit & \omit & \omit & \omit &\cr
\noalign{\hrule}
height2pt & \omit & \omit & \omit & \omit &\cr
 & $\tau_1$
  & $(0,1,1,0,0,0,0)$
   & $\langle y_1^2y_2y_3y_4y_6, y_1^2y_3y_4y_6^3, y_2y_5 \rangle $
    & $\langle y_1y_3y_4y_6,  y_2y_5 \rangle $
     & \cr
height2pt & \omit & \omit & \omit & \omit &\cr
 & $\tau_2$
  & $(0,0,2,0,0,1,0)$
   & $\langle y_1^2y_2y_3y_4y_6,y_0^6y_5,  y_2y_5 \rangle $
    & $\langle y_1y_2y_3y_4y_6,  y_0y_5, y_2y_5 \rangle $
     & \cr
height2pt & \omit & \omit & \omit & \omit &\cr
 & $\tau_3$
  & $(1,0,0,0,0,1,0)$
   & $\langle y_2y_5,  y_0^6 \rangle $
    & $\langle y_2y_5,  y_0 \rangle $
     & \cr
height2pt & \omit & \omit & \omit & \omit &\cr
 & $\tau_4$
  & $(1,0,0,0,0,0,1)$
   & $\langle y_0^2,  y_6^2 \rangle $
    & $\langle y_0,  y_6 \rangle $
     & \cr
height2pt & \omit & \omit & \omit & \omit &\cr
 & $\tau_5$
  & $(1,0,0,0,0,0,7)$
   & $\langle y_1^2y_2y_3y_4y_6,y_6^2,y_0^6y_6, y_0^{12} \rangle $
    & $\langle y_0,  y_6 \rangle $
     & \cr
height2pt & \omit & \omit & \omit & \omit &\cr
 & $\tau_6$
 & $(0,0,0,0,0,1,3)$
  & $\langle y_6^2,y_0^6y_6,y_1^2y_2y_3y_4y_6,y_1^4y_2^3y_3^2y_4^2y_5 \rangle $
    & $\langle y_6, y_1y_2y_3y_4y_5 \rangle $
     & \cr
height2pt & \omit & \omit & \omit & \omit &\cr
 & $\tau_7$
 & $(0,1,0,0,0,1,1)$
  & $\langle y_6^2,y_1^2y_2y_3y_4y_6,y_1^2y_2^2y_3y_4y_5 \rangle $
    & $\langle y_6, y_1y_2y_3y_4y_5 \rangle $
     & \cr
height2pt & \omit & \omit & \omit &  \omit & \cr }
\hrule}}
\vskip0.2cm
\vbox{\baselineskip=12pt
\noindent
{\tabfont
Table 2. Gr\"obner cones with typical weights for $\IP(4,3,2,2,1)$.
         \par\noindent
The first cone $\tau_1$ corresponds to the maximal
triangulation $T_0$ in the text.  }  }
\vskip0.2cm

\vskip0.5cm

\vfill\eject

\newsec{ Applications to Mirror Symmetry and Mirror Map }

In application to mirror symmetry, the secondary fan
can be regarded as a collection of different phases of a type II string theory
compactified on a Calabi-Yau manifold (see for example \ref\wittenPhase{
E.Witten, Nucl.Phys.{\bf B403}(1993)159. }
\ref\morrisonPlesser{D.Morrison and Plesser,
Nucl.Phys,{\bf B404}(1995)279. }).
The triangulations of $\Ds$ which induce different subdivisions of the
fan $\SDs$, and their corresponding cones in the secondary fan are known
to have a clear physical meaning in terms of orbifold as well as the
smooth Calabi-Yau manifold. Among them, the maximal triangulations of
$\Ds$ or the finest refinements of the fan $\SDs$ constitute the Calabi-Yau
phase. In this phase we have the large radius limit of the smooth
Calabi-Yau manifold where the non-perturbative instanton corrections
are suppressed exponentially. The structure observed in \Ksol\ is
consistent with the quantum cohomology ring near the
large radius limit.

In this section, we use several models to show how our general
framework applies.

\subsec{Quantum cohomology ring}

Quantum cohomology ring is one of the nontrivial consequences of
the local operator algebra of the type II string theory compactified
on a Calabi-Yau threefold. In N=2 string theory,
 two different kinds of the local topological
operator algebras, called $(a,c)$- and
$(c,c)$-ring, correspond respectively to the $H^{1,1}$-type cohomology and the
$H^{2,1}$-type cohomology in the topological
$\sigma$-model
\ref\AspinwallMorr\ref{
Aspinwall and D. Morrison, Commun.Math.Phys.{\bf 151}245.}
\ref\topWitten{E.Witten, {\it Mirror Manifolds and Topological Field Theory}
in {\it ``Essays on Mirror Symmetry''}, ed. S.-T.Yau, (1992)
International Press.}. On physical ground, the
$(a,c)$-ring receives quantum corrections from
$\sigma$-model instantons whereas the $(c,c)$-ring does not
\DistlerGreene.
Mirror symmetry which exchanges the two provides a powerful
hypothesis to determine the quantum cohomology ring in terms
of the $(c,c)$-ring:
\eqn\qiso{
\oplus_{i=0}^3 H_q^{i,i}(X_\D) \cong
\oplus_{i=0}^3 H^{3-i,i}(X_{\Ds,a}) \;\;, }
where $q$ in the left hand side represents the quantum deformation
and $a$ in the right hand side represents the classical deformation
of the mirror hypersurface in \defeq. More precisely, we may regard
the right hand side as the Jacobian ring of the mirror hypersurface
and we can use the theory of variation of Hodge structures to study
this side. The isomorphism can then be realized in
terms of the flat coordinates on moduli spaces. This map is
called  the mirror map.
It is known that the mirror map has many remarkable
properties such as modular
property, integrality in the $q$-expansion e.t.c.\LYa.

In the classical limit, the instanton
 corrections in the quantum cohomology ring are exponentially suppressed.
The monodromy of the periods near the limit is
maximally unipotent \Morrison.
This is the property we established in general
in \Ksol\ for any maximal triangulation of the polyhedron $\Ds$ of type I
or II.  It is found
in \HKTYI\HKTYII\ that if we define the local variables via the basis
$\{l^{(k)}\}$ of the Mori cone by
\eqn\xvar{
x_k=(-1)^{l^{(k)}_0} a^{l^{(k)}} \;\;,
}
then we may express the mirror map as
\eqn\mirrormap{
t_j={1\over 2\pi i} { \drho{j} w_0(x,0) \over w_0(x,0) } \;\; ,
}
where $q_j = {\rm e}^{2\pi i t_j}$.
The inverse map is denoted as $x_k=x_k(q)$.
 The quantum couplings are related to the
geometrical couplings $K_{ijk}(x):=\int \Omega(x) \wedge
\pd_i\pd_j\pd_k \Omega(x)$
-- $\Omega(x)$ being the holomorphic threeform -- by
\eqn\qK{
K_{ijk}(q)= {1\over (2\pi i)^3 } \( {1\over w_0(x)} \)^2 \sum_{l,m,n}
K_{l m n}(x) {dx_l \over d t_i}{dx_m \over d t_j}{dx_n \over d t_k}
\vert_{x_k=x_k(q)} \;\;.
}
Special geometry in the $H^{2,1}$-moduli space
enables us to express the same couplings using
the so-called prepotential $F(t)$
\ref\stroming{A.Strominger, Commun.Math.Phys.{\bf 133}(1990)163.}:
\eqn\qKF{
K_{ijk}(q)= {1\over (2\pi i)^3 }
{\pd \; \over \pd t_i} {\pd \;\over \pd t_j}{\pd \;\over \pd t_k} F(t) \;\; .
}

For the prepotential, there is a concise formula given in \HKTYII\
based on the local structure \Ksol:
\eqn\PotF{
F(t)={1\over 2}\({1\over w_0(x)}\)^2 \big\{
w_0(x)D^{(3)}w_0(x) + \sum_l D_l^{(1)}w_0(x) D_l^{(2)}w_0(x) \big\}
\vert_{x_k=x_k(q)} \;\; ,
}
where we define
\eqn\solsD{
D^{(1)}_l= \drho{l} \;,\;
D^{(2)}_l= {1\over2}
\sum_{m,n}K_{lmn}^{cl}\drho{m}\drho{n} \;,\;
D^{(3)}=-{1\over6}
\sum_{l,m,n}K_{lmn}^{cl}\drho{l}\drho{m}\drho{n} \;.
}
It is also observed that the prepotential defined above has the following
asymptotic form with topological data in the leading terms, i.e.,
\eqn\Fasysmpt{
F(t)={1\over6}\sum_{i,j,k} K^{cl}_{ijk}t_i t_j t_k
     -\sum_k {(c_2\cdot J_k) \over 24} t_k
     -i {\zeta(3) \over 16 \pi^3 } \chi  + {\cal O}(q) \;\;,
}
where $\chi$ is the Euler number of $X_\D$ and the ${\cal O}(q)$-
terms represent the quantum corrections.
The first example understood was the
case of the quintic in $\IP^4$ studied by
Candelas et al\CandelasETAL.
We denote by $N^r(d)$ the predicted number of $\sigma$-model instantons with
multi degree $(d_1,\cdots,d_{h^{1,1}})$.
The genus one (string 1 loop) topological
amplitude\ref\BCOV{
M.Bershadsky, S.Cecotti, H.Ooguri and C.Vafa (with an appendix by S.Katz),
Nucl.Phys.{\bf B405}(1993)279.},
$F_1^{top}$ has the form
\eqn\FI{
F_1^{top}={\rm log}\big\{ ({1\over w_0})^{5-\chi/12}
{\pd (x_1,\cdots,x_{h^{1,1}}) \over
 \pd (t_1,\cdots,t_{h^{1,1}}) }
\prod_j dis_j^{r_j} \prod_i x_i^{s_i} \big\} + {\rm const.}
}
where the $dis_j$ are irreducible parts of the discriminant of the
hypersurface $X_\D$ and $r_j$ and $s_i$ are some parameters to be fixed
by the asymptotic form of the topological amplitude. It is known that
the amplitude has an expansion of the form
\eqn\FIexpand{
\eqalign{
F^{top}_1=& {\rm const.}-{2\pi i \over 12}\sum_k (c_2\cdot J_k)t_k  \cr
 & -\sum_{d} \big\{
2 N^e(d) {\rm log}\( \eta( q^d ) \)
+{1\over 6} N^r(d) {\rm log}(1-q^d) \big\} \;\;, \cr}
}
where $q^d=q_1^{d_1}\cdots q_{h^{1,1}}^{d_{h^{1,1}}}$ and the number
$N^e(d)$ is the prediction for the number of 1 loop
instantons, i.e., elliptic curves
in the Calabi-Yau manifold $X_\D$ with multi degree $n$.

In the following, based on our general observation \Ksol, we analyze
the large radius limit. In this paper, we will be mainly concerned with
the determination of the Picard-Fuchs operators from which we can
determine the quantum corrections in a straightforward way. For example
we can calculate the quantum corrected yukawa couplings \qK\ using the
Mathematica program INSTANTON appended to \HKTYII.
The required input
data come from the Picard-Fuchs operators and the classical couplings
given here in
appendix C. For interested reader, a complete list of the
Picard-Fuchs operators for the Calabi-Yau hypersurfaces with $h^{1,1}\leq
3$ is appended in the source file of this text \HLY.
The determination of the numbers $N^e(d)$ is a little involved because we need
to know the form of the discriminants of the hypersurfaces and need to
fix unknown parameters $r_i$ and $s_i$ in \FI. We will list,
in the appendix to the source file, the form of the discriminants for some
of our models.
However the detailed analysis, together
with the analysis of the conifold singularities where one Calabi-Yau model
may be connected to another (cf.
\ref\stromingII{A.Strominger,{\it Massless Black Holes and Conifolds in
String Theory}, hep-th/9504090.}
\ref\GMS{B.Greene, D.Morrison and A.Strominger,
       {\it Black Hole Condensation and the Unification of String Vacua},
       hep-th/9504145.}) will be presented elsewhere.

\vskip50pt

\subsec{Selected Examples}

\leftline{\undertext{ $X_{9}(3,2,2,1,1)^2_{-168}$ }}
This is a singular model of type III.
The polyhedron $\D(w)$ for this model has the vertices
\eqn\singIImus{
\eqalign{
&
\nu_1=(2,-1,-1,-1)  \;\;,\;\;
\nu_2=(-1,3,-1,-1) \;\;,\;\;
\nu_3=(-1,3,-1,0) \;\;, \cr
&
\nu_4=(-1,-1,3,-1) \;\;,\;\;
\nu_5=(-1,-1,3,0)  \;\;,\;\;
\nu_6=(-1,-1,-1,8) \;\;, \cr
&
\nu_7=(-1,-1,-1,-1) \;\;,\;\;
\nu_8=(0,2,-1,-1)   \;\;,\;\;
\nu_9=(0,-1,2,-1)   \;\;, \cr }
}
with respect to the basis $\{\Lambda_1,\cdots,\Lambda_4\}$ for the lattice
$H(w)$ defined after \exImus.  Then the vertices of the dual polyhedron
$\Ds(w)$ are given by
\eqn\singIInus{
\eqalign{
&
\ns{1}=(1,0,0,0) \;\;,\;\;
\ns{2}=(0,1,0,0) \;\;,\;\;
\ns{3}=(0,0,1,0) \;\;, \cr
&
\ns{4}=(0,0,0,1) \;\;,\;\;
\ns{5}=(-3,-2,-2,-1) \;\;,\;\;
\ns{6}=(-1,-1,-1,0) \;\;. \cr}
}
There are no integral points inside the polyhedron except the origin. For
the maximal triangulation of $\Ds(w)$, we obtain
\eqn\singIImT{
\eqalign{
T_0  =\{ &
\langle 0,3,4,5,6  \rangle , \langle 0,2,4,5,6 \rangle ,
\langle 0,1,3,5,6  \rangle , \langle 0,1,2,5,6 \rangle ,
\langle 0,1,3,4,6  \rangle ,  \cr &
\langle 0,1,2,4,6  \rangle , \langle 0,2,3,4,5  \rangle ,
\langle 0,1,2,3,5  \rangle ,
\langle 0,1,2,3,4  \rangle \}. \cr}
}
This triangulation induce the fan $\SDs$, however the resulting fan is
singular because the simplex $\langle 0,2,3,4,5 \rangle$ has volume
three. In fact we find two integral points
$\ns{7}={ (\ns{2}+\ns{3}) + 2(\ns{4}+\ns{5}) \over 3 }=(-2,-1,-1,0) \;,\;
\ns{8}={ 2(\ns{2}+\ns{3})+(\ns{4}+\ns{5}) \over 3 }=(-1,0,0,0)$ which are
inside the
cone spanned by $\ns{2},\ns{3},\ns{4}$ and $\ns{5}$ but outside the
polyhedron, indicating that this model is  of type III. As described in
the previous section, we subdivide the fan $\SDs$ by $\ns{7}$ and $\ns{8}$
to obtain a regular fan $\SDs_{reg}$.
The intersection ring $A^*(\IP_{\SDs_{reg}},\ZZ)$ is described by the ideal
$(i)$ in \intring\ with generators
\eqn\singIISRA{
\eqalign{
D_1D_7 & \;,\; D_1D_8 \;,\; D_6D_7 \;,\; D_6D_8 \;,\;
D_1D_4D_5 \;,\; D_2D_3D_6 \;,\; \cr
&
D_2D_3D_7 \;,\; D_4D_5D_8 \;,\;D_2D_3D_4D_5 \;, \cr}
}
and the linear relations $(ii)$ among the integral points
$\ns{1},\cdots,\ns{8}$. The hypersurface divisor
$[X_\D]=D_1+\cdots + D_6+2D_7+2D_8$ determines the ideal quotient ${\cal
I}_{quot}$. It is generated by
\eqn\singIIquot{
D_4D_5-D_4D_6+4D_3D_6 \;\;,\;\; D_1D_4D_5 \;\;,\;\; D_7 \;\;,\;\; D_8 \;,
}
together with the linear relations.
Starting from  those operators $\sqbox_l$ whose leading terms
match \singIISRA (under the correspondence
$\ta{i} \leftrightarrow D_i (i=1,\cdots,p)$), we can derive
the Picard-Fuchs operators via some nontrivial factorizations.

We first note that the generator $D_1D_4D_5$ induces
$l_{\{1,4,5\}}=(-1,1,0,0,1,1,$ $-2,0,0)$ in $L$. From this we immediately see
that the operator
\eqn\singIIdI{
\sqbox_{l_{\{1,4,5\}}}=\da{1}\da{4}\da{5}-\da{0}\(\da{6}\)^2
}
is one of the Picard-Fuchs operators. To find the other, we need to derive
the following relations from the analysis of the Jacobian ring of the
hypersurface;
\eqn\singIIjac{
\eqalign{
\da{0}\da{7}
&=
-{3a_1 \over a_0}\(\da{0}\)^2-{a_0\over a_6}\da{0}\da{8} \;, \cr
\da{0}\da{8}
&=
-{3 a_1a_4a_5 \over 16 a_2a_3a_6 }\(\da{6}\)^2
-{a_0a_4a_5 \over 16 a_2a_3a_6}\da{4}\da{5}   \cr
& \quad
+{a_0a_4 \over 16 a_2a_3}\da{4}\da{6}
-{a_0 \over 4a_2} \da{3}\da{6} \;. \cr}
}
The derivation of the above relations may be done most efficiently
by representing the hypersurface in terms of the homogeneous coordinate of
$\IP(3,2,2,1,1)$;
\eqn\singIIdefp{
W=z_1^3+z_2^4z_4+z_3^4z_5+z_4^9+z_5^9 \;\;.
}
The mirror of this hypersurface, whose period we are analyzing, can be
constructed by the transposition argument of Berglund and H\"ubsch
\BH. We consider the orbifold $\hat W/\ZZ_4\times\ZZ_9$ with
the transposed
potential $\hat W$. Relating to our toric description, we may write
the transposed
potential
\eqn\singIIdefphat{
\hat W=a_1z_1^3+a_2z_2^4+a_3z_3^4+a_4z_2z_4^9+a_5z_3z_5^9
+a_0z_1z_2z_3z_4z_5+a_6z_1z_4^4z_5^4 \;,
}
which is regarded as the degree 12 hypersurface in $\IP(4,3,3,1,1)$. Then
the integral points $\ns{7},\ns{8}$ are mapped, respectively,
to degree 24 (charge two)
monomials $z_2^2z_3^2z_4^6z_5^6$ and $z_2^3z_3^3z_4^3z_5^3$
under the monomial-divisor map\AGM.
The equations in \singIIjac\ represent the relations among
the charge two monomials in the Jacobian ring $\IC[z_1,\cdots,z_5]/( \pd
\hat W)$.

We now focus on the generators $D_1D_8$ and $D_6D_7$ which correspond to the
primitive collections whose primitive relations are
$l_{\{1,8\}}=(-2,1,0,0,0,0,0,0,1)$ and
$l_{\{6,7\}}=(0, 0, 0, 0, -1, -1, 1, 1, 0)$ in $L$. We find
that the operator
\eqn\singIIfactor{
{\cal O}=3{a_1a_2a_3a_6 \over a_0}\da{0}\sqbox_{l_{\{1,8\}}}-
{a_2a_3a_6^3 \over a_0^2 }\da{0}\sqbox_{l_{\{6,7\}}}
}
has the property that
$a_0{\cal O}\Pi(a)=(\ta{0}-2){\cal D}_2\tilde \Pi(a)$. Using this we
obtain a complete set of the Picard-Fuchs operators in the $\tilde
\Pi(a)$ gauge
\eqn\singIIPF{
\eqalign{
{\cal D}_1&=
\ty\(\ty-\tx\)^2-y\(3\tx+\ty+1\)\(3\tx-2\ty-1\)\(3\tx-2\ty\) \;,
\cr
{\cal D}_2&=
\(\tx-\ty\)^2+\(\tx-\ty\)\(3\tx-2\ty\)+4\tx\(3\tx-2\ty\) \cr
&
-48xy\(3\tx-2\ty-1\)\(3\tx+\ty+1\)-3y\(3\tx-2\ty-1\)\(3\tx-2\ty\) \cr
&
-48xy\(3\tx+\ty+3\)\(3\tx+\ty+1\)-16x\(\tx-\ty\)^2 \;,\cr}
}
where $x$ and $y$ are defined by \xvar\ using the basis
$l^{(1)}=(-3,0,1,1,-1,-1,3,0,0)$ and $l^{(2)}=(-1,1,0,0,1,1,-2,0,0)$
generating the Mori cone in the reduced lattice
$\{ l\in L \; \vert \;l_7=l_8=0 \;\}$.

Using the hypersurface divisor $[X_\D]=D_1+\cdots+D_6+2D_7+2D_8$
we determine the following topological data
\eqn\singIItopK{
\eqalign{
K_{xxx}^{cl}=6   &
\;\;,\;\; K_{xxy}^{cl}=9 \;\;,\;\;
K_{xyy}^{cl}=13 \;\;,\;\; K_{yyy}^{cl}=17 \;\;,\;\; \cr
&
c_2 \cdot J_x=48 \;\;, c_2\cdot J_y=74 \;.\cr }
}
According to the general form \Ksol, these topological data determine the
local solutions of the Picard-Fuchs equations \singIIPF\ near $x=y=0$.
We notice that this model has the same Hodge numbers as
the model $X_8(2,2,2,1,1)^2_{-168}$. However there is no linear
transformation which relates the
topological data: the cubic and linear forms of the two manifolds.
By Wall's theorem
\ref\wall{C.T.Wall, Invent.Math.{\bf 1}(1966)355.}
we see that the two manifolds are topologically distinct.

\vskip0.3cm

\leftline{\undertext{ Non LG model related to $X_{9}(3,2,2,1,1)^2_{-168}$ }}

For the model analyzed in the last subsection, we can verify that the
polyhedron $\Dps={\rm Conv.}(\{\ns{1},\cdots,\ns{8}\})$
is reflexive and the complete fan $\SDps$ for a triangulation of
$\Dps$ (,i.e., the
triangulation $T_A$ below) coincides with $\SDs_{reg}$. Therefore
we have another family of Calabi-Yau hypersurfaces $X_\Dp$ in the same ambient
space $\IP_{\SDps}=\IP_{\SDs_{reg}}$. The hypersurface represents
the divisor class
\eqn\singIIDp{
[X_\Dp]=D_1+\cdots+D_6+D_7+D_8 \;\;.
}

The dual polyhedron $\Dps$ is the
convex hull of the
points $\ns{1},\cdots, \ns{8}$. The polyhedron $\Dp$ has vertices
$\nu_2,\nu_3,\cdots,\nu_9$ (the corner $\nu_1$ is deleted) and
\eqn\singIIaddvert{
\eqalign{
&
\nu_{10}=(1,-1,-1,2) \;\;,\;\;
\nu_{11}=(1,0,-1,0)  \;\;,\;\;
\nu_{12}=(1,0,-1,-1) \;\;, \cr
&
\nu_{13}=(1,-1,-1,-1) \;\;,\;\;
\nu_{14}=(1,-1,0,0) \;\;,\;\;
\nu_{15}=(1,-1,0,-1) \;\;.}
}
By the formula \Hodge, we know that Hodge numbers of $X_\Dp$ are
$h^{1,1}=4$ and $h^{2,1}=85$. It turns out that this model is not in the
list of \KlemmSchimmrigk.  Also this model
gives an example of a topology change due to flop operations\AGM.

There are 37 triangulations for the polyhedron $\Dps$ and among them two
triangulations give us different resolutions of the ambient space. The
first one is the triangulation corresponding to
the subdivision $\SDs_{reg}$:
\eqn\singIITA{
\eqalign{
T_A & =\{
\langle 0,3,5,7,8 \rangle ,
\langle 0,2,5,7,8 \rangle ,
\langle 0,3,4,7,8 \rangle ,
\langle 0,2,3,5,8 \rangle ,
\langle 0,2,3,4,8 \rangle , \cr &
\langle 0,1,3,5,6 \rangle ,
\langle 0,1,2,5,6 \rangle ,
\langle 0,1,3,4,6 \rangle ,
\langle 0,1,2,4,6 \rangle ,
\langle 0,1,2,3,5 \rangle , \cr &
\langle 0,1,2,3,4 \rangle ,
\langle 0,3,4,5,7 \rangle ,
\langle 0,2,4,5,7 \rangle ,
\langle 0,3,4,5,6 \rangle ,
\langle 0,2,4,5,6 \rangle  \} .  \cr }
}
The second triangulation $T_B$ is $T_A$ but with the last four simplices
replaced by
\eqn\singIITB{
\langle 0,3,4,6,7 \rangle ,
\langle 0,2,5,6,7 \rangle ,
\langle 0,3,4,6,7 \rangle ,
\langle 0,2,4,6,7 \rangle \;\; . }
We verify that the difference in the two triangulations is due to
two different triangulation of
the two dimensional face (square) $\langle
\ns{4},\ns{7},\ns{5},\ns{6} \rangle$. They are
$\{
\langle \ns{4},\ns{5},\ns{6} \rangle \;,\;
\langle \ns{4},\ns{5},\ns{7} \rangle \} $ for $T_A$, and
$\{
\langle \ns{4},\ns{6},\ns{7} \rangle \;,\;
\langle \ns{5},\ns{6},\ns{7} \rangle \} $ for $T_B$.

For the triangulation $T_A$, we have in \singIISRA\ the generators of the
Stanley-Reisner ideal. Each generator $D_{i_1}D_{i_2}\cdots D_{i_k}$
determines uniquely the element $l_{\{i_1,i_2,\cdots,l_k\}}$ in the lattice
$L$, and in turn the operator  $\sqbox_{l_{\{i_1,i_2,\cdots,i_k\}}}$.
We observe that some combinations of the operators factorize to give a
complete set of Picard-Fuchs operators.
The principal parts of these operators generate
the ideal ${\cal I}_{quot}$ as in \singIIquot.
In appendix A, we list the resulting
Picard-Fuchs operators in terms of the local coordinates $x,y,z$ and $w$
defined by
$l_A^{(1)}=(-1,1,0,0,1,1,-2,0,0)$,
$l_A^{(2)}=(-1,0,0,0,1,1,0,-2,0)$,
$l_A^{(3)}=(-1,0,1,1,0,0,0,1,-2)$ and
$l_A^{(4)}=(0,0,0,0,-1,-1,1,1,0)$, respectively.
The intersection ring \Atoric\ determines the topological data as
follows;
\eqn\singIItopA{
\eqalign{ &
K_{xxx}^{A,cl}=17 \;,\;
K_{xxy}^{A,cl}=26 \;,\;
K_{xyy}^{A,cl}=36 \;,\;
K_{yyy}^{A,cl}=46 \;,\;
K_{xxz}^{A,cl}=13 \;,\; \cr
&
K_{xyz}^{A,cl}=18 \;,\;
K_{yyz}^{A,cl}=23 \;,\;
K_{xzz}^{A,cl}=9 \;,\;
K_{yzz}^{A,cl}=11 \;,\;
K_{zzz}^{A,cl}=4 \;,\; \cr
&
K_{xxw}^{A,cl}= 39\;,\;
K_{xyw}^{A,cl}= 54\;,\;
K_{yyw}^{A,cl}= 72\;,\;
K_{xzw}^{A,cl}= 27\;,\;
K_{yzw}^{A,cl}= 36\;,\;  \cr
&
K_{zzw}^{A,cl}= 18\;,\;
K_{xww}^{A,cl}= 81\;,\;
K_{yww}^{A,cl}= 108\;,\;
K_{zww}^{A,cl}= 54\;,\;
K_{www}^{A,cl}= 162\;,\;  \cr
&
\hskip0.5cm
c_2\cdot J_x^A=74 \;,\;
c_2\cdot J_y^A=100 \;,\;
c_2\cdot J_z^A=52 \;,\;
c_2\cdot J_w^A=144 \;.\;  \cr}
}

For the triangulation $T_B$, we find the generators for the Stanley-Reisner
ideal
\eqn\singIISRB{
D_1D_7\;,\; D_1D_8 \;,\; D_4D_5 \;,\; D_6D_8 \;,\;
D_2D_3D_6 \;,\; D_2D_3D_7 \;.
}
We observe again some (less trivial) factorizations among the operators
$\{ \sqbox_{l_{\{i_1,i_2,\cdots,i_k\}}} \}$ and their combinations
in order to obtain
the Picard-Fuchs operators listed in appendix A.
The local coordinates $x,y,z$ and $w$ for this
triangulation are defined by
$l_B^{(1)}=l_A^{(1)}+l_A^{(4)}$,
$l_B^{(2)}=l_A^{(2)}+l_A^{(4)}$,
$l_B^{(3)}=l_A^{(3)}+l_A^{(4)}$ and $l_B^{(4)}=-l_A^{(4)}$, respectively.
Then the topological data turns out to be
\eqn\singIItopB{
\eqalign{
K_{xxw}^{B,cl}&= 4\;,\;
K_{xyw}^{B,cl}= 8\;,\;
K_{yyw}^{B,cl}= 10\;,\;
K_{xzw}^{B,cl}= 4\;,\;
K_{yzw}^{B,cl}= 5\;,\;
K_{zzw}^{B,cl}= 2\;,\;  \cr
&
c_2\cdot J_x^B=74 \;,\;
c_2\cdot J_y^B=100 \;,\;
c_2\cdot J_z^B=52 \;,\;
c_2\cdot J_w^B=24 \;,\; \cr}
}
where the cubic couplings among $J_x,J_y,J_z$ are the same as in
\singIItopA\ and $K^{B,cl}_{www}=K^{B,cl}_{*ww}=0 \; (*=x,y,z)$.
As we observe in \singIItopB, the topological data for the phase B indicate
that the Calabi-Yau hypersurface has the property of a K3 fibration\KLM. In
fact, we verify that $\Ds_{K3}:={\rm Conv.}(\{ \ns{1},\ns{2},\ns{3},
\ns{6},\ns{7},\ns{8} \})$ is a three dimensional reflexive polyhedron.
We observe that $c_2\cdot J_i=24$ for some $i$ (cf.\singIItopB)
is necessary for the Calabi-Yau hypersurface to contain a $K3$.
We also remark that  the existence of
a three dimensional reflexive polyhedron $\Ds_{K3}$ in $\Ds$ does not
always yield the above topological condition. We will return to this
point later in the final section.

\vskip0.3cm

\leftline{\undertext{ $X_{14}(7,3,2,1,1)^2_{-260}$ }}

This model provides us an example in which we have two different
resolutions of point singularities in the ambient space, however the
difference of the two resolutions does not affect the topology of the
Calabi-Yau hypersurface. This model has also been solved in \BKK.

\goodbreak\midinsert
\centerline{\epsfxsize 3.5truein\epsfbox{fig2.eps}}
\smallskip
\noindent{\tabfont
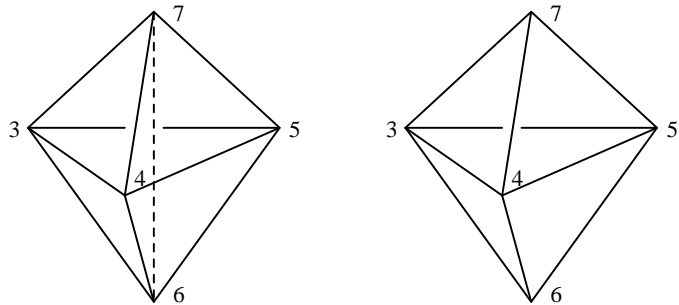
Fig.2 Two different triangulations $T_A$ (left) and $T_B$ (right) for the
models $X_7(7,3,2,1,1)$
\par\noindent
In the left, we see three 3-simplices whereas in the right we
see two 3-simplices. This results in different regular fans
$\SDs^A_{reg}$ and $\SDs^B_{reg}$ for the different desingularizarions of
the ambient space. However the Calabi-Yau hypersurfaces in them have the same
topology.
 }
\smallskip\endinsert

Let us summarize the toric data for this model. The reflexive polyhedron we
consider is given by the convex hull of the following integral points;
\eqn\singIIImus{
\eqalign{
&
\nu_1=(1,-1,-1,-1) \;,\;
\nu_2=(-1,3,0,-1)  \;,\;
\nu_3=(-1,3,-1,-1) \;,\; \cr
&
\nu_4=(-1,3,-1,1)  \;,\;
\nu_5=(-1,-1,6,-1)  \;,\;
\nu_6=(-1,-1,-1,13) \;,\; \cr
&
\nu_7=(-1,-1,-1,-1) \;, \cr}
}
with respect to the basis $\{ \Lambda_1,\cdots,\Lambda_4 \}$ of $H(w)$
given after\exImus. Then the vertices of the dual polyhedron $\Ds(w)$
are
\eqn\singIIInus{
\eqalign{
&
\ns{1}=(1,0,0,0) \;\;,\;\;
\ns{2}=(0,1,0,0) \;\;,\;\;
\ns{3}=(0,0,1,0) \;\;, \cr
&
\ns{4}=(0,0,0,1) \;\;,\;\;
\ns{5}=(-7,-3,-2,-1) \;\;,\;\;
\ns{6}=(-2,-1,0,0) \;\;. \cr}
}
We will find one point $\ns{7}=(-1,0,0,0)$ on a codimension one face $\langle
\ns{2},\ns{3},\ns{4},\ns{5},\ns{6} \rangle$.
If we triangulate the polyhedron $\Ds(w)$, we will find
the following two different
triangulations $T_A$ and $T_B$ which induce the
complete fans $\SDs^A$ and $\SDs^B$ respectively;
\eqn\singIIITA{
\eqalign{
T_A & =\{
\langle 0,4,5,6,7 \rangle ,
\langle 0,3,5,6,7 \rangle ,
\langle 0,3,4,6,7 \rangle ,
\langle 0,2,4,5,7 \rangle , \cr &
\langle 0,2,3,5,7 \rangle ,
\langle 0,2,3,4,7 \rangle ,
\langle 0,1,4,5,6 \rangle ,
\langle 0,1,3,5,6 \rangle , \cr &
\langle 0,1,3,4,6 \rangle ,
\langle 0,1,2,4,5 \rangle ,
\langle 0,1,2,3,5 \rangle ,
\langle 0,1,2,3,4 \rangle    \} .  \cr }
}
$T_B$ can be obtained from $T_A$ by replacing
the first three simplices of $T_A$ by
$\langle 0,3,4,5,7 \rangle$ and $\langle 0,3,4,5,6 \rangle$. The difference
between $T_A$ and $T_B$ are depicted in fig.2. Since it turns out that some of
the cones in the fan $\SDs$ are singular for both triangulations, we need to
subdivide them. In the case of $T_A$, we find the following integral
points make the cones regular;
\eqn\singIIIaddA{
\eqalign{
&
\ns{8}={1\over2}(\ns{4}+\ns{5}+\ns{6}+\ns{7}) \;,\;
\ns{9}={1\over2}(\ns{2}+\ns{4}+\ns{5}+\ns{7}) \;,\;  \cr
&
\ns{10}={1\over2}(\ns{1}+\ns{4}+\ns{5}+\ns{6}) \;,\;
\ns{11}={1\over2}(\ns{1}+\ns{2}+\ns{4}+\ns{5}) \;,\; \cr}
}
and for $T_B$ we find
\eqn\singIIIaddB{
\eqalign{
&
\ns{8}={1\over3}(\ns{3}+\ns{7})+{2\over3}(\ns{4}+\ns{5}) \;,\;
\ns{9}={1\over2}(\ns{2}+\ns{4}+\ns{5}+\ns{7}) \;,\;  \cr
&
\ns{10}={1\over2}(\ns{1}+\ns{4}+\ns{5}+\ns{6}) \;,\;
\ns{11}={1\over2}(\ns{1}+\ns{2}+\ns{4}+\ns{5}) \;,\; \cr
&
\ns{12}={1\over3}(\ns{4}+\ns{5})+{2\over3}(\ns{3}+\ns{7}) \;.\; \cr}
}
Subdividing $\SDs^A$ and $\SDs^B$ by these integral points
results in the regular fans $\SDs^A_{reg}$ and $\SDs^B_{reg}$,
respectively, both of which do not come from any triangulation of the
polyhedron $\Dps$ -- the convex hull of all the integral points.
Using each of the two
 regular fans, we determine the basis for the K\"ahler cone,
and the Mori cone of the ambient space.  We summarize in
appendix B the bases
$\{\eta_A^1,\cdots,\eta_A^8\}$ and $\{\eta_B^1,\cdots,\eta_B^9 \}$
of the Mori cones
for $\SDs^A_{reg}$ and $\SDs^B_{reg}$ respectively.
We observe that the Mori cones for both ambient spaces are not simplicial,
implying that neither are the K\"ahler cones of the ambient spaces.

The divisor $[X_\D]$ of the form \ndiv\ determines the same
intersection ring for the two hypersurfaces, and for both cases
we find that the divisors $D_i \; (i\geq 7)$ decouple. In
fact the ideal ${\cal I}_{quot}$ is generated by
\eqn\singIIIquot{
D_2D_6 \;,\; D_3D_4D_5 \;,\; D_i \; (i\geq 7) \;,
}
together with the linear relations $(ii)$ in \intring.
We remark that in this model both $\{\ns{2},\ns{6}\}$ and
$\{\ns{3},\ns{4},\ns{5}\}$ are the primitive collections of
$\SD_{reg}^A$ and $\SD_{reg}^B$.
The reduced lattices which we denote ${L_A}', {L_B}'$
in the two cases are generated by
\eqn\singIIIsubLA{
\eqalign{
&
l_A^{(1)}=
3\eta_A^2+\eta_A^3+2\eta_A^4+2\eta_A^5+4\eta_A^6 =
(-4,2,1,0,0,0,1,0,\cdots,0) ,\cr
&
l_A^{(2)}=
2\eta_A^1+\eta_A^2+\eta_A^3+\eta_A^7+2\eta_A^8  =
(-2,1,0,2,1,1,-3,0,\cdots,0) , \cr}
}
for ${L_A}'$ and
\eqn\singIIIsubLB{
\eqalign{
&
l_B^{(1)}=
3\eta_B^2+6\eta_B^3+7\eta_B^4+\eta_B^5+2\eta_B^6 =
(0,0,1,-4,-2,-2,7,0,\cdots,0), \cr
&
l_B^{(2)}= \qquad\qquad
2\eta_B^1+\eta_B^7 \qquad\qquad =
(-2,1,0,2,1,1,-3,0,\cdots,0) , \cr}
}
for ${L_B}'$.
We remark that the Mori cone for for the ambient spaces are not
simplicial but their intersection with ${L}_\IR'$'s are.
 We also note that two restricted cones for $\SDs^A_{reg}$
and $\SDs^B_{reg}$ have an
intersection, in fact the former is included in the latter since
$l_A^{(1)}=l_B^{(1)}$ and $l_A^{(2)}=l_B^{(2)}+2l_B^{(1)}$.
We draw in fig.3 the restricted K\"ahler cones in the secondary fan
for the polyhedron $\Ds$, more precisely in the secondary fan for
the point configurations $\ns{0},\ns{1},\cdots,\ns{6}$
( we delete the point $\ns{7}$ corresponding to automorphisms).
Since ${\cal I}_{quot}$ is the same for both $\SDs^A_{reg}$ and
$\SDs^B_{reg}$,
we expect that the two triangulations define the same Calabi-Yau hypersurface
in different ambient spaces, i.e., the only difference is in the topology
of the ambient space which is irrelevant to the hypersurface.

\goodbreak\midinsert
\centerline{\epsfxsize 2truein\epsfbox{fig3.eps}}
\smallskip
\noindent{\tabfont
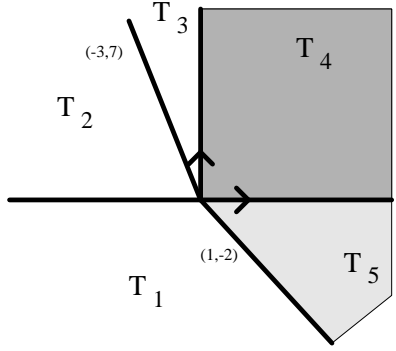
Fig.3 The secondary fan for the polyhedron $\Ds(7,3,2,1,1)$
\par\noindent
The K\"ahler cones of the smooth ambient spaces $\IP_{\SDs^A_{reg}}$ and
$\IP_{\SDs^B_{reg}}$ have different restrictions to the secondary fan.
The restricted K\"ahler cone for the former space is given
by the union of the cones parametrized by $T_4$ and $T_5$, while for the
latter it is given by the cone parametrized $T_4$.
 }
\smallskip\endinsert

Now we derive the Picard-Fuchs operators based on the triangulation
$T_A$.  We note that this model is of type III with non-trivial
automorphisms. This is the most general situation. The point $\ns{7}$ on a
codimension one face is a root vector in \rootsys\ for the fan $\SD$.
According
to \lieaction, this results in the following  linear operator annihilating
the periods:
\eqn\singIIIZp{
{\cal Z}'=2a_1\da{0}+a_0\da{7} \;\;.
}
Now we look at the operators which correspond to the first two generators
in \singIIIquot,
\eqn\singIIIsqbox{
\eqalign{
\sqbox_{l_{\{2,6\}}}&=\da{2}\da{6}-\(\da{7}\)^2 \;\;,\;\; \cr
\sqbox_{l_{\{3,4,5\}}}&=
\da{3}\da{4}\da{5}-\da{0}\da{6}\da{8} \;\;. \cr}
}
Starting with these operators, we derive the Picard-Fuchs operators
for the period restricted to the sublattice \singIIIsubLA. It is easy
to see the first operator $\sqbox_{l_{\{2,6\}}}$ combined with the
linear operator \singIIIZp\ results in a second order differential
operator. For the operator $\sqbox_{l_{\{3,4,5\}}}$, we need to look
into the structure of the Jacobian ring of the hypersurface.
For this, as in the previous example, it would be most efficient
to express the hypersurface in terms of the homogeneous coordinates:
\eqn\singIIIW{
W=z_1^2+z_2^4z_3+z_3^7+z_4^{14}+z_5^{14} \;\;,
}
in $\IP(7,3,2,1,1)$. Then the transposition argument in \BH\ applied to
this hypersurface indicates that the mirror is given by the orbifold
$\hat W/\ZZ_2\times\ZZ_{14}$ with
\eqn\singIIIhatW{
\eqalign{
\hat W=a_1z_1^2 & +a_2z_2^4+a_3z_2z_3^7+a_4z_4^{14}+a_5z_5^{14}
+a_0z_1z_2z_3z_4z_5  \cr
&
+a_6z_3^4z_4^4z_5^4+a_7z_2^2z_3^2z_4^2z_5^2 \;,\cr}
}
in $\IP(14,7,3,2,2)$. We note that, in this form, the automorphism used for
\singIIIZp\ is identified with
\eqn\singIIIauto{
z_1 \mapsto z_1+\varepsilon \; z_2z_3z_4z_5 \;\;, \;\;
z_i \mapsto z_i \;\; (i \geq 2) \;\;,
}
in infinitesimal form. The deformation parameters $a_8,\cdots,a_{11}$
corresponds to the degree 56 (charge two) monomials
$z_2z_3^3z_4^{10}z_5^{10} \,,\, z_2^3z_3z_4^8z_5^8 \,,\,
z_1z_3^2z_4^9z_5^9$ and $z_1z_2^2z_4^7z_5^7$, respectively.
Since we can verify $\( 1-{64 a_1^2a_2a_6 \over a_0^4}\)
z_2z_3^3z_4^{10}z_5^{10} =
{112a_1^2 a_2a_3 \over a_0^4}z_2^2z_3^6z_4^6z_5^6 +
{2 a_1a_3 \over a_0^2} z_3^8z_4^8z_5^8 $ modulo terms in the Jacobian
ring $(\pd \hat W )$ which vanish inside the period integral, we have
the relation
\eqn\singIIIrel{
\(1-{64 a_1^2a_2a_6 \over a_0^4}\)\da{0}\da{8}=
-224{a_1^3a_2a_3 \over a_0^5} \da{6}\da{0} +
2{a_1a_3 \over a_0^2 } \(\da{6}\)^2 \;, }
where we use \singIIIZp\ in the derivation. If we combine \singIIIrel\
with the operator $\sqbox_{l_{\{3,4,5\}}}$ in \singIIIsqbox, we will obtain
a third order differential operator. Thus we obtain the Picard-Fuchs
operators which determine the local solutions with the property \Ksol;
\eqn\singIIIPF{
\eqalign{
{\cal D}_1&= \tx\(\tx-3\ty\)-4x\(2\ty+4\tx+3\)\(2\ty+4\tx+1\)\cr
{\cal D}_2&= \(1-64x\)^2\ty^3 \cr
&
-64\{112 x^2y\(\tx-3\ty\)\(2\ty+4\tx+1\)+
x y\(\tx-3\ty-1\)\(\tx-3\ty\) \}   \cr
&
-\(1-64x\)\{112xy\(\tx-3\ty-1\)\(\tx-3\ty\)\(2\ty+4\tx+1\)   \cr
& \qquad \qquad \qquad
+y\(\tx-3\ty-2\)\(\tx-3\ty-1\)\(\tx-3\ty\) \} \cr}
}
with
$x=a^{l_A^{(1)}}={a_1a_2a_6 \over a_0^4 } , \;
 y=a^{l_A^{(2)}}={a_1a_3^2a_4a_5 \over a_0^2 a_6^3 } $.
The topological data for the local solutions about $x=y=0$ are given by
\eqn\singIIItop{
\eqalign{
K_{xxx}^{cl,A}= 9 \;,\;
&
K_{xxy}^{cl,A}= 3 \;,\;
K_{xyy}^{cl,A}= 1 \;,\;
K_{yyy}^{cl,A}= 0 \;,\; \cr
&
c_2\cdot J_x^A=66 \;\;,\;\;
c_2\cdot J_y^A=24 \;. \cr}
}
The analysis for $T_B$ is the same as the above and the Picard-Fuchs
operators are given by \singIIIPF\ with the variables $(x_A,y_A):=(x,y)$
changed to $(x_B,y_B)$ under the relations $ x_A=x_By_B^2, y_A=y_B$. The
topological data are connected by the linear relations which results from
these relations.

\vfill\eject

\newsec{Conclusion and Discussions}

We have analyzed the GKZ hypergeometric system -- which we
call $\Ds$-hypergeometric -- for a reflexive
polyhedron. The characteristic feature of
this system in mirror symmetry
is that it is $T$-resonant in general. Especially, for a
maximal triangulation $T_0$ of the polyhedron $\Ds$, the
monodromy of this system becomes
maximally unipotent. We have found close relationships between the
Stanley-Reisner ideal for the triangulation $T_0$ and the ring of the
leading terms of the $\Ds$-hypergeometric system at the maximally
unipotent point. For the models of type I and II,
 we have proved these two
ideals are actually equal,
using the general theory of toric ideals. We have found
 a closed formula for the local solutions near the
maximally unipotent point, in terms of the intersection form.
As was observed in \HKTYI\HKTYII, the $\Ds$-hypergeometric system is
reducible. If we extract the irreducible part of the system by
factoring out the operator $\ta{0}$, the resulting system gives a
sufficient set of differential operators to determine the quantum geometry
of moduli space.
We have verified our observations for the Calabi-Yau hypersurfaces in weighted
projective spaces up to $h^{1,1}\leq 3$, including models of type III.

In the table of appendix C, we have summarized the topological data for
each models.  There we can see several isomorphisms or relations
between different models. For example we have
\X{14}{\II}(7,2,2,2,1){2}{-240} $\cong$ \X{8}{\II}(3,1,1,1,1){2}{-240},
\X{15}{\III}(5,3,3,3,1){3}{-144} $\cong$ \X{10}{\III}(3,2,2,2,1){3}{-144}
and \X{18}{\II}(9,4,2,2,1){3}{-240} $\cong$ \X{9}{\I}(4,2,1,1,1){3}{-240},
all of which can be explained by a
fractional change of the variables
\ref\frac{M.Lynker and R.Schimmrigk, Phys.Lett.{B249}(1990)237.}.
Also there can be a reflexive polyhedron
$\Ds(w')$ in another reflexive polyhedron $\Ds(w)$\footnote{$^\dagger$}{
This observation has also been made in ref.\BKK.}.
For example, by listing all integral points in the polyhedra, we see
$\Ds(2,2,2,1,1) \subset \Ds(3,3,3,2,1) \;,\; \Ds(6,2,2,1,1) \subset
\Ds(9,3,3,2,1)$ and $\Ds(1^5), \Ds(2,1^4) \subset \Ds(3,2,2,2,1) \subset
\Ds(5,3,3,3,1) $. Since all integral points in $\Ds(w')$ is are contained in
$\Ds(w)$, the inclusion relation $\Ds(w') \subset \Ds(w)$ implies that the fan
$\Sigma(\Ds(w))$ is a refinement of the fan $\Sigma(\Ds(w'))$.  This reminds
us the cases we encountered in the singular models of type III, in which we
found that topologically different Calabi-Yau hypersurfaces can sit in the
same ambient space.

To see the details, let us consider the case
$\Ds(6,2,2,1,1) \subset \Ds(9,3,3,2,1)$.
The integral points in $\Ds(9,3,3,2,1)$ with respect to the basis
given after \exImus\ are
$
\ns{0}=(0,0,0,0), \;
\ns{1}=(1,0,0,0), \;
\ns{2}=(0,1,0,0), \;
\ns{3}=(0,0,1,0), \;
\ns{4}=(0,0,0,1), \;
\ns{5}=(-9,-3,-3,-2), \;
\ns{6}=(-6,-2,-2,-1), \;
\ns{7}=(-3,-1,-1, 0)$ and $\ns{8}=(-1,0,0,0)$, where the last points
$\ns{8}$ is on a codimension one face of the polyhedron. The
polyhedron $\Ds(6,2,2,1,1)$ has integral points
$\ns{0},\ns{1},\ns{2},\ns{3},\ns{4},\ns{6},\ns{7},\ns{8}$,
where the point $\ns{8}$ is also on a codimension one face.
 Therefore $\Sigma(\Ds(9,3,3,2,1))$ is a refinement of
the fan $\Sigma(\Ds(6,2,2,1,1))$, and we will have two different Calabi-Yau
hypersurfaces in the same ambient space $\IP_{\Sigma(\Ds(9,3,3,2,1))}$
\footnote{$^\ddagger$}{ Since the ambient space is still singular, we need
further subdivisions of some cones. However the following arguments
are valid for the fully resolved ambient space.}. According to \ndiv,
the divisor for the hypersurface is given by
\eqn\divI{
[X_{\D(w)}]=D_1+D_2+D_3+D_4+D_5+D_6+D_7+D_8 \;,
}
for the model \X{18}{{\,}}(9,3,3,2,1){3}{-186}; and
\eqn\divII{
[X_{\D(w')}]=D_1+D_2+D_3+D_4+2 D_5+D_6+D_7+D_8 \;,
}
for the model \X{12}{{\,}}(6,2,2,1,1){2}{-252}. This can also
be understood by the fractional transformation on the defining polynomial.
The polynomial $W(z)=\hat W(z)$ for the mirror of
\X{18}{{\,}}(9,3,3,2,1){3}{-186} is
\eqn\potI{
\hat
W=a_1z_1^2 + a_2z_2^6 + a_3z_3^6 + a_4z_4^9 + a_5z_5^{18}
 + a_0z_1z_2z_3z_4z_5 + a_6z_4^3z_5^{12} + a_7z_4^6z_5^6
\;\;,
}
in $\IP(9,3,3,2,1)/(\ZZ_6)^2$, where the deformation by $a_8$, which
corresponds to the divisor $D_8$, is eliminated using the automorphism. Now
consider the transformation
$\xi_i=z_i \; (i=1,2,3), \xi_4=z_4^{3/4}, \, \xi_5=z_4^{1/4}z_5$. Then the
potential becomes, if we set $a_5=0$,
\eqn\potII{
\hat W(\xi)=
a_1\xi_1^2 + a_2 \xi_2^6 + a_3 \xi_3^6 + a_4 \xi_4^{12}
+a_0 \xi_1\xi_2\xi_3\xi_4\xi_5 +a_6 \xi_5^{12} + a_7 \xi_4^6\xi_5^6 \;,
}
which can be regarded as a hypersurface in $\IP(6,2,2,1,1)/(\ZZ_6^2\times
\ZZ_{12})$, the mirror of \X{12}{{\,}}(6,2,2,1,1){2}{-252}.
The additional quotient by $\ZZ_{12}$ comes from the identification
$(\xi_4,\xi_5)\equiv (\alpha^4\xi_4,\alpha\xi_5)$ with $\alpha^4=1$ (, see
\HKTYI\ for the detailed form of the actions for $\ZZ_6^2$).
The Mori cone of each model may be obtained by restricting the Mori cone
of the ambient space to the sublattice $L'$, namely $l\in L$ with $l_8=0$
for $\Ds(9,3,3,2,1)$ and $l_5=l_8=$ for $\Ds(6,2,2,1,1)$. Thus the
inclusion of the dual polyhedron, $\Ds(w') \subset \Ds(w)$, implies an
embedding of the (quantum) K\"ahler moduli of $X_{\D(w')}$ to that of
$X_{\D(w)}$, or equivalently under mirror symmetry,
the complex structure moduli for the mirror
$X_{\Ds(w')}$ to that of $X_{\Ds(w)}$.

As a different kind of inclusion relation, we also observe that the dual
polyhedron $\Ds_{K3}(w')$ for some $K3$ hypersurface\ref\yonemura{
T.Yomemura, T\^ohoku Math. J. {\bf 42}(1990),351.} sits inside
the polyhedron $\Ds(w)$ for a Calabi-Yau hypersurface. It has also been
observed that if, in addition, we have the following specific form of the
topological data; $c_2\cdot K =24, \; J\cdot K \cdot K =
K \cdot K \cdot K =0$ for some divisor class $K$, then
the following ``CY-K3 correspondence'' occurs: the Picard-Fuchs
operators for the Calabi-Yau manifold specialize to those for a $K3$-model
under a suitable limit of the variables. In our list, the following models
shows this specific properties;
{\X{8}{{\,}}(2,2,2,1,1){2}{-168}}, {\X{12}{{\,}}(6,2,2,1,1){2}{-252}},
{\X{12}{{\,}}(3,3,3,2,1){3}{-126}},
{\X{18}{{\,}}(9,3,3,2,1){3}{-186}}, {\X{24}{{\,}}(12,8,2,1,1){3}{-480}},
{\X{10}{{\,}}(4,2,2,1,1){3}{-192}} and
{\X{16}{{\,}}(8,3,3,1,1){3}{-256}}.
Also our non Landau-Ginzburg model found in relation to
\X{9}{{\,}}(3,2,2,1,1){2}{-168} shows this property as well. The K3 polyhedron
$\Ds_{K3}$ contained in the reflexive polyhedron $\Dps$
provides an example of non Landau-Ginzburg K3 hypersurface. We have noticed
that the specific form of the topological data depends on how we
triangulate the polyhedron, namely in this example,
the CY-K3 correspondence occurs only in the phase B(see \singIItopB).
  Some of the models where the CY-K3 correspondence occurs
 has been studied extensively, and has provided
strong evidence for the so-called heterotic-type II string
duality\KachruVafa\KLM\KKLMV.
We believe that our general framework outlined here
will provide powerful techniques for the
studying questions in heterotic-type II duality.

\vskip0.5cm

\noindent{\bf Acknowledgements:} S.H. would like to
thank the Department of Mathematics, Harvard University for the financial
support for his visit, during which this work has been completed.
S.H. is supported in part by Grant-in-Aid for Scientific Research on Priority
Area 231 ``Infinite Analysis'', Japan Ministry of Education. B.H.L.
thanks J. Dalbec and G. Zuckerman for discussions. S.H. and B.H.L thank
A.Klemm for discussions and communications.

\vskip0.5cm

\vfill\eject

\appendix{A}{Picard-Fuchs equations for the model in section 3.}

This non Landau-Ginzburg model is defined by the reflexive
polyhedron $\Dps$ which has the
property $\SDps=\SDs_{reg}$ for $\Ds=\Ds(3,2,2,1,1)$.
There are two Calabi-Yau phases, phase A and phase B, which
are connected by flop operations.

\leftline{\undertext{Phase A:} }
\eqn\pfA{
\eqalign{
{\cal D}_1 &=
\tx\( \tz-2\ty+\tw\)  - xw\(\tx+\tz+\ty+11\)\(\tw-2\tx\) \cr
{\cal D}_2 &=
\tx\(\ty-2\tz\)-xyw^2\(\tx+\tz+\ty+2\)\(\tx+\tz+\ty+1\) \cr
{\cal D}_3 &=
\(\tw-2\tx\)\(\tz-2\ty+\tw\)-w\(\tx+\ty-\tw\)\(\tx+\ty-\tw\) \cr
{\cal D}_4 &=
\(2\tx-\tw\)\(2\tz-\ty\)-yw\(\tx+\tz+\ty+1\)\(\tz-2\ty+\tw\) \cr
{\cal D}_5 &=
\(\tx+\ty-\tw\)^2 \tx-x\(\tx+\tz+\ty+1\)\(\tz-2\ty+\tw-1\)
\(\tz-2\ty+\tw\) \cr
{\cal D}_6 &=(\tx+\ty-\tw)^2(\ty-2\tz)  \cr
&
  -y(\tx+\tz+\ty+1)(\tz-2\ty+\tw-1)(\tz-2\ty+\tw) \cr
{\cal D}_7 &=
\tz\(\tx+\ty-\tw\)^2+3yz\(\tx+\tz+\ty+1\)\(\tz-\ty+\tw\)\(\ty-2\tz\) \cr
&
-y\tz\(\tz-2\ty+\tw-1\)\(\tz-2\ty+\tw\)-x\tz\(\tw-2\tx\)\(\tw-2\tx-1\) \cr
{\cal D}_8 &=
9\tx^2-18\tx\tw+25\tx\ty-41\tx\tz+16\tz\tw  \cr
&
-48yzw\(\tx+\tz+\ty+1\)\(\ty-2\tz\)
-9x\(\tw-2\tx\)\(\tw-2\tx-1\) \cr
&
+yw\(\tx-2\tz\)\(\tz-2\ty+\tw\)
+4xyw^2\(\ty+1\)\(\tx+\tz+\ty+1\)  \cr
&
+xw\(9\tz+10\ty+9\tw+9\)\(\tw-2\tx\) \cr
{\cal D}_9 &=
3\tx\ty-6\tx\tz-6\tw\ty+3\tw\ty+3\ty^2-9\ty\tz+13\tw\tz+3\tz^2 \cr
&
-3y(\tz-2\ty+\tw-1)(\tz-2\ty+\tw)-3z(2\tz-\ty+1)(2\tz-\ty)   \cr
&
-xw\tz(\tw-2\tx)
+yw(5\tz+3\tw+3)(\tz-2\ty+\tw) \cr
&
+xyw^2(\tx+\tz+\ty+1)(8\tz+6\tw+12)
\cr}
}

\vfill\eject

\leftline{\undertext{Phase B:} }
\eqn\pfB{
\eqalign{
{\cal D}_1 &=\tx\(\tx-\ty+\tz-\tw\)-x\(\tx+\ty+\tz+1\)\(\ty-\tx-\tw\) \cr
{\cal D}_2 &=\tx\(\ty-2\tz\)-xy\(\tx+\ty+\tz+2\)\(\tx+\ty+\tz+1\)  \cr
{\cal D}_3 &=\tw^2-w\(\ty-\tx-\tw\)\(\tx-\ty+\tz-\tw\) \cr
{\cal D}_4 &=\(\ty-\tx-\tw\)\(\ty-2\tz\)
             -y\(\tx+\ty+\tz+1\)\(\tx-\ty+\tz-\tw\)\cr
{\cal D}_5 &= 9\tx\tw-2\tx\ty-16\tz\tw-16\tx\tz+16\ty\tz
-48yz\(\tx+\ty+\tz+1\)\(\ty-2\tz\) \cr
&
-16y\tz\(\tx-\ty+\tz-\tw\)-9xw\(\ty-\tx-\tw-1\)\(\ty-\tx-\tw\)\cr
&
-4xy\(2\tx-3\ty+2\tw-1\)\(\tx+\ty+\tz+1\)-x\(9\tw-10\ty\)\(\ty-\tx-\tw\) \cr
{\cal D}_6 &=
3\tx\tw+8\ty\tw+6\tx\ty-24\tz\tw-16\tx\tz+8\tz^2
-8z\(\ty-2\tz-1\)\(\ty-2\tz\)  \cr
&
-16yz\(\tx+\ty+\tz+1\)\(\ty-2\tz\)-8yw\(\tx-\ty+\tz-\tw-1\)\(\tx-\ty+\tz-\tw\)
\cr &
-3xw\(\ty-\tx-\tw-1\)\(\ty-\tx-\tw\)-8y\tw\(\tx-\ty+\tz-\tw\) \cr
&
-4xy\(5\tw+2\tx+\ty+3\)\(\tx+\ty+\tz+1\)-x\(3\tw-2\ty\)\(\ty-\tx-\tw\) \cr}}

\vfill\eject

\appendix{B}{Basis of the Mori cone for $\IP(7,3,2,1,1)$}
For this weighted projective space, we have two different desingularization
of the ambient space, $\IP_{\SDs_{reg}^A}$ and $\IP_{\SDs_{reg}^B}$ in the
text. For each desingularization, we obtain the basis of the Mori cone
following\OdaPark. We see the Mori cone for $\SDs^B_{reg}$ is not simplicial.

\noindent
For the regular fan $\SDs^A_{reg}$;

\eqn\TAbase{
\eqalign{
&\eta_A^1=(\;1,0,0,1,0,0, -2, -1,1,\;0,\;0,\;0) ,\; \cr
&\eta_A^3=( -2,0,1,0,1,1,\;0,\;1,0, -2,\;0,\;0)  ,\;  \cr
&\eta_A^5=( -2,1,0,0,1,1,\;1,\;0,0,\;0, -2,\;0) ,\; \cr
&\eta_A^7=( -2,1,1,0,1,1,\;0,\;0,0,\;0,\;0,-2) ,\;  \cr}
\eqalign{
&\eta_A^2=( -2,0,\;0,0,\;1,\;1,\;1,\;1, -2,0,0,0)  , \cr
&\eta_A^4=(\;2,0,\;0,0, -1, -1,\;0, -2,\;1,1,0,0) ,  \cr
&\eta_A^6=(\;1,0,\;0,0, -1, -1, -1,\;0,\;1,0,1,0) , \cr
&\eta_A^8=(\;1,0, -1,0, -1, -1,\;0,\;0,\;0,1,0,1). \cr}
}
For the regular fan $\SDs^B_{reg}$;
\eqn\TBbase{
\eqalign{
&\eta_B^1=(\;0,0,\;0,1,\;0,\;0, -2,0,\;0,\;0,\;1,0,0) ,\; \cr
&\eta_B^3=( -1,0,\;0,0,\;1,\;1,\;0,0, -2,\;0,\;0,0,1) ,\; \cr
&\eta_B^5=( -2,0,\;1,0,\;1,\;1,\;0,1,\;0, -2,\;0,0,0) ,\;\cr
&\eta_B^7=( -2,1,\;0,0,\;1,\;1,\;1,0,\;0,\;0, -2,0,0) ,\;\cr
&\eta_B^9=(\;1,0, -1,0, -1, -1,\;0,0,\;0,\;1,\;0,1,0) . \cr}
\eqalign{
&\eta_B^2=( -1,0,0,\;1,\;0,\;0,0,\;1,1,0,0,\;0, -2) , \cr
&\eta_B^4=(\;1,0,0, -1, -1, -1,1,\;0,1,0,0,\;0,\;0) , \cr
&\eta_B^6=(\;2,0,0,\;0, -1, -1,0, -2,1,1,0,\;0,\;0) , \cr
&\eta_B^8=( -2,1,1,\;0,\;1,\;1,0,\;0,0,0,0, -2,\;0) , \cr
&  \cr}
}

\vskip1cm

\appendix{C}{Topological data for models with $h^{1,1}\leq 3$}

We list the topological couplings for the Calabi-Yau models with
$h^{1,1}\leq 3$. We follow the conventions in \HKTYI\HKTYII, i.e., $8
J_1^3+4 J_1^2 J_2$ for the coupling means $K^{cl}_{x_1x_1x_1}=8,
K^{cl}_{x_1x_1x_2}=4$ and others are zero. The superscript in each model
shows the type of the model defined in \models. The divisors $J_k$ and the
variables $x^{(k)}=(-1)^{l_0^{(k)}}a^{l^{(k)}}$ are connected by the
identification $J_k=m(\theta_{x^{(k)}})$ made in \coupling\ and \yukcp.
According to Wall's theorem cited in sect.4, the topological type of the
Calabi-Yau manifolds are classified by the classical Yukawa couplings (
cubic form) and the invariant $c_2\cdot J_k$ (linear form) on
$H^{1,1}(X,\ZZ)$.

For interested reader we list the concrete basis $\{ l^{(k)} \}$ for
the Mori cone in the file appended to \HLY. The basis for the Mori cone
and the topological couplings in this list determine the prepotential
$F(t)$ in \PotF.

\vskip0.3cm
\def\cI{
6\,{{J_{1}}^3} + 9\,{{J_{1}}^2}\,J_{2} +
   13\,J_{1}\,{{J_{2}}^2} + 17\,{{J_{2}}^3}
}
\def\cII{
14\,{{J_{1}}^3} + 7\,{{J_{1}}^2}\,J_{2} +
   3\,J_{1}\,{{J_{2}}^2}
}
\def\cIII{
36\,{{J_{1}}^3} + 12\,{{J_{1}}^2}\,J_{2} +
   4\,J_{1}\,{{J_{2}}^2} + {{J_{2}}^3}
}
\def\cIV{
63\,{{J_{1}}^3} + 21\,{{J_{1}}^2}\,J_{2} +
   7\,J_{1}\,{{J_{2}}^2} + 2\,{{J_{2}}^3}
}
\def\cV{
 9\,{{J_{1}}^3} + 3\,{{J_{1}}^2}\,J_{2} + \,J_{1}\,{{J_{2}}^2}
}
\def\cVIa{
 8\,{{J_{1}}^3} + 14\,{{J_{1}}^2}\,J_{2} + 24\,J_{1}\,{{J_{2}}^2} +
   37\,{{J_{2}}^3}
}
\def\cVIb{
+ 4\,{{J_{1}}^2}\,J_{3} + 7\,J_{1}\,J_{2}\,J_{3} +
   10\,{{J_{2}}^2}\,J_{3}
}
\def\cVIc{
+ 2\,J_{1}\,{{J_{3}}^2} + 2\,J_{2}\,{{J_{3}}^2}
}
\def\cVIIa{
90\,{{J_{1}}^3} + 30\,{{J_{1}}^2}\,J_{2} +
   10\,J_{1}\,{{J_{2}}^2} + 3\,{{J_{2}}^3}
}
\def\cVIIb{
+ 45\,{{J_{1}}^2}\,J_{3} + 15\,J_{1}\,J_{2}\,J_{3} + 5\,{{J_{2}}^2}\,J_{3}
}
\def\cVIIc{
+ 15\,J_{1}\,{{J_{3}}^2} + 5\,J_{2}\,{{J_{3}}^2} + 5\,{{J_{3}}^3}
}
\def\cVIIIa{
15\,{{J_{1}}^3} + 20\,{{J_{1}}^2}\,J_{2} +
   26\,J_{1}\,{{J_{2}}^2}
}
\def\cVIIIb{
+ 32\,{{J_{2}}^3} + 10\,{{J_{1}}^2}\,J_{3} + 13\,J_{1}\,J_{2}\,J_{3}
}
\def\cVIIIc{
+ 16\,{{J_{2}}^2}\,J_{3} + 6\,J_{1}\,{{J_{3}}^2} + 6\,J_{2}\,{{J_{3}}^2}
}
\def\cIXa{
18\,{{J_{1}}^3} + 12\,{{J_{1}}^2}\,J_{2} + 8\,J_{1}\,{{J_{2}}^2} +
   5\,{{J_{2}}^3}
}
\def\cIXb{
+ 9\,{{J_{1}}^2}\,J_{3} + 6\,J_{1}\,J_{2}\,J_{3} +
   4\,{{J_{2}}^2}\,J_{3}
}
\def\cIXc{
+ 3\,J_{1}\,{{J_{3}}^2} + 2\,J_{2}\,{{J_{3}}^2} +
   {{J_{3}}^3}
}
\def\cXa{
40\,{{J_{1}}^3} + 20\,{{J_{1}}^2}\,J_{2} + 10\,J_{1}\,{{J_{2}}^2}
+ 4\,{{J_{2}}^3}
}
\def\cXb{
+ 10\,{{J_{1}}^2}\,J_{3} + 5\,J_{1}\,J_{2}\,J_{3} + 2\,{{J_{2}}^2}\,J_{3}
}
\def\cXc{  }
\def\cXIa{
36\,{{J_{1}}^3} + 12\,{{J_{1}}^2}\,J_{2} + 4\,J_{1}\,{{J_{2}}^2} +
   {{J_{2}}^3}
}
\def\cXIb{
 + 18\,{{J_{1}}^2}\,J_{3} + 6\,J_{1}\,J_{2}\,J_{3} +
   2\,{{J_{2}}^2}\,J_{3}
}
\def\cXIc{
+ 6\,J_{1}\,{{J_{3}}^2} + 2\,J_{2}\,{{J_{3}}^2} +
   2\,{{J_{3}}^3}
}
\def\cXIIa{
50\,{{J_{1}}^3} + 30\,{{J_{1}}^2}\,J_{2} +
   18\,J_{1}\,{{J_{2}}^2} + 9\,{{J_{2}}^3}
}
\def\cXIIb{
+ 60\,{{J_{1}}^2}\,J_{3} + 36\,J_{1}\,J_{2}\,J_{3} + 21\,{{J_{2}}^2}\,J_{3}
}
\def\cXIIc{
+ 72\,J_{1}\,{{J_{3}}^2} + 43\,J_{2}\,{{J_{3}}^2} + 86\,{{J_{3}}^3}
}
\def\cXIIIa{
 8\,{{J_{1}}^3} + 18\,{{J_{1}}^2}\,J_{2} + 36\,J_{1}\,{{J_{2}}^2}
 + 72\,{{J_{2}}^3}
}
\def\cXIIIb{
 + 4\,{{J_{1}}^2}\,J_{3} + 9\,J_{1}\,J_{2}\,J_{3}+ 18\,{{J_{2}}^2}\,J_{3}
}
\def\cXIIIc{
+ 2\,J_{1}\,{{J_{3}}^2} + 4\,J_{2}\,{{J_{3}}^2}
}
\def\cXIVa{
 72\,{{J_{1}}^3} + 18\,{{J_{1}}^2}\,J_{2} + 4\,J_{1}\,{{J_{2}}^2}
}
\def\cXIVb{
+ 36\,{{J_{1}}^2}\,J_{3} + 9\,J_{1}\,J_{2}\,J_{3} + 2\,{{J_{2}}^2}\,J_{3}
}
\def\cXIVc{
+ 18\,J_{1}\,{{J_{3}}^2} + 4\,J_{2}\,{{J_{3}}^2} + 8\,{{J_{3}}^3}
}
\def\cXVa{
6\,{{J_{1}}^3} + 16\,{{J_{1}}^2}\,J_{2} + 42\,J_{1}\,{{J_{2}}^2} +
   104\,{{J_{2}}^3}
}
\def\cXVb{
+ 2\,{{J_{1}}^2}\,J_{3} + 5\,J_{1}\,J_{2}\,J_{3} +
   10\,{{J_{2}}^2}\,J_{3}
}

\def\cXVIa{
50\,{{J_{1}}^3} + 10\,{{J_{1}}^2}\,J_{2} + 2\,J_{1}\,{{J_{2}}^2}
}
\def\cXVIb{
+ 80\,{{J_{1}}^2}\,J_{3} + 16\,J_{1}\,J_{2}\,J_{3} + 3\,{{J_{2}}^2}\,J_{3}
}
\def\cXVIc{
+ 128\,J_{1}\,{{J_{3}}^2} + 25\,J_{2}\,{{J_{3}}^2} + 203\,{{J_{3}}^3}
}
\def\fI{
8\,{{J_{1}}^3} + 4\,{{J_{1}}^2}\,J_{2}
}
\def\fII{
 4\,{{J_{1}}^3} + 2\,{{J_{1}}^2}\,J_{2}
}
\def\fIII{
 2\,{{J_{1}}^3} + 3\,{{J_{1}}^2}\,J_{2} + 3\,J_{1}\,{{J_{2}}^2} +
   3\,{{J_{2}}^3}
}
\def\fIV{
 2\,{{J_{1}}^3} + 7\,{{J_{1}}^2}\,J_{2} + 21\,J_{1}\,{{J_{2}}^2} +
   63\,{{J_{2}}^3}
}
\def\fV{
9\,{{J_{1}}^3} + 3\,{{J_{1}}^2}\,J_{2} + J_{1}\,{{J_{2}}^2}
}
\def\fVIa{
 18\,{{J_{1}}^3} + 6\,{{J_{1}}^2}\,J_{2} + 2\,J_{1}\,{{J_{2}}^2} }
\def\fVIb{
+ 18\,{{J_{1}}^2}\,J_{3} + 6\,J_{1}\,J_{2}\,J_{3} + {{J_{2}}^2}\,J_{3}
}
\def\fVIc{
+ 18\,J_{1}\,{{J_{3}}^2} + 3\,J_{2}\,{{J_{3}}^2} + 9\,{{J_{3}}^3}
}
\def\fVIIa{
6\,{{J_{1}}^3} + 4\,{{J_{1}}^2}\,J_{2} + 8\,{{J_{1}}^2}\,J_{3}
+ 4\,J_{1}\,J_{2}\,J_{3}
}
\def\fVIIb{
+ 8\,J_{1}\,{{J_{3}}^2} + 4\,J_{2}\,{{J_{3}}^2} +    8\,{{J_{3}}^3}
}

\def\fVIIIa{
  3\,{{J_{1}}^3} + 5\,{{J_{1}}^2}\,J_{2} + 5\,J_{1}\,{{J_{2}}^2} +
   5\,{{J_{2}}^3} }
\def\fVIIIb{
  + 10\,{{J_{1}}^2}\,J_{3} + 15\,J_{1}\,J_{2}\,J_{3} +
    15\,{{J_{2}}^2}\,J_{3}}
\def\fVIIIc{
  + 30\,J_{1}\,{{J_{3}}^2} + 45\,J_{2}\,{{J_{3}}^2} + 90\,{{J_{3}}^3}
}
\def\fIXa{
3\,{{J_{1}}^3} + 2\,{{J_{1}}^2}\,J_{2} + 4\,{{J_{1}}^2}\,J_{3}
+ 2\,J_{1}\,J_{2}\,J_{3}
}
\def\fIXb{
+ 4\,J_{1}\,{{J_{3}}^2} + 2\,J_{2}\,{{J_{3}}^2} + 4\,{{J_{3}}^3}
}

\def\fXa{
 8\,{{J_{1}}^3} + 2\,{{J_{1}}^2}\,J_{2} + 4\,{{J_{1}}^2}\,J_{3}}
\def\fXb{
+  J_{1}\,J_{2}\,J_{3} + 2\,J_{1}\,{{J_{3}}^2}
}

\def\X#1#2(#3)#4#5{ {$X_{#1}^#2(#3)_{#5}^{#4}$} }

\def\I{{\rm I}}
\def\II{{{\rm I}\hskip-1.5pt{\rm I}} }
\def\III{{{\rm I}\hskip-1.5pt{\rm I}\hskip-1.5pt{\rm I}} }


\vbox{\tabskip=0pt \offinterlineskip
\def\tabrule{\noalign{\hrule}}
\def\vspace#1{\omit& height #1 & \omit && \omit && \omit &\cr}
\def\vtabrule{\vspace{0.2cm}\tabrule\vspace{0.2cm}}
\def\second#1{ && \omit && $ #1 $ && \omit &\cr}
\halign
{\strut#& \vrule# \tabskip=0.5em  &
 \hfil# & \vrule# & \hfil#  & \vrule# &
 \hfil# & \vrule#
 \tabskip=0pt \cr\tabrule
&& \multispan5\hfil
   Fermat type Calabl-Yau hypersurfaces \hfil&\cr\tabrule
&& \omit\hidewidth model     \hidewidth  &&
   \omit\hidewidth topological couplings    \hidewidth &&
   \omit\hidewidth $c_2\cdot \vec J $ \hidewidth  &\cr\tabrule
\vspace{0.2cm}
&& \X{8}{\I}(2,2,2,1,1){2}{-168} && $\fI$ \hfill && $(56,24)$ &\cr\vtabrule
&& \X{12}{\II}(6,2,2,1,1){2}{-252} && $\fII$ \hfill&& $(52,24)$   &\cr\vtabrule
&& \X{12}{\III}(4,3,2,2,1){2}{-144} && $\fIII$ \hfill && $(32,42)$
&\cr\vtabrule
&& \X{14}{\II}(7,2,2,2,1){2}{-240} &&  $\fIV$ \hfill &&  $(44,126)$
&\cr\vtabrule
&& \X{18}{\II}(9,6,1,1,1){2}{-540} && $\fV$ \hfill && $(102,36)$ &\cr\vtabrule
&& \X{12}{\II}(6,3,1,1,1){3}{-344} && $\fVIa$  \hfill && $(96,36,102)$ &\cr
                 \second{\fVIb}
                 \second{\fVIc} \vtabrule
&& \X{12}{\III}(3,3,3,2,1){3}{-126} && $\fVIIa$  \hfill && $(48,24,56)$ &\cr
                 \second{\fVIIb}  \vtabrule
&& \X{15}{\III}(5,3,3,3,1){3}{-144} && $\fVIIIa$  \hfill && $(42,50,120)$ &\cr
                 \second{\fVIIIb}
                 \second{\fVIIIc} \vtabrule
&& \X{18}{\III}(9,3,3,2,1){3}{-186} && $\fIXa$ \hfill && $(42,24,52)$ &\cr
                 \second{\fIXb}   \vtabrule
&& \X{24}{\II}(12,8,2,1,1){3}{-480} && $\fXa$  \hfill && $(92,24,48)$ &\cr
                 \second{\fXb}   \tabrule
\noalign{\smallskip} } }

\vfill\eject


\vbox{\tabskip=0pt \offinterlineskip
\def\tabrule{\noalign{\hrule}}
\def\vspace#1{\omit& height #1 & \omit && \omit && \omit &\cr}
\def\vtabrule{\vspace{0.2cm}\tabrule\vspace{0.2cm}}
\def\second#1{ && \omit && $ #1 $ && \omit &\cr}
\halign
{\strut#& \vrule# \tabskip=0.5em  &
 \hfil# & \vrule# & \hfil#  & \vrule# &
 \hfil# & \vrule#
 \tabskip=0pt \cr\tabrule
&& \multispan5\hfil
   Non-Fermat type Calabl-Yau hypersurfaces \hfil&\cr\tabrule
&& \omit\hidewidth model     \hidewidth  &&
   \omit\hidewidth topological couplings    \hidewidth &&
   \omit\hidewidth $c_2\cdot \vec J $ \hidewidth &\cr\tabrule
\vspace{0.2cm}
&& \X{9}{\III}(3,2,2,1,1){2}{-168} && $\cI$ \hfill && $(48,74)$ &\cr\vtabrule
&& \X{7}{\I}(2,2,1,1,1){2}{-186} && $\cII$ \hfill&& $(68,36)$  &\cr\vtabrule
&& \X{8}{\III}(3,2,1,1,1){2}{-202} && $\cIII$ \hfill && $(96,34)$ &\cr\vtabrule
&& \X{8}{\I}(3,1,1,1,1){2}{-240} &&  $\cIV$ \hfill &&  $(126,44)$ &\cr\vtabrule
&& \X{14}{\III}(7,3,2,1,1){2}{-260} && $\cV$ \hfill && $(66,24)$  &\cr\vtabrule
&& \X{15}{\III}(5,4,3,2,1){3}{-126} && $\cVIa$  \hfill && $(44,82,24)$ &\cr
                 \second{\cVIb}
                 \second{\cVIc} \vtabrule
&& \X{10}{\III}(3,2,2,2,1){3}{-144} && $\cVIIa$  \hfill && $(120,42,50)$ &\cr
                 \second{\cVIIb}
                 \second{\cVIIc} \vtabrule
&& \X{10}{\III}(3,3,2,1,1){3}{-168} && $\cVIIIa$  \hfill && $(66,92,48)$ &\cr
                 \second{\cVIIIb}
                 \second{\cVIIIc} \vtabrule
&& \X{20}{\III}(10,4,3,2,1){3}{-192} && $\cIXa$ \hfill && $(72,50,34)$ &\cr
                 \second{\cIXb}
                 \second{\cIXc} \vtabrule
&& \X{10}{\I}(4,2,2,1,1){3}{-192} && $\cXa$  \hfill && $(100,52,24)$ &\cr
                 \second{\cXb} \vtabrule
&& \X{16}{\III}(8,3,2,2,1){3}{-200} && $\cXIa$  \hfill && $(96,34,44)$ &\cr
                 \second{\cXIb}
                 \second{\cXIc} \vtabrule
&& \X{12}{\III}(5,3,2,1,1){3}{-204}&& $\cXIIa$ \hfill && $(104,66,128)$ &\cr
                 \second{\cXIIb}
                 \second{\cXIIc} \tabrule
\noalign{\smallskip} } }

\vfill\eject

\vbox{\tabskip=0pt \offinterlineskip
\def\tabrule{\noalign{\hrule}}
\def\vspace#1{\omit& height #1 & \omit && \omit && \omit &\cr}
\def\vtabrule{\vspace{0.2cm}\tabrule\vspace{0.2cm}}
\def\second#1{ && \omit && $ #1 $ && \omit &\cr}
\halign
{\strut#& \vrule# \tabskip=0.5em  &
 \hfil# & \vrule# & \hfil#  & \vrule# &
 \hfil# & \vrule#
 \tabskip=0pt \cr
\multispan7 table cont'd   \hfill \cr
\multispan7   \hfill  \cr\tabrule
\vspace{0.2cm}
&& \X{18}{\II}(9,4,2,2,1){3}{-240}&& $\cXIIIa$ \hfill && $(68,132,36)$ &\cr
                 \second{\cXIIIb}
                 \second{\cXIIIc} \vtabrule
&& \X{9}{\I}(4,2,1,1,1){3}{-240} && $\cXIVa$ \hfill && $(132,36,68)$ &\cr
                 \second{\cXIVb}
                 \second{\cXIVc} \vtabrule
&& \X{16}{\II}(8,3,3,1,1){3}{-256} && $\cXVa$  \hfill && $(60,164,24)$ &\cr
                 \second{\cXVb}
                 \second{\cXc} \vtabrule
&& \X{16}{\II}(8,5,1,1,1){3}{-456}&& $\cXVIa$ \hfill&& $(164,36,266)$ &\cr
                 \second{\cXVIb}
                 \second{\cXVIc} \tabrule
\noalign{\smallskip} } }


\listrefs

\bye

\input harvmac
\overfullrule=0pt
\baselineskip=12pt
\def\undertext#1{$\underline{\smash{\hbox{#1}}}$}
\def\X#1#2(#3)#4#5{ {$X_{#1}^{#2}(#3)^{#4}_{#5}$} }
\def\model#1{\noindent\undertext{#1}}
\def\I{{\rm I}}
\def\II{{{\rm I}\hskip-1.5pt{\rm I}} }
\def\III{{{\rm I}\hskip-1.5pt{\rm I}\hskip-1.5pt{\rm I}} }
\def\({ \left(  }
\def\){ \right) }

\def\IP{{\bf P}}
\def\cD#1{{\cal D}_{#1}}
\def\z#1{z_{#1}}
\def\a#1{a_{#1}}
\def\ls#1{l^{(#1)}}
\def\tx{\theta_x}
\def\ty{\theta_y}
\def\tz{\theta_z}
\def\tw{\theta_w}
\def\bx{\bar x}
\def\by{\bar y}
\def\bz{\bar z}

\centerline{\bf Picard-Fuchs differential operators
for Calabi-Yau hypersurfaces with $h^{1,1}\leq 3$
\footnote{$^\dagger$}
{Comments and/or questions should be sent to hosono@sci.toyama-u.ac.jp}}
\centerline{\bf --- Optional appendix
to "GKZ-Generalized Hypergeometric Systems }
\centerline{\bf in Mirror Symmetry of Calabi-Yau
Hypersurfaces"(alg-geom/9511001)}
\rightline{\bf by S.Hosono, B.Lian and S.-T.Yau}

\vskip0.5cm

In this appendix, we will list a complete set of the Picard-Fuchs
differential operators for the large radius limit. From a consideration
of the intersection ring, we can see that we should have two differential
operators, one is second order and the other is
third order, for $h^{1,1}=2$. For the models with $h^{1,1}=3$, we
have five differential operators in general, three of them are second
order and two of them are third order.

The listing also contains the form of the defining equation in the weighted
projective space ( and its transposed potential for the non-Fermat models),
the basis of the Mori cone and the irreducible parts of the discriminants.
Since for most of the models with $h^{1,1}=3$, the form of the discriminants
becomes too complicated, we have not determined all of them in this list.

\vskip0.3cm

For the two moduli models ($h^{1,1}=2$), the Picard-Fuchs operators
have been obtained in {\it "Mirror Symmetry and Mirror Map and Applications to
Calabi-Yau Hypersurfaces"}(S.Hosono,A.Klemm,S.Theisen and S.-T.Yau, Commun.
Math.Phys.{\bf 167}(1995)30) and
{\it "Mirror Symmetry and the Moduli Space for Generic Hypersurfaces
in Toric Varieties"} (P.Berglund,A.Klemm and S.Katz,
hepth-9506091)). For completeness, we have included these models in our
list.

\vskip0.5cm

\leftline{ \undertext{{\bf Fermat models}} }
\vskip0.5cm

\model{\X{8}{\I}(2,2,2,1,1){2}{-168}:}
$$
W=\a1\z1^4+\a2\z2^4+\a3\z3^4+\a4\z4^8+\a5\z5^8
+\a0\z1\z2\z3\z4\z5+\a6\z4^4\z5^4
$$
$$
\ls1=(-4,1,1,1,0,0,1) \;\;,\;\;  \ls2=(0,0,0,0,1,1,-2)
$$
$$
\eqalign{
\cD1&=\tx^2(\tx-2\ty)-4x(4\tx+3)(4\tx+2)(4\tx+1) \cr
\cD2&=\ty^2-y(2\ty-\tx+1)(2\ty-\tx) \cr}
$$
$$
\bx=2^8x \quad,\quad \by=4y
$$
$$
dis_0=(1-\bx)^2-\bx^2\by \;\;,\;\;
dis_1=1-\by
$$
\vskip0.3cm

\model{\X{12}{\II}(6,2,2,1,1){2}{-252}:}
$$
W=\a1\z1^2+\a2\z2^6+\a3\z3^6+\a4\z4^{12}+\a5\z5^{12}
+\a0\z1\z2\z3\z4\z5+\a6\z4^6\z5^6
$$
$$
\ls1=(-6,3,1,1,0,0,1) \;\;,\;\;  \ls2=(0,0,0,0,1,1,-2)
$$
$$
\eqalign{
\cD1&=\tx^2(\tx-2\ty)-8x(6\tx+5)(6\tx+3)(6\tx+1) \cr
\cD2&=\ty^2-y(2\ty-\tx+1)(2\ty-\tx)  \cr}
$$
$$
\bx=2^63^3 x \;\;,\;\; \by=4y
$$
$$
dis_0=(1-\bx)^2-\bx^2\by  \;\;,\;\;
dis_1=1-\by
$$
\vskip0.3cm

\model{\X{12}{\III}(4,3,2,2,1){2}{-144}:}
$$
W=\a1\z1^3+\a2\z2^4+\a3\z3^6+\a4\z4^6+\a5\z5^{12}
+\a0\z1\z2\z3\z4\z5+\a6\z2^2\z5^6
$$
$$
\ls1=(-6,2,0,1,1,-1,3) \;\;,\;\;  \ls2=(0,0,1,0,0,1,-2)
$$
$$
\eqalign{
\cD1&=\tx^2(3\tx-2\ty)-36x(6\tx+5)(6\tx+1)\{\ty-\tx+2y(1+6\tx-2\ty)\} \cr
\cD2&=(\ty-\tx)\ty-y(3\tx-2\ty-1)(3\tx-2\ty) \cr}
$$
$$
\bx=2^3 3^3 x \;\;,\;\; \by=2^2 3 y
$$
$$
\eqalign{
& dis_0=1+2\bx-6\bx\by-9\bx^2\by+6\bx^2\by^2-\bx^2\by^3
\cr
& dis_1=3-\by \cr}
$$

\vskip0.3cm

\model{\X{14}{\II}(7,2,2,2,1){2}{-240}:}
$$
W=\a1\z1^2+\a2\z2^7+\a3\z3^7+\a4\z4^7+\a5\z5^{14}
+\a0\z1\z2\z3\z4\z5+6\z1\z5^7
$$
$$
\ls1=(-7,0,1,1,1,-3,7) \;\;,\;\;  \ls2=(0,1,0,0,0,1,-2)
$$
$$
\eqalign{
\cD1&=\tx^2(7\tx-2\ty)-7x\{y(28\tx-4\ty+18)+\ty-3\tx-2\} \cr
&\times\{y(28(\tx-4\ty+10)+\ty-3\tx-1\}\{y(28\tx-4\ty+2)\ty-3\tx\} \cr
\cD2&=(\ty-3\tx)\ty-y(7\tx-2\ty-1)(7\tx-2\ty) \cr}
$$
$$
\bx=x \;\;,\;\; \by=7y
$$
$$
\eqalign{
& dis_0=1+27\bx-63\bx\by+56\bx\by^2-112\bx\by^3-(7-4\by)^4\bx^2\by^3
\cr
& dis_1=7-4\by \cr}
$$
\vskip0.3cm

\model{\X{18}{\III}(9,6,1,1,1){2}{-540}:}
$$
W=\a1\z1^2+\a2\z2^3+\a3\z3^{18}+\a4\z4^{18}+\a5\z5^{18}
+\a0\z1\z2\z3\z4\z5+\a6\z3^6\z4^6\z5^6
$$
$$
\ls1=(-6,3,2,0,0,0,1) \;\;,\;\;  \ls2=(0,0,0,1,1,1,-3)
$$
$$
\eqalign{
\cD1&=\tx(\tx-3\ty)-12x(6\tx+5)(6\tx+1) \cr
\cD2&=\ty^3-y(\tx-3\ty-2)(\tx-3\ty-1)(\tx-3\ty) \cr}
$$
$$
\bx=2^43^3x \;\;,\;\; \by=3^3 y
$$
$$
dis_0=(1-\bx)^3-\bx^3\by \;\;,\;\; dis_1=1+\by
$$
\vskip0.3cm

\model{\X{12}{\II}(6,3,1,1,1){3}{-344}}
\footnote{\hskip-5pt$^\dagger$}
{This model has a contribution from the twisted sector.
Here we have included  the twisted sector by modifying the polyhedron
$\Delta(w)$. The incorporation of the twisted sector through the
modification of the polyhedron $\Delta(w)$ for this model
has been settled in {\it "Mirror Symmetry and the Moduli Space for
Generic Hypersurfaces in Toric Varieties"} by P.Berglund, A.Klemm and S.Katz
(hepth-9506091). Here we introduced the deformation "$\a7$"
following their modified polyhedron $\Delta(w)$. (We cannot express the
deformation parameter $\a7$ in the potential $W$ in $\IP(w)$,
however we can see it in the toric description $f^*(a)$ in the text.) }:
$$
W=\a1\z1^2+\a2\z2^4+\a3\z3^{12}+\a4\z4^{12}+\a5\z5^{12}
+\a0\z1\z2\z3\z4\z5 +\a6\z3^4\z4^4\z5^4
$$
$$
\eqalign{
& \ls1=(0,1,1,0,0,0,-2,0) \;\;,\;\;  \ls2=(0,0,0,1,1,1,0,-3) \cr
& \ls3=(-4,-1,0,0,0,0,4,1)   \cr}
$$
$$
\eqalign{
\cD1&=(\tx-\tz)\tx-x(2\tx-4\tz+1)(2\tx-4\tz) \cr
\cD2&=(\tx-2\tz)(3\ty-\tz)  \cr
& +2z(\tz+1)\{4x(8x+1)(\tx-2\tz)-(8x+1)(\tx-\tz)-6x(4\tx+1)\}\cr
\cD3&=(\tx-\tz)(\tz-3\ty)-z\{2x(2\tx-4\tz-3)-(\tx-\tz-1)\} \cr
  &\hskip2.5cm \times
      \{4x(8x+1)(\tx-2\tz)-(8x+1)(\tx-\tz)-6x(4\tz+1)\}\cr
\cD4&=\ty^3-y(\tz-3\ty-2)(\tz-3\ty-1)(\tz-3\ty) \cr
\cD5&=(\tx-2\tz)\ty^2 -2z\{3y(\tz-3\ty-1)(\tz-3\ty)+\ty^2\} \cr
  &\hskip2.5cm \times
      \{4x(8x+1)(\tx-2\tz)-(8x+1)(\tx-\tz)-6x(4\tz+1)\}\cr }
$$
\vskip0.3cm

\model{\X{12}{\III}(3,3,3,2,1){3}{-126}:}
$$
W=\a1\z1^4+\a2\z2^4+\a3\z3^4+\a4\z4^6+\a5\z5^{12}
+\a0\z1\z2\z3\z4\z5+\a6\z4^4\z5^4+\a7\z4^2\z5^8
$$
$$
\eqalign{
& \ls1=(-4,1,1,1,0,-1,0,2) \;\;, \ls2=(0,0,0,0,1,0,-2,1)\;\;  \cr
& \ls3=(0,0,0,0,0,1,1,-2)  \cr}
$$
$$
\eqalign{
\cD1&=\ty(\tz-\tx)-yz(\tz-2\ty)(2\tx-2\tz+\ty) \cr
\cD2&=\ty(2\tx-2\tz+\ty)-y(\tz-2\ty-1)(\tz-2\ty) \cr
\cD3&=(\tz-\tx)(\tz-2\ty)-z(2\tx-2\tz+\ty-1)(2\tx-2\tz+\ty) \cr
\cD4&=\tx^2(\tz-2\ty)-4xz(4\tx+3)(4\tx+2)(4\tx+1) \cr
\cD5&=\tx^2(2\tx-2\tz+\ty)  \cr
 &  -4x(4\tx+3)(4\tx+1)\{6yz(\tz-2\ty)+z(2\tx-2\tz+\ty)+2(\tz-\tx)\} \cr}
$$
$$
\bx=2^6 x \;\;,\;\; \by=y \;\;,\;\; \bz=z
$$
$$
\eqalign{
&
dis_0=(1 + \bx - 4\bx\bz)(1 - 4\bx\bz)^2
    - 8\bx^2\bz^2(9 + 8\bx - 36\bx\bz)\by
    -432\bx^3\bz^4\by^2   \cr
&
dis_1 =1 - 4\by - 4\bz + 18\by\bz - 27\by^2\bz^2  \cr}
$$

\vskip0.3cm

\model{\X{15}{\III}(5,3,3,3,1){3}{-144}:}
$$
W=\a1\z1^3+\a2\z2^5+\a3\z3^5+\a4\z4^5+\a5\z5^{15}
+\a0\z1\z2\z3\z4\z5+\a6\z1\z5^{10}+\a7\z1^2\z5^5
$$
$$
\eqalign{
& \ls1=(-5,0,1,1,1,-3,5,0) \;\;,\;\;  \ls2=(0,1,0,0,0,0,1,-2,0)  \cr
& \ls3=(0,0,0,0,0,1,-2,1)  \cr}
$$
$$
\eqalign{
\cD1&=(3\tx-\tz)(2\ty-\tz)-z(5\tx+\ty-2\tz-1)(5\tx+\ty-2\tz) \cr
\cD2&=\ty(5\tx+\ty-2\tz)-y(2\ty-\tz+1)(2\ty-\tz) \cr
\cD3&=\ty(3\tx-\tz)-yz(5\tx+\ty-2\tz)(2\ty-\tz)  \cr
\cD4&=\tx^2(\tz-2\ty)-15xyz^3(5\tx+4)(5\tx+2)(5\tx+1) \cr
  &+10xz^2(5\tx+4)(5\tx+2)\{(3-10z)\tx+(6yz-2z)\ty-(1+3yz-4z)\tz\} \cr
  &-5xz(5\tx+4)\{(3-10z)\tx+(6yz-2z)\ty-(1+3yz-4z)\tz+1-2z\} \cr
  & \hskip2cm  \times
               \{(3-10z)\tx+(6yz-2z)\ty-(1+3yz-4z)\tz \} \cr
\cD5&=\tx^2(5\tx+\ty-2\tz) \cr
  &+15xyz^2(5\tx+3)(5\tx+2)\{(3-10z)\tx+(6yz-2z)\ty-(1+3yz-4z)\tz\} \cr
  &-10xz(5\tx+3)\{(3-10z)\tx+(6yz-2z)\ty-(1+3yz-4z)\tz+1-2z \} \cr
  & \hskip2cm \times\{(3-10z)\tx+(6yz-2z)\ty-(1+3yz-4z)\tz\} \cr
  &+5x\{(3-10z)\tx+(6yz-2z)\ty-(1+3yz-4z)\tz+2-4z\} \cr
  & \quad \times \{(3-10z)\tx+(6yz-2z)\ty-(1+3yz-4z)\tz+1-2z\} \cr
  & \quad \times \{(3-10z)\tx+(6yz-2z)\ty-(1+3yz-4z)\tz\} \cr}
$$
\vskip0.3cm

\model{\X{18}{\III}(9,3,3,2,1){3}{-186}:}
$$
W=\a1\z1^2+\a2\z2^6+\a3\z3^6+\a4\z4^9+\a5\z5^{18}
+\a0\z1\z2\z3\z4\z5+\a6\z4^3\z5^{12}+\a7\z4^6\z5^6
$$
$$
\eqalign{
& \ls1=(-6,3,1,1,0,-1,2,0) \;\;,\;\;  \ls2=(0,0,0,0,1,0,1,-2)  \cr
& \ls3=(0,0,0,0,0,1,-2,1)  \cr}
$$
$$
\eqalign{
\cD1&=(\tz-\tx)(\tz-2\ty)-z(2\tx+\ty-2\tz-1)(2\tx+\ty-2\tz) \cr
\cD2&=\ty(2\tx+\ty-2\tz)-y(2\ty-\tz-1)(2\ty-\tz)  \cr
\cD3&=\ty(\tz-\tx)+yz(2\tx+\ty-2\tz)(2\ty-\tz)    \cr
\cD4&=\tx^2(2\tx+\ty-2\tz)  \cr
& -24x(6\tx+5)(6\tx+1)\{(4z-1)\tx+(2z-6yz)\ty+(3yz-4z+1)\tz \} \cr
\cD5&=\tx^2(\tz-2\ty)-8xz(6\tx+5)(6\tx+3)(6\tx+1) \cr}
$$
$$
\bx=2^4 3^3 x \;\;,\;\; \by=y \;\;,\;\; \bz=z
$$
$$
\eqalign{
&
dis_0= (1 + \bx - 4\bx\bz)(1 - 4\bx\bz)^2
      - 8\bx^2\bz^2(9 + 8\bx - 36\bx\bz)\by
      - 432\bx^3\bz^4\by^2                     \cr
&
dis_= (1 - 4\by - 4\bz + 18\by\bz - 27\by^2\bz^2) \cr}
$$

\vskip0.3cm

\model{\X{24}{\II}(12,8,2,1,1){3}{-480}:}
$$
W=\a1\z1^2+\a2\z2^3+\a3\z3^{12}+\a4\z4^{24}+\a5\z5^{24}
+\a0\z1\z2\z3\z4\z5+\a6\z3^6\z4^6\z5^6+\a7\z4^{12}\z5^{12}
$$
$$
\eqalign{
& \ls1=(-6,3,2,0,0,0,1,0) \;\;,\;\;  \ls2=(0,0,0,0,1,1,0,-2) \cr
& \ls3=(0,0,0,1,0,0,-2,1)  \cr}
$$
$$
\eqalign{
\cD1&=\tx(\tx-2\tz)-12x(6\tx+5)(6\tx+1) \cr
\cD2&=\ty^2-y(2\ty-\tz+1)(2\ty-\tz) \cr
\cD3&=\tz(\tz-2\ty)-z(2\tz-\tx+1)(2\tz-\tx) \cr}
$$
$$
\bx=2^43^3x\;\;,\;\; \by=2^2y \;\;,\;\; \bz=2^2z
$$
$$
\eqalign{
& dis_0=(1-\bx)^4-2\bx^2\bz(1-\bx)^2+\bx^4\bz^2(1-\by) \cr
& dis_1=(1-2\bz+\bz^2-\by\bz^2) \;\;,\;\;
  dis_2=1-\by \cr}
$$

\vskip1cm

\leftline{ \undertext{{\bf Non-Fermat models}} }
\vskip0.5cm

\model{\X{9}{\III}(3,2,2,1,1){2}{-168}:}
$$
\eqalign{ &
W=\z1^3+\z2^4\z4+\z3^4\z5+\z4^9+\z5^9
\cr&
\hat W=\a1\z1^3+\a2\z2^4+\a3\z3^4+\a4\z2\z4^9+\a5\z3\z5^9
+\a0\z1\z2\z3\z4\z5+\a6\z1\z4^4\z5^4
\cr}
$$
$$
\ls1=(-3,0,1,1,-1,-1,3) \;\;,\;\;  \ls2=(-1,1,0,0,1,1,-2)
$$
$$
\eqalign{
\cD1&=
\ty\(\ty-\tx\)^2-y\(3\tx+\ty+1\)\(3\tx-2\ty-1\)\(3\tx-2\ty\)
\cr
\cD2&=
\(\tx-\ty\)^2+\(\tx-\ty\)\(3\tx-2\ty\)+4\tx\(3\tx-2\ty\) \cr
&
-48xy\(3\tx-2\ty-1\)\(3\tx+\ty+1\)-3y\(3\tx-2\ty-1\)\(3\tx-2\ty\) \cr
&
-48xy\(3\tx+\ty+3\)\(3\tx+\ty+1\)-16x\(\tx-\ty\)^2   \cr}
$$
$$
\eqalign{
dis_0=
&
(1 - 4y)^2 - (1 + 531y - 4335y^2 + 8748y^3)x           \cr
&
+(512 + 74496y - 1539648y^2 + 8503056y^3 - 14348907y^4)x^2y
-65536x^3y^2   \cr
}
$$

\vskip0.3cm

\model{\X{7}{\I}(2,2,1,1,1){2}{-186}:}
$$
\eqalign{ &
W=\z1^3\z4+\z2^3\z5+\z3^7+\z4^7+\z5^7
\cr&
\hat W=\a1\z1^3+\a2\z2^3+\a3\z3^7+\z1\z4^7+\z2\z5^7
+\a0\z1\z2\z3\z4\z5+\a6\z3^3\z4^3\z5^3
\cr}
$$
$$
\ls1=(-3,1,1,0,0,0,1) \;\;,\;\;  \ls2=(-1,0,0,1,1,1,-2)
$$
$$
\eqalign{
\cD1&=
\tx^2\(\tx-2\ty\)-x\(3\tx+\ty+3\)\(3\tx+\ty+2\)\(3\tx+\ty+1\)  \cr
\cD2&=
9\tx^2-21\tx\ty+7\ty^2 \cr
&
-27x(3\tx+\ty+2)(3\tx+\ty+1)-7y(\tx-2\ty-1)(\tx-2\ty) \cr
}
$$
$$
\eqalign{
dis_0  =
&(1 - 27x)^3 - (8 - 675x + 71442x^2)y    \cr
&
-(16 - 1372x + 453789x^2)y^2 - 823543x^2y^3
}
$$

\vskip0.3cm

\model{\X{8}{\III}(3,2,1,1,1){2}{-202}:}
$$
\eqalign{ &
W=\z1^2\z2+\z2^4+\z3^8+\z4^8+\z5^8
\cr&
\hat W=\a1\z1^2+\a2\z1\z2^4+\a3\z3^8+\a4\z4^8+\a5\z5^8
+\a0\z1\z2\z3\z4\z5+\a6\z2^2\z3^2\z4^2\z5^2
\cr}
$$
$$
\ls1=(-2,1,0,0,0,0,1)  \;\;,\;\;  \ls2=(-2,0,2,1,1,1,-3)
$$
$$
\eqalign{
\cD1&=\tx(\tx-3\ty)-x(2\tx+2\ty+2)(2\tx+2\ty+1)  \cr
\cD2&=(1-4x)^2\{ 4\ty^3-\ty^2\tx+2x\ty^2(2\tx+2\ty+1) \cr
&\hskip2cm      -4y(\tx-3\ty-2)(\tx-3\ty-1)(\tx-3\ty) \}  \cr
&\hskip0.5cm -16(1-4x)xy(\tx-3\ty-1)(\ty+5\tx+1)(\tx-3\ty) \cr
&\hskip0.5cm -64x^2y(\ty+5\tx+1)(\tx-3\ty) \cr
}
$$
$$
\bx=4x \;\;,\;\; \by=y
$$
$$
dis_0= ( 1 - \bx )^5 +
   ( 27 - 144\bx + 320\bx^2 - 2816\bx^3 - 512\bx^4 ) \by - 65536\bx^3 \by^2
$$
\vskip0.3cm

\model{\X{7}{\I}(3,1,1,1,1){2}{-240}:}
$$
\eqalign{ &
W=\z1^2\z5+\z2^7+\z3^7+\z4^7+\z5^7
\cr&
\hat W=\a1\z1^2+\a2\z2^7+\a3\z3^7+\a4\z4^7+\a5\z1\z5^7
+\a0\z1\z2\z3\z4\z5+\a6\z2^2\z3^2\z4^2\z5^2
\cr}
$$
$$
\ls1=(-2,1,0,0,0,0,1) \;\;,\;\;  \ls2=(-1,0,1,1,1,1,-3)
$$
$$
\eqalign{
\cD1&=\tx(3\ty-\tx)+x(\ty+2\tx+2)(\ty+2\tx+1)  \cr
\cD2&=\ty^2(7\ty-2\tx) \cr
&
+ 7y(3\ty-\tx+2)(3\ty-\tx+1)(3\ty-\tx)
+4x\ty^2(\ty+2\tx+1) \cr
}
$$
$$
dis_0=(1 - 4y)^4 + (27 - 441y + 2744y^2 - 38416y^3)x - 823543x^2y^3
$$

\vskip0.3cm

\model{\X{14}{\III}(7,3,2,1,1){2}{-260}:}
$$
\eqalign{ &
W=\z1^2+\z2^4\z3+\z3^7+\z4^{14}+\z5^{14}
\cr&
\hat W=\a1\z1^2+\a2\z2^4+\a3\z2\z7^7+\a4\z4^{14}+\a5\z5^{14}
+\a0\z1\z2\z3\z4\z5+\a6\z3^4\z4^4\z5^4
\cr}
$$
$$
\ls1=(-4,2,1,0,0,0,1)  \;\;,\;\; \ls2=(-2,1,0,2,1,1,-3)
$$
$$
\eqalign{
\eqalign{
\cD1&= \tx(\tx-3\ty)-4x(2\ty+4\tx+3)(2\ty+4\tx+1)\cr
\cD2&= (1-64x)^2\ty^3 \cr
&
-64\{112 x^2 y(\tx-3\ty)(2\ty+4\tx+1)+
xy (\tx-3\ty-1)(\tx-3\ty) \}   \cr
&
-(1-64x)\{112xy(\tx-3\ty-1)(\tx-3\ty)(2\ty+4\tx+1)   \cr
& \qquad \qquad \qquad
+y(\tx-3\ty-2)(\tx-3\ty-1)(\tx-3\ty) \} \cr}
}
$$
$$
\bx= 16 x \;\;,\;\; \by=y
$$
$$
dis_0=( 1 - 4\bx )^4+
( 27 - 441\bx + 2744\bx^2 - 38416\bx^3 ) \by - 823543\bx^3 \by^2
$$

\vskip0.3cm


\model{\X{15}{\III}(5,4,3,2,1){3}{-126}:}
$$
\eqalign{ &
W=\z1^3+\z2^3\z3+\z3^5+\z4^7\z5+\z5^{15}
\cr&
\hat W=\a1\z1^3+\a2\z2^3+\a3\z2\z3^5+\a4\z4^7+\a5\z4\z5^{15}
+\a0\z1\z2\z3\z4\z5+\a6\z1\z3^2\z5^7+\a7\z3^3\z4^3\z5^3
\cr}
$$
$$
\eqalign{
& \ls1=(-2,0,1,-1,0,-1,2,1) \;\;,\;\;  \ls2=(-1,1,0,1,0,1,-2,0) \;\; \cr
& \ls3=(0,-1,0,0,1,-1,3,-2)  \cr}
$$
$$
\eqalign{
\cD1&=\tx(\ty-\tz)(\tx-2\tz)-xy(2\tx+\ty+3)(2\tx+\ty+2)(2\tx+\ty+1) \cr
\cD2&=(\ty-\tz)(\ty-\tx)(\ty-\tx-\tz)   \cr
& -y(2\tx+\ty+1)(2\tx-2\ty+3\tz-1)(2\tx-2\ty+3\tz) \cr
\cD3&=(1-27xy)^2\tz(2\tx-2\ty+3\tz) \cr&
-(1-27xy)\{yz(\tx-2\tz-1)(\tx-2\tz)+45xy^2z(2\tx+\ty+1)(\tx-2\tz) \} \cr&
-3(1+54xy)xyz(\ty-\tx-\tz)(\tx-2\tz) - 405x^2y^2z(\ty-\tx)(\ty-\ty-\tz) \cr& 47

-27xy^2z(2\tx-2\ty+3\tz)(\tx-2\tz)-162x^2y^2z(2\tx+\ty+1)(\ty-\tx-\tz)  \cr&
-1215x^2y^3z(2\tx+\ty+1)(2\tx-2\ty+3\tz) \cr
\cD4&=
\tx(\tx-2\tz)-3xy(2\tx+\ty+2)(2\tx+\ty+1)-x(\tx-\ty)(\tx-\ty+\tz) \cr&
-3xy(2\tx+\ty+1)(2\tx-2\ty+3\tz) \cr
\cD5&=
(1-27xy)^2\{\ty^2+\tz^2-4\tx\ty-2\ty\tz+11\tz\tx \cr
& \hskip0.5cm
+9xy(2\tx+\ty+2)(2\tx+\ty+1)
-y(2\tx-2\ty+3\tz-1)(2\tx-2\ty+3\tz) \}  \cr
&
-21(1-27xy)xy^2z\{3(\tx-2\tz-1)(\tx-2\tz)+5(2\tx+\ty+1)(\tx-2\tz) \} \cr
&
-7xyz(1+54xy)(\ty-\tx-\tz)(\tx-2\tz)-945x^2y^2z(\ty-\tx)(\ty-\tx-\tz)
\cr
&
-63xy^2z(2\tx-2\ty+3\tz)(\tx-2\tz)+378x^2y^2z(2\tx+\ty+1)(\tx-\ty+\tz)\cr
&
-2835x^2y^3z(2\tx+\ty+1)(2\tx-2\ty+3\tz) \cr
}
$$
\vskip0.3cm

\model{\X{10}{\III}(3,2,2,2,1){3}{-144}:}
$$
\eqalign{ &
W=\z1^3\z4+\z5^5+\z3^5+\z4^5+\z5^{10}
\cr&
\hat W=\a1\z1^3+\a2\z2^5+\a3\z3^5+\a4\z4^5+\a5\z1\z5^{10}
+\a0\z1\z2\z3\z4\z5+\a6\z1^2\z5^5+\a7\z2\z3\z4\z5^6
\cr}
$$
$$
\eqalign{
&\ls1=(-1,0,0,0,0,-1,1,1) \;\;,\;\; \ls2=(-2,0,1,1,1,2,0,-3)  \;\; \cr
&\ls3=(0,1,0,0,0,1,-2,0) }
$$
$$
\eqalign{
\cD1&=
\(\tx-2\tz\)\(\tx-3\ty\) -3xz\(\tx-2\tz\)\(\tx-2\tz-1\) \cr
&
-2x\(\tx-2\tz+1\)\(2\ty-\tx+\tz\)
-x\(2\ty-\tx+\tz-1\)\(2\ty-\tx+\tz\) \cr
\cD2&=
\tz\(2\ty-\tx+\tz\)-z\(\tx-2\tz-1\)\(\tx-2\tz\)
\cr
\cD3&=
\tz\(\tx-3\ty\) -3xz\(\tz+1\)\(\tx-2\tz\) \cr
&
-2x\tz\(2\ty-\tx+\tz\)
-xz\(2\ty-\tx+\tz\)\(\tx-2\tz\) \cr
\cD4&=
\ty\(5\ty-2\tx\)\(2\tz-\tx\) -x\ty\(5\ty+2\tx+2\)\(2\ty-\tx+\tz\) \cr
&
-15x^3yz\(2\ty+\tx+4\)\(2\ty+\tx+2\)\(2\ty+\tx+1\)  \cr
&
-5xy\(2\ty+\tx+2\)\(3\ty-\tx+1\)\(3\ty-\tx\) \cr
&
+10x^2y\(2\ty+\tx+3\)\(2\ty+\tx+1\)\(3\ty-\tx\) \cr
\cD5&=
\ty^2\(5\ty-2\tx+\tz\)
+15x^2yz\(3\ty-\tx\)\(2\ty+\tx+2\)\(2\ty+\tx+1\) \cr
&
+5y\(3\ty-\tx+2\)\(3\ty-\tx+1\)\(3\ty-\tx\)  \cr
&
-10xy\(2\ty+\tx+1\)\(3\ty-\tx+1\)\(3\ty-\tx\) \cr
&
+3xz\ty^2\(\tx-2\tz\) + 2x \ty^2\(2\ty-\tx+\tz\) \cr
}
$$
\vskip0.3cm

\model{\X{10}{\III}(3,3,2,1,1){3}{-168}:}
$$
\eqalign{ &
W=\z1^3\z4+\z2^3\z5+\z3^5+\z4^{10}+\z5^{10}
\cr&
\hat W=\a1\z1^3+\a2\z2^3+\a3\z3^5+\a4\z1\z4^{10}+\a5\z2\z5^{10}
+\a0\z1\z2\z3\z4\z5+\a6\z3\z4^6\z5^6+\a7\z3^3\z4^3\z5^3
\cr}
$$
$$\eqalign{
& \ls1=(-2,1,1,0,-1,-1,2,0) \;\;,\;\; \ls2=(-1,0,0,0,1,1,-2,1,0)  \;\; \cr
& \ls3=(0,0,0,1,0,0,1,-2,0)  \cr}
$$
$$
\eqalign{
\cD1&=\tz(2\tx-2\ty+\tz)-z(\ty-2\tz-1)(\ty-2\tz) \cr
\cD2&=(4\tx-\ty)(\ty-2\tz)  \cr
& -9xy(2\tx+\ty+2)(2\tx+\ty+1)
  +y(2\tx-2\ty+\tz)(2\tx-2\ty+\tz-1) \cr
\cD3&=(\tx-\ty)^2+(\tx-\ty)(2\tx-2\ty+\tz)+3\tx(2\tx-2\ty+\tz) \cr
&
-5yz(2\tx-2\ty+\tz)(\ty-2\tz)-3y(2\tx-2\ty+\tz-1)(2\tx-2\ty+\tz) \cr
&
-45xyz(2\tx+\ty+1)(\ty-2\tz)
+27xy(2\tx+\ty+1)(2\tx-2\ty+\tz)\cr
&
-9x(\tx-\ty)^2  \cr
\cD4&=\tz(\tx-\ty)^2-yz(2\tx+\ty+1)(2\tx-2\ty+\tz)(\ty-2\tz) \cr
\cD5&=(\tx-\ty)^2(\ty-2\tz)-y(2\tx+\ty+1)(2\tx-2\ty+\tz)(2\tx-2\ty+\tz-1)\cr
}
$$
\vskip0.3cm

\model{\X{20}{\III}(10,4,3,2,1){3}{-192}:}
$$
\eqalign{ &
W=\z1^2+\z2^5+\z3^6\z4+\z4^{10}+\z5^{20}
\cr&
\hat W=\a1\z1^2+\a2\z2^5+\a3\z3^6+\a4\z3\z4^{10}+\a5\z5^{20}
+\a0\z1\z2\z3\z4\z5+\a6\z2\z4^6\z5^6+\a7\z3^3\z5^{10}
\cr}
$$
$$
\eqalign{
& \ls1=(-2,1,0,0,-1,-1,2,1) \;\;,\;\;  \ls2=(-2,1,1,0,2,1,-3,0) \;\;  \cr
& \ls3=(0,0,0,1,0,1,0,-2)   \cr}
$$
$$
\eqalign{
\cD1&=\tz(\ty+\tz-\tx)-z(\tx-2\tz-1)(\tx-2\tz) \cr
\cD2&=\tz(2\tx-3\ty)  \cr
&
+2xz(\tx-2\ty)(\tx-2\tz)-40x^2yz(2\tx+2\ty+3)(2\tx+2\ty+1)\cr
\cD3&=(2\tx-3\ty)(\tx-2\tz)-2x(\tx-2\ty)(\tx-\ty-\tz)     \cr
& -20xy(2\tx-3\ty-1)(2\tx-3\ty)+120x^2y(\tx-\ty-\tz)(2\tx+2\ty+1) \cr
& -240x^2yz(\tx-2\tz)(2\tx+2\ty+1) \cr
\cD4&=\ty(\tx-2\ty)(\tx-2\tz)+4xy(2\tx-3\ty)(2\tx+2\ty+3)(2\tx+2\ty+1)\cr
\cD5&=\ty(\tx-2\ty)(\tx-\ty-\tz)-2y(2\tx-3\ty-2)(2\tx-3\ty-1)(2\tx-3\ty)\cr
&
+12xy(2\tx-3\ty)(\tx-\ty-\tz)(2\tx+2\ty+1) \cr
&
-24xyz(2\tx-3\ty)(\tx-2\tz)(2\tx+2\ty+1) \cr
 }
$$
\vskip0.3cm

\model{\X{10}{\I}(4,2,2,1,1){3}{-192}:}
$$
\eqalign{ &
W=\z1^2\z3+\z2^5+\z3^5+\z4^{10}+\z5^{10}
\cr&
\hat W=\a1\z1^2+\a5\z2^5+\a3\z1\z3^5+\a4\z4^{10}+\a5\z5^{10}+
\a0\z1\z2\z3\z4\z5+\a6\z2^2\z3^2\z4^2\z5^2+\a7\z4^5\z5^5
\cr}
$$
$$
\eqalign{
& \ls1=(-2,1,0,0,0,0,1,0) \;\;,\;\;  \ls2=(-1,0,1,1,0,0,-2,1) \;\; \cr
& \ls3=(0,0,0,0,1,1,0,-2)  \cr}
$$
$$
\eqalign{
\cD1&=\tx(\tx-2\ty)-x(2\tx+\ty+2)(2\tx+\ty+1)  \cr
\cD2&=\tz^2-z(\ty-2\tz-1)(\ty-2\tz)  \cr
\cD3&=4\tx^2+4\tx\tz-10\tx\ty-10\tz\ty+5\ty^2 \cr
&
-8x(\tx+\tz+1)(2\tx+\ty+1) -5y(\tx-2\ty-1)(\tx-2\ty) \cr
 }
$$
$$
\bx=x \;\;,\;\; \by=y \;\;,\;\; \bz=4 z
$$
$$
\eqalign{
& dis_0=
   ((1 - 4\bx)^3 - 2(2 - 25\bx + 500\bx^2)\by  + 3125\bx^2\by^2)^2  \cr
& - 2\by^2(8 - 200\bx + 2125\bx^2 - 12500\bx^3 + 350000\bx^4 + 200000\bx^5 -
     12500\bx^2\by  \cr
& + 156250\bx^3\by - 3125000\bx^4\by + 9765625\bx^4\by^2)\bz \
   + 9765625\bx^4\by^4\bz^2  \cr
&
dis_1=1-\bz  \cr}
$$

\vskip0.3cm

\model{\X{16}{\III}(8,3,2,2,1){3}{-200}:}
$$
\eqalign{ &
W=\z1^2+\z2^5\z5+\z3^8+\z4^8+\z5^{16}  \cr
\cr&
\hat W=\a1\z1^2+\a2\z2^5+\a3\z3^8+\a4\z4^8+\a5\z2\z5^{16}
+\a0\z1\z2\z3\z4\z5+\a6\z3^2\z4^2\z5^{10}+\a7\z2^3\z5^8
\cr}
$$
$$
\eqalign{
& \ls1=(-2,1,0,0,0,-1,1,1) \;\;,\;\;  \ls2=(-2,1,0,1,1,2,-3,0) \;\; \cr
& \ls3=(0,0,1,0,0,1,0,-2)   \cr}
$$
$$
\eqalign{
\cD1&= \tz(\tx-2\ty-\tz)+z(\tx-2\tz-1)(\tx-2\tz) \cr
\cD2&=\tz(\tx-3\ty)-2xz(2\tx+2\ty+1)(\tx-2\tz) \cr
\cD3&=(\tx-3\ty)(\tx-2\tz)+2x(2\tx+2\ty+1)(\tx-2\ty-\tz)  \cr
\cD4&=\ty^2(2\ty-\tx+\tz)-2y(\tx-3\ty-2)\{20x^2z(2\tx+2\ty+3)(2\tx+2\ty+1) \cr
& \hskip2cm
+(\tx-3\ty-1)(\tx-3\ty)+6x(\tx-3\ty)(2\tx+2\ty+1) \}  \cr
\cD5&=\ty^2(\tx-2\tz)-4xy(2\tx+2\ty+3)\{20x^2z(2\tx+2\ty+3)(3\tx+2\ty+1) \cr
& \hskip2cm
+6x(2\tx+2\ty+1)(\tx-3\ty)+(\tx-3\ty-1)(\tx-3\ty) \}
 }
$$
\vskip0.3cm

\model{\X{12}{\III}(5,3,2,1,1){3}{-204}:}
$$
\eqalign{ &
W=\z1^2\z3+\z2^4+\z3^6+\z4^{12}+\z5^{12}
\cr&
\hat W=\a1\z1^2+\a2\z2^4+\a3\z1\z3^6+\a4\z4^{12}+\a5\z5^{12}
+\a0\z1\z2\z3\z4\z5+\a6\z3^4\z4^4\z5^4+\a7\z2^2\z3^2\z4^2\z5^2
\cr}
$$
$$
\eqalign{
& \ls1=(0,1,0,-2,-1,-1,3,0) \;\;,\;\;  \ls2=(0,0,1,0,0,0,1,-2,0) \;\; \cr
& \ls3=(-2,0,0,2,1,1,-3,1)   \cr}
$$
$$
\eqalign{
\cD1&= \tx(\tz-2\ty)-xz(2\tz+2)(2\tz+1)  \cr
\cD2&= \ty(3\tx+\ty-3\tz)-y(\tz-2\ty-1)(\tz-2\ty)  \cr
\cD3&=\tx(3\tx+\ty-3\tz)-2xz(2\tz+1)(3\tx+\ty-3\tz)+4xyz(2\tz+1)(2\ty-\tz) \cr
\cD4&=\{64x^2yz^2-(1-4xz)^2\}^2\{(3\tx\tz-3\tx^2-\tz^2)(\tz-2\ty)
                                  +2xz(\tx+1)^2(2\tz+1)\} \cr
&
+64x^2yz^2\{(1-4xz)z(3\tx+\ty-3\tz-1)(3\tx+\ty-3\tz) \cr
& \hskip1cm
+8xyz^2(3\tx+\ty-3\tz)(\tz-2\ty)
+96x^2yz^3(3\tx+\ty-3\tz)(2\tz+1) \} \cr
&
-\{64x^2yz^2-(1-4xz)^2\}  \cr
& \hskip1.5cm \times
 \{(1-4xz)z(3\tx+\ty-3\tz-2)(3\tx+\ty-3\tz-1)(3\tx+\ty-3\tz) \cr
&  \hskip2.5cm
+8xyz^2(3\tx+\ty-3\tz-1)(3\tx+\ty-3\tz)(\tz-2\ty)    \cr
&  \hskip2.5cm
+96x^2yz^3(3\tx+\ty-3\tz-1)(3\tx+\ty-3\tz)(2\tz+1)  \} \cr
\cD5&=\{64x^2yz^2-(1-4xz)^2\}(\tx-\ty)^3 \cr
&
+64x^2yz^2\{(1-4xz+16xyz)z(3\tx+\ty-3\tz-1)(3\tx+\ty-3\tz) \cr
& \hskip1.5cm
+2yz(3\tx+\ty-3\tz)(\tz-2\ty)+24xyz^2(2\tz+1)(3\tx+\ty-3\tz) \}\cr
&
-\{64x^2yz^2-(1-4xz)^2\} \cr
& \hskip1.5cm \times
\{(1-4xz)z(3\tx+\ty-3\tz-2)(3\tx+\ty-3\tz-1)(3\tx+\ty-3\tz) \cr
& \hskip2.5cm
+16xyz^2(3\tx+\ty-3\tz-1)^2(3\tx+\ty-3\tz)     \cr
& \hskip2.5cm
+2yz(3\tx+\ty-3\tz-1)(3\tx+\ty-3\tz)(\tz-2\ty) \cr
& \hskip2.5cm
+24xyz^2(3\tx+\ty-3\tz-1)(2\tz+1)(3\tx+\ty-3\tz) \}
 }
$$

\vskip0.3cm

\model{\X{18}{\II}(9,4,2,2,1){3}{-240}:}
$$
\eqalign{ &
W=\z1^2+\z2^4\z4+\z3^9+\z4^9+\z5^{18}
\cr&
\hat W=\a1\z1^2+\a2\z2^4+\a3\z3^9+\a4\z2\z4^9+\a5\z5^{18}
+\a0\z1\z2\z3\z4\z5+\a6\z3^4\z4^4\z5^4+\a7\z1\z5^9
\cr}
$$
$$
\eqalign{
& \ls1=(-4,0,1,0,0,-2,1,4) \;\;,\;\;  \ls2=(0,1,0,0,0,1,0,-2) \;\; \cr
& \ls3=(-1,0,0,1,1,0,-2,1)  \cr}
$$
$$
\eqalign{
\cD1&= \ty(\ty-2\tx)-y(4\tx-2\ty+\tz-1)(4\tx-2\ty+\tz) \cr
\cD2&= \tx(\tx-2\tz)
-x\{2(1-8y)\tx-(1-4y)\ty-4y\tz+1-10y\} \cr
& \hskip2.5cm \times
  \{2(1-8y)\tx-(1-4y)\ty-4y\tz-2y \} \cr
\cD3&=(4\tx-9\tz)(4\tx-2\ty+\tz)-9z(\tx-2\tz-1)(\tx-2\tz) \cr
&
+16x\{2(4y-1)\tx+(1-4y)\ty+2y\tz+4y-1\}\cr
& \hskip1cm \times
    \{2(1-8y)\tx+(4y-1)\ty-4y\tz-2y\} \cr
\cD4&=\tz^2(\ty-2\tx)+z(\tx-2\tz-1)(\tx-2\tz)\{2(1-4y)\tx-(1-4y)\ty-2y\tz\} \cr
\cD5&=\tz^2(4\tx-2\ty+\tz)-z(4\tx+\tz+1)(\tx-2\tz-1)(\tx-2\tz) \cr
 }
$$
\vskip0.3cm

\model{\X{9}{\I}(4,2,1,1,1){3}{-240}:}
$$
\eqalign{ &
W=\z1^2\z5+\z2^4\z4+\z3^9+\z4^9+\z5^9
\cr&
\hat W=\a1\z1^2+\a2\z2^4+\a3\z3^9+\a4\z2\z4^9+\a5\z1\z5^9
+\a0\z1\z2\z3\z4\z5+\a6\z3^4\z4^4\z5^4+\a7\z2^2\z3^2\z4^2\z5^2
\cr}
$$
$$
\eqalign{
& \ls1=(-2,1,0,0,0,0,0,1) \;\;,\;\;  \ls2=(-1,0,0,1,1,1,-2,0) \;\; \cr
& \ls3=(0,0,1,0,0,0,1,-2)   \cr}
$$
$$
\eqalign{
\cD1&= \tx(\tx-2\tz)-x(2\tx+\ty+2)(2\tx+\ty+1) \cr
\cD2&= \tz(\tz-2\ty)-z(\tx-2\tz-1)(\tx-2\tz)  \cr
\cD3&= 9\ty^2-18\tx\ty-64\tx\tz+36\tx^3-9y(2\ty-\tz+1)(2\ty-\tz) \cr
& -x(72\tx+16\tz+72)(2\tx+\ty+1)-32xz(2\tx+\ty+1)(\tx-2\tz) \cr
 }
$$

\vskip0.3cm

\model{\X{16}{\II}(8,3,3,1,1){3}{-256}:}
$$
\eqalign{ &
W=\z1^2+\z2^5\z4+\z3^5\z5+\z4^{16}+\z5^{16}
\cr&
\hat W=\a1\z1^2+\a2\z2^5+\a3\z3^5+\a4\z2\z4^{16}+\a5\z3\z5^{16}
+\a0\z1\z2\z3\z4\z5+\a6\z4^{10}\z5^{10}+\a7\z1\z4^5\z5^5
\cr}
$$
$$
\eqalign{
& \ls1=(-5,0,1,1,0,0,-2,5) \;\;,\;\;  \ls2=(0,1,0,0,0,0,1,-2) \;\; \cr
& \ls3=(-1,0,0,0,1,1,-2,1)    \cr}
$$
$$
\eqalign{
\cD1&= \ty(\ty-2\tx-2\tz)-y(5\tx-2\ty+\tz-1)(5\tx-2\ty+\tz)  \cr
\cD2&=\tz^2-2yz(\ty-2\tx-2\tz)(5\tx-2\ty+\tz)  \cr
&  -z(\ty-2\tx-2\tz-1)(\ty-2\tx-2\tz) \cr
\cD3&=(5\tx-\tz)(5\tx-2\ty+\tz) +z(\ty-2\tx-2\tz-1)(\ty-2\tx-2\tz) \cr
&     -25x\{(2-20y)\tx+(4y-1)\ty+(2-4y)\tz+1-10y\} \cr
&\hskip0.5cm \times
          \{(2-20y)\tx+(4y-1)\ty+(2-4y)\tz \}       \cr
\cD4&=\tx^2(\ty-2\tx-2\tz)  \cr
&-x\{(2-20y)\tx+(4y-1)\ty+(2-4y)\tz+2-14y\} \cr
& \hskip0.5cm
   \times
   \{(2-20y)\tx+(4y-1)\ty+(2-4y)\tz+1-6y\}  \cr
& \hskip0.5cm
   \times
   \{(10y-2)\tx+(1-4y)\ty+(2y-2)\tz \}  \cr
\cD5&=\tz^2(5\tx-2\ty+\tz)-z(5\tx+\tz+1)(\ty-2\tx-2\tz-1)(\ty-2\tx-2\tz) \cr
 }
$$

\vskip0.3cm

\model{\X{16}{\II}(8,5,1,1,1){3}{-456}:}
$$
\eqalign{ &
W=\z1^2+\z2^3\z5+\z3^{16}+\z4^{16}+\z5^{16}
\cr&
\hat W=\a1\z1^2+\a2\z2^3+\a3\z3^{16}+\a4\z4^{16}+\a5\z2\z5^{16}
+\a0\z1\z2\z3\z4\z5+\a6\z3^6\z4^6\z5^6+\a7\z1\z3^3\z4^3\z5^3
\cr}
$$
$$
\eqalign{
& \ls1=(-3,0,1,0,0,0,-1,3) \;\;,\;\;  \ls2=(-1,0,0,1,1,1,-3,1) \;\; \cr
& \ls3=(0,1,0,0,0,0,1,-2)  \cr}
$$
$$
\eqalign{
\cD1&= \tz(\tz-3\ty-\tx)-z(3\tx+\ty-2\tz-1)(3\tx+\ty-2\tz) \cr
\cD2&= \tx(3\tx+\ty-2\tz)
-2xz(3\tx+\ty+3)(6\tx+2\ty-2\tz+1) \cr
& \hskip3cm
  + x(3\tx+\ty+3)(\tx+3\ty-\tz) \cr
\cD3&= \tx(\tz-\tx-3\ty)
   -x\{(6z-1)\tx+(2z-3)\ty-(4z-1)\tz+4z-1\}  \cr
& \hskip3cm
     \{(12z-1)\tx+(4z-3)\ty-(4z-1)\tz+2z\}\cr
\cD4&=\ty^3-2yz(\tx+3\ty-\tz+1)(\tx+3\ty-\tz)(3\tx+\ty-2\tz) \cr
& \hskip1cm +y(\tx+3\ty-\tz+2)(\tx+3\ty-\tz+1)(\tx+3\ty-\tz) \cr
\cD5&=\ty^2(3\tx+\ty-2\tz)
-6xz\ty^2(6\tx+2\ty-2\tz+1) +3x\ty^2(\tx+3\ty-\tz) \cr
&
+y(\tx+3\ty-\tz+2)(\tx+3\ty-\tz+1)(\tx+3\ty-\tz) \cr
 }
$$

\bye